\documentclass[11pt,leqno]{article}
\usepackage{amsmath,amssymb}
\usepackage{natbib}

\renewcommand{\thesection}{\arabic{section}}

\newcommand{\email}[1]{{\it  E-mail address:}\  \textsf{#1}\\*\protect}
\newcommand{\address}[1]{\small \textsc{#1}\\*\protect}

\DeclareMathOperator{\sgn}{sgn}
\title{Pathwise Uniqueness for Stochastic Heat Equations with H\"older Continuous Coefficients: the White Noise Case}
\bigskip
\author{Leonid Mytnik $\mbox{}^{1}$ \hspace{1cm} 
 Edwin Perkins $\mbox{}^{2}$  }
\date{}
\begin{document}
\maketitle
\begin{center}
\address{Faculty of Industrial Engineering and Management,}
\address{Technion -- Israel Institute of Technology, Haifa 32000, Israel}
\email{leonid@ie.technion.ac.il}
\mbox{}\\
\address{ Department of Mathematics, The University of British Columbia,}
\address{ 1984 Mathematics Road, Vancouver, B.C., Canada V6T 1Z2}
\email{perkins@math.ubc.ca}
\end{center}

\vskip 0.3in 
{\bf Abstract.} We prove pathwise uniqueness for solutions of parabolic stochastic pde's with multiplicative white noise if the coefficient is H\"older continuous of index $\gamma>3/4$.  The method of proof is an infinite-dimensional version of the Yamada-Watanabe argument for ordinary stochastic differential equations. 
\vfill
{\footnoterule
\footnotesize
\noindent 
\today\\
AMS 2000 {\it subject classifications}. Primary 60H15. 
Secondary 60G60, 60H10, 60H40, 60K35, 60J80. 
  \\
{\it Running head}. Pathwise uniqueness for SPDE's\\
1. Supported in part by the Israel Science Foundation 
(grant No. 1162/06). \\ 
2. Supported by an NSERC Research grant.\\}
\pagebreak
\vspace{1cm}

\newcommand{\mysection}[1]{\section{#1}
\markright{{\protect\footnotesize Chapter \thechapter, Section \thesection}}}
\newcommand{\resc}[2]{{#1}^{n}(\qn)={#2}
\lbr{\displaystyle \frac{\qn}{n}} \rbr}
\newcommand{\rescb}[2]{\bar{#1}^{n}(\qn)={#2}
\lbr{\displaystyle \frac{\qn}{n}} \rbr}
\newcounter{archapter}[section]
\newtheorem{theorem}{Theorem}[section]
\newtheorem{lemma}[theorem]{Lemma}
\newtheorem{definition}[theorem]{Definition}
\newtheorem{remark}[theorem]{Remark}
\newtheorem{proposition}[theorem]{Proposition}
\newtheorem{case}{Case}[archapter]
\newtheorem{corollary}[theorem]{Corollary}
\newtheorem{fact}[theorem]{Fact}
\newtheorem{assumption}[theorem]{Assumption}
\newtheorem{convention}[theorem]{Convention}
\newtheorem{conjecture}[theorem]{Conjecture}
\newcommand{\gdm}{\hfill\vrule  height5pt width5pt \vspace{.1in}}
\newcommand{\itemstr}[1]{\item[{\bf (#1)}]}
\newfont{\msb}{msbm10 scaled \magstep1}
\newfont{\msbh}{msbm7 scaled \magstep1}
\newfont{\msbhh}{msbm5 scaled \magstep1}
\renewcommand{\Re}{\mbox{\msb R}}
\newcommand{\IR}{\mbox{\msb R}}
\newcommand{\IRsm}{\mbox{\msbh R}}
\newcommand{\Z}{\mbox{\msb Z}}
\newcommand{\Ze}{\mbox{\msb Z}}
\newcommand{\N}{\mbox{\msb N}}
\newcommand{\Zesm}{\mbox{\msbh Z}}
\newcommand{\Ce}{\mbox{\msb C}}
\newcommand{\Cesm}{\mbox{\msbh C}}
\newcommand{\Te}{\mbox{\msb T}}
\newcommand{\Tesm}{\mbox{\msbh T}}
\newcommand{\Le}{\mbox{\msb L}}
\newcommand{\Lesm}{\mbox{\msbh L}}
\newcommand{\Ee}{\mbox{\msb E}}
\newcommand{\Pe}{\mbox{\msb P}}
\newcommand{\Resm}{\mbox{\msbh R}} 
\newcommand{\Resmm}{\mbox{\msbhh R}} 
\newcommand{\Rp}{\Re_{+}}
\newcommand{\Rpf}{\Re_{+}}
\newcommand{\Rpsm}{\Resm_{+}}
\newcommand{\Rd}{\Re^{d}}
\newcommand{\Rdsm}{\Resm^{d}} 
\newcommand{\Rdsmm}{\Resmm^{d}} 
\newcommand{\Rdn}[1]{\Re^{#1}}
\newcommand{\Rdnsm}[1]{\Resm^{#1}}
\newcommand{\Zcal}{{\cal Z}}
\newcommand{\esssup}{\mathrm{ess\, sup}} 
\newcommand{\lsc}{\left\langle}
\newcommand{\rsc}{\right\rangle}
\newcommand{\lbr}{\left(}
\newcommand{\rbr}{\right)}
\newcommand{\lfi}{\left\{}
\newcommand{\rfi}{\right\}}
\newcommand{\lsbr}{\left[}
\newcommand{\rsbr}{\right]}
\newcommand{\rmed}{\right|}
\newcommand{\lmed}{\left|}
\newcommand{\lnor}{\left\|}
\newcommand{\rnor}{\right\|}
\newcommand{\lnorm}{\lmed \lmed}
\newcommand{\rnorm}{\rmed \rmed}
\newcommand{\cndp}{\stackrel{\rm P}{\rightarrow}}
\newcommand{\mf}{M_{F}(\Rd)}
\newcommand{\mfo}{{{\cal M}_{F}}}
\newcommand{\scf}[2]{\lsc #1,\,#2 \rsc} 
\newcommand{\Nat}{{\cal N}} 
\newcommand{\rmp}{{\mathrm{p}}}
\newcommand{\ep}{\epsilon}
\newcommand{\ve}{\varepsilon}
\newcommand{\la}{\langle}
\newcommand{\ra}{\rangle}
\newcommand{\cF}{{\cal F}}
\newcommand{\cB}{{\cal B}}
\newcommand{\cI}{{\cal I}}
\newcommand{\tcI}{\tilde{\cal I}}
\newcommand{\nn}{\nonumber}
\newcommand{\om}{\overline m}
\newcommand{\oln}{\overline{l_n}}
\newcommand{\Ft}[1]{\mbox{${\cal F}_{#1}$}}
\newcommand{\mn}{\medskip\noindent}
\newcommand{\vsp}{\vspace{4.5 mm}} 
\newcommand{\mvsp}[1]{\vspace{#1 mm}} 
\newcommand{\hsp}{\mbox{\hspace{10 mm}}} 
\newcommand{\mhsp}[1]{\mbox{\hspace{#1 mm}}}
\newcommand{\intobl}[1]{
 \int\!\!\int_{{}_{{}_{\!\!\!\!\!\!\!\!\!\mbox{\scriptsize $#1$}}}}\;}
\section{Introduction}
\setcounter{equation}{0}
\setcounter{theorem}{0}

Let $\sigma:\Rp\times\Re^2\to \Re$ and consider the stochastic heat equation
\begin{equation}
\label{spde}
\frac{\partial}{\partial t}
X(t,x) = \frac{1}{2} \Delta X(t,x) dt + \sigma(t,x,X(t,x)) \dot{W}(x,t)+b(t,x,X(t,x)).
\end{equation}
Here $\Delta$ denotes the Laplacian and $\dot{W}$ is space-time white
noise on $\Rp\times \Re$. 
If $\sigma(t,x,X)$ and $b(t,x,X)$ are Lipschitz continuous in $X$ it is well-known that there are pathwise unique solutions to \eqref{spde} (see \cite{bib:wal86}).  When $\sigma(t,x,X)=\sqrt{f(t,x,X)X}$ such equations arise naturally as the scaling limits of critical branching particle systems where the branching rate at $(t,x)$ is given by $f(t,x,X(t,x))$ and $X(t,x)$ is a measure of the local particle density at $(t,x)$.  Such coefficients are not Lipschitz continuous and pathwise uniqueness remains open even in the case where $f\equiv 1$, $b\equiv 0$ and $X$ is the density of super-Brownian motion (see Section III.4 of \cite{bib:per02}).  In this case, and more generally for $f=X^p$ for $p>0$, uniqueness in law is known by duality arguments (see \cite{bib:myt99}).  The duality arguments are highly non-robust, however, and pathwise uniqueness, if true, would typically hold for a much less restrictive set of coefficients.  Our goal in this work is to show pathwise uniqueness holds for solutions to \eqref{spde} if $\sigma(t,x,\cdot)$ is H\"older continuous of order $\gamma$ and $\gamma>3/4$.  The attentive reader will have already noted that the motivating example given above does not satisfy this condition.

The above equation does have the advantage 
of having a diagonal form--that is, when viewed as a continuum-dimensional
stochastic differential equation there are no off-diagonal terms in the
noise part of the equation and the diffusion coefficient for the $x$ coordinate is a
function of that coordinate alone.  For finite-dimensional sde's this was the setting for
Yamada and Watanabe's extension \cite{bib:YW71} of It\^o's pathwise uniqueness results to
H\"older ($1/2$) continuous coefficients, and so our plan will be to carry over their approach
to our infinite dimensional setting.  This programme was already carried out in the context of coloured noise in \cite{bib:mps}, but the methods used there when specialized to white noise given nothing beyond the classical Lipschitz uniqueness. In fact for coloured noise in higher dimensions the results in \cite{bib:mps} did not even come close to the known results on pathwise uniqueness for Lipschitz continuous coefficients \cite{bib:dal99}--see the discussion after Remark~1.5 in \cite{bib:mps}.  This is what led to our belief that there was room for substantial improvement in the methods of \cite{bib:mps} and hence to the present work. 

We introduce a growth condition, a H\"older continuity condition on $\sigma$ and the standard Lipschitz condition on $b$:
\begin{align}
\label{growthcond}
\hbox{ there exists a }&\hbox{constant }c_{\ref{growthcond}}\hbox{ such that for all }(t,x,X) \in \Rp\times\Re^2,\\
&\nonumber|\sigma(t,x,X)| +|b(t,x,X)|\leq c_{\ref{growthcond}} (1 + |X|),
\end{align}
\begin{align}
\label{holder}
&\hbox{for some $\gamma>3/4$ there are $R_1,R_2>0$ and for all $T>0$ there is an $R_0(T)$}\\
 \nonumber&\hbox{so that for all $t\in [0,T]$ and all $(x,X,X')\in\Re^3$,}\\ 
\nonumber &|\sigma(t,x,X)-\sigma(t,x,X')|\le R_0(T)e^{R_1|x|}(1+|X|+|X'|)^{R_2}|X-X'|^\gamma,
\end{align}
and
\begin{align}\label{bLip}
\hbox{there is a }B>0 \hbox{ s.t. for all }(t,x,X,X')\in \Rp\times\Re^3,\ |b(t,x,X)-b(t,x,X')|\le B|X-X'|.
\end{align}

We assume $W$ is a white noise on the filtered probability space $(\Omega,\cF,\cF_t,\Pe)$, where $\cF_t$ satisfies the usual hypotheses.  This means $W_t(\phi)$ is an $\cF_t$-Brownian motion with 
variance $\Vert\phi\Vert_2^2$ for each $\phi\in \Le^2(\Re,dx)$
 and $W_t(\phi)$ and $W_t(\psi)$ are independent if $\int\phi(x)\psi(x)dx=0$.  We set $p_t(x)=(2\pi 
t)^{-1/2}\exp\{-x^2/2t\}$, $P_tf(x)=\int f(y)p_t(y-x)dy$, and let $\cF^W_t\subset\cF_t$ be the filtration generated by $W$ satisfying the usual hypotheses.  
A stochastic process $X: \Omega \times \Rp\times \Re
\rightarrow \Re$, which is jointly measurable and ${\cal F}_{t}$-adapted,
is said to be a solution to the stochastic heat equation (\ref{spde})
on $(\Omega,\cF,\cF_t,\Pe)$ with initial condition $X_{0}:\Re\to\Re$, if for each $t\ge 0$, and $x\in\Re$, \begin{eqnarray}
\label{wspde}
X(t,x) &=& \int_{\Re} p_{t}(y-x) X_{0}(y) dy
+ \int_{0}^{t}  \int_{\Re} p_{t-s}(y-x) \sigma(s,y,X(s,y)) W(ds, dy)\\
\nn&&\phantom{\int_{\Re} p_{t}(y-x) X_{0}(y) dy
}+\int_0^t\int_{\Re}p_{t-s}(y-x)b(s,y,X(s,y))dyds\ a.s.
\end{eqnarray}
\noindent

To state the main results we introduce some notation, which will
be used throughout this work:  If $E\subset \Re^d$, we write $C(E)$ for the space of continuous
functions on $E$. A superscript $k,$ respectively $\infty$, indicates that
functions are in addition $k$ times, respectively infinitely often, continuously
differentiable. A subscript $b,$ respectively $c,$ indicates that they are also
bounded, respectively have compact support. We also define
\begin{equation*}
||f||_{\lambda}:=\sup_{x \in \Re} |f(x)| e^{-\lambda |x|},
\end{equation*}
set $C_{tem}:=\{f \in C(\Re) , ||f||_{\lambda}
< \infty  \text{ for any } \lambda >0\}$ and endow it with the topology induced
by the norms $||\cdot ||_{\lambda}$ for $\lambda>0.$ That is,
$f_n\to f$ in $C_{tem}$ iff $d(f,f_n)=\sum_{k=1}^\infty2^{-k}(\Vert f-f_n\Vert_{1/k}\wedge 1)\to 0$ as $n\to\infty$.
Then $(C_{tem},d)$ is a Polish space.  By identifying the white noise $W$, with the associated Brownian sheet, we may view $W$ as a stochastic processes with sample paths in $C(\Rp,C_{tem})$.  Here as usual, $C(\Rp,C_{tem})$ is given the topology of uniform convergence on compacts.

A stochastically weak solution to
(\ref{spde}) is a solution on some filtered space with respect to some
noise $W$, i.e., the noise and space are not specified in advance.

With this notation we can state the following standard existence
result whose proof is a minor modification of Theorem~1.2 of \cite{bib:mps} and is given in the next Section.
\begin{theorem}
\label{theorem:nonLipschitzexistence}
Let $X_{0} \in C_{tem},$ and  let $b, \sigma:\Rp\times\Re^2\to\Re$ satisfy (\ref{growthcond}), \eqref{holder}, and \eqref{bLip}.  Then there exists a stochastically
weak solution to (\ref{spde}) with sample paths a.s. in
$C(\Rp, C_{tem})$.
\end{theorem}

We say pathwise uniqueness holds for solutions of (\ref{spde}) in
$C(\Rp,C_{tem})$ if for every $X_0\in C_{tem}$, any two solutions to
(\ref{spde}) with sample paths a.s. in $C(\Rp,C_{tem})$ must be equal
with probability
$1$.  For Lipschitz continuous $\sigma$, this follows from Theorem~2.2 of \cite{bib:shiga94}.  
Here then is our main result:

\begin{theorem}
\label{theorem:unique}
Assume that $b,\sigma:\Rp\times\Re^2\to\Re$ satisfy 
\eqref{growthcond}, \eqref{holder} and \eqref{bLip}. 
Then pathwise uniqueness holds for 
solutions of \eqref{spde} in $C(\Rp,C_{tem})$.
\end{theorem}

As an immediate consequence of Theorems~\ref{theorem:nonLipschitzexistence} and \ref{theorem:unique} we get existence and uniqueness of strong solutions and joint uniqueness in law of $(X,W)$.

\begin{theorem}
\label{strong}
Assume that $b,\sigma:\Rp\times\Re^2\to\Re$ satisfy 
\eqref{growthcond}, \eqref{holder} and \eqref{bLip}.  Then for any $X_0\in C_{tem}$ there is a solution $X$ to
\eqref{spde} on $(\Omega,\cF^W_\infty,\cF^W_t,\Pe)$ with sample paths a.s. in $C(\Rp, C_{tem})$. If $X'$ is any other solution to \eqref{spde} on $(\Omega,\cF,\cF_t,\Pe)$ with sample paths a.s. in $C(\Rp, C_{tem})$, then $X(t,x)=X'(t,x)$ for all $t,x$ a.s. The joint law $P_{X_0}$ of $(X,W)$ on $C(\Rp, C_{tem})$ is uniquely determined by $X_0$ and is Borel measurable in $X_0$.  
\end{theorem}
\paragraph{Proof.} The Borel measurability of the law is proved as in Exercise 6.7.4 in \cite{bib:SV79}.  
We now apply Theorem 3.14 of \cite{bib:kur07}, with the Polish state spaces $S_1$ and $S_2$ for the driving process ($W$) and solution ($X$) in that work both equal to $C(\Rp,C_{tem})$.  Theorems~\ref{theorem:nonLipschitzexistence} and \ref{theorem:unique} imply the hypotheses of weak existence and pointwise uniqueness of (a) of that result.  The conclusions of an $\cF^W_t$-adapted (strong) solution and uniqueness in law of $(X,W)$ follow from the conclusions in Theorem~3.14 (b) (of \cite{bib:kur07}) of a strong compatible solution and joint uniqueness in law, respectively.  (Note that Lemma~3.11 of the above reference shows that a strong, compatible solution must be $\cF^W_t$-adapted.)
\gdm

\begin{remark}\label{holderhyp} (a) When assuming \eqref{growthcond}, it suffices to assume \eqref{holder} for $|X-X'|\le 1$.  Indeed this condition is immediate from \eqref{growthcond} for $|X-X'|\ge 1$ with 
$R_1=0$ and $R_2=1$. 

\medskip

\noindent (b) \eqref{holder} implies the local H\"older condition:
\begin{align}
\label{locholder}
&\hbox{for some $\gamma>3/4$ for all $K>0$ there is an $L_K$ so that for all }t\in [0,K]\\
&\nonumber\hbox{ and }x,X_1,X_2\in[-K,K],\ |\sigma(t,x,X_1)-\sigma(t,x,X_2)|\le L_K|X_1-X_2|^\gamma.
\end{align}
In fact it prescribes the growth rate of the H\"older constants $L_K$ (polynomial in $X$ and exponential in $x$). 
\end{remark}

In order to give a bit of intuition for Theorem~\ref{theorem:unique}, we recall
the result from~\cite{bib:mps} which
dealt with the stochastic heat equation driven by {\it coloured} noise. Let $\dot W(t,x)$ be the mean zero Gaussian noise on $\Re_+\times\Re^d$ with covariance given by 
\begin{eqnarray}
E\left[ \dot W(t,x)\dot W(s,y)\right]=\delta_0(t-s) k(x-y),
\end{eqnarray}
where
\begin{eqnarray}\label{kform}
k(x-y)&\leq& c|x-y|^{-\alpha},
\end{eqnarray}
for some $\alpha\in (0,d\wedge 2)$. Note that the white noise considered in this paper is the case $k(x)=\delta_0(x)$. It formally corresponds to  $\alpha=1$ in 
dimension $d=1$. Now let $X$ satisfies the SPDE:
\begin{equation}
\label{spde_c}
\frac{\partial}{\partial t}
X(t,x) = \frac{1}{2} \Delta X(t,x) dt + \sigma(X(t,x)) \dot{W}(x,t),
\end{equation}
with $\dot W$ being the coloured noise just described. Then the following result was proved in~\cite{bib:mps}.
\begin{theorem}[\cite{bib:mps}]
\label{thr:mps06}
For $\alpha<2\gamma -1$, pathwise uniqueness holds for~(\ref{spde_c}).
\end{theorem}
 Let 
 $\tilde{u}=X^1-X^2$ be the difference of two solutions to 
(\ref{spde_c}).  The the proof of Theorem~\ref{thr:mps06} relied on a study of the 
H\"older continuity of $\tilde u(t,\cdot)$ at points where $\tilde{u}(t,x)$ is ``small".
Let $\xi$ be the H\"{o}lder exponent of $\tilde{u}(t,\cdot)$ at such points.  The following connection 
between parameter $\xi$ and the pathwise uniqueness was shown in 
 \cite{bib:mps} (see condition (41) in the proof of Theorem 4.1 there): If 
\begin{eqnarray}
\label{equt:2}
\alpha<\xi(2\gamma-1),
\end{eqnarray}
then pathwise uniqueness holds for~(\ref{spde_c}). Hence, the better the regularity one 
has for $\tilde{u}$ near its zero set, the ``weaker" the hypotheses required for pathwise uniqueness. 
It was shown in~\cite{bib:mps} that 
at the
points $x$ where $\tilde{u}(t,x)$ is ``small", $\tilde{u}(t,\cdot)$  is H\"{o}lder continuous with any  exponent $\xi$ such that 
\begin{eqnarray}
\label{equt:1}
\xi< \frac{1-\frac{\alpha}{2}}{1-\gamma}\wedge 1. 
\end{eqnarray}
(For the precise statement of this result see Theorem~\ref{theorem:collipmod} 
in the next section.)  Note that in the case when $\alpha<2\gamma-1$, 
 (\ref{equt:1}) turns into the following condition
\begin{eqnarray}
\xi <1,
\end{eqnarray}
and this together with~(\ref{equt:2})  imply Theorem~\ref{thr:mps06}. 

Now assume $\dot W$ is white noise on $\Re_+\times \Re$ and $d=1$. This formally corresponds to the $\alpha=1$, and in this case the conditions~(\ref{equt:2}), (\ref{equt:1}) 
 can be written as 
\begin{eqnarray}
\label{equt:5}
1 &<& \xi(2\gamma-1),\\
\label{equt:3}
\xi&<& \frac{1}{2(1-\gamma)}\wedge 1. 
\end{eqnarray}
 For $\gamma\geq 1/2$, we have $\frac{1}{2(1-\gamma)}\geq 1$ and hence 
 one can take  $\xi<1$ arbitrarily close to $1$ (a proof of this is given in Theorem~\ref{theorem:lipmod} below),
 substitute it into~(\ref{equt:5}), and get  a vacuous condition for 
 pathwise  uniqueness, namely
\begin{eqnarray*}
\gamma>1.
\end{eqnarray*}
To improve on this we will need to get more refined information 
on the difference, $u$,  of two solutions to~(\ref{spde}) near the 
points $x_0$ where $u(t,x_0)\approx 0$. 
To be more precise, 
suppose one is  able to show that  
\begin{eqnarray}\label{locubehav}
|u(t,x)|\leq c|x-x_0|^{\xi}
\end{eqnarray}
for any 
\begin{eqnarray}
\label{equt:4}
\xi< \frac{1}{2(1-\gamma)}\wedge 2.
\end{eqnarray}
By substituting the upper bound for
 $\xi$ from~(\ref{equt:4}) into~(\ref{equt:5}) and doing a bit of arithmetic one 
 gets the following condition for pathwise uniqueness
\begin{eqnarray}
\gamma>3/4, 
\end{eqnarray}
which is  the result claimed in Theorem~\ref{theorem:unique}. We will in fact verify a version of \eqref{locubehav} under \eqref{equt:4} and $\gamma>3/4$.  A more detailed description
of our approach,   
is given in Section~\ref{secmainres}.

The above discussion allows us to conjecture 
 a stronger result on pathwise uniqueness for the case of 
 equations driven by a coloured noise:
\begin{conjecture}
 If
\begin{eqnarray}
\label{equt:7}
\alpha<2(2\gamma -1),
\end{eqnarray}
then pathwise uniqueness holds for~(\ref{spde_c}).
\end{conjecture}
The reasoning for this conjecture is similar to that 
 for the white noise case. Let $\tilde u$ again be the difference of two solutions 
to~(\ref{spde_c}).
Suppose that if $\tilde u(t,x_0)\approx 0$ then at the points nearby we have
\begin{eqnarray}
\tilde u(t,x)\leq c|x-x_0|^{\xi}
\end{eqnarray}
for any 
\begin{eqnarray}
\label{equt:6}
\xi< \frac{1-\frac{\alpha}{2}}{1-\gamma}\wedge 2.
\end{eqnarray}
 By substituting the upper bound for
 $\xi$ from~(\ref{equt:6}) into~(\ref{equt:2}) and simple algebra one 
 gets~(\ref{equt:7}) as  a condition for pathwise uniqueness
 for~(\ref{spde_c}). Note that~(\ref{equt:7}) can be equivalently written as 
\begin{eqnarray}
\gamma>\frac{1}{2} +\frac{\alpha}{4}.  
\end{eqnarray}

\mbox{}\\

In the next Section we give a quick proof of Theorem~\ref{theorem:nonLipschitzexistence} and then turn to the main result, Theorem~\ref{theorem:unique}.  Following the natural analogue of the Yamada-Watanabe argument for stochastic pde's, as in \cite{bib:mps}, the problem quickly reduces to one of showing that the analogue of the local time term is zero (Proposition~\ref{prop:Inbound}).  
As described above, the key ingredient here will be tight control on the spatial behaviour of the difference of two solutions, when this difference is very small, that is, when the solutions separate.  Roughly speaking, as in Yamada and Watanabe's argument we first show that solutions must separate in a gentlemanly manner and therefore cannot separate at all.  
Section~\ref{secmainres} includes a heuristic description of the method and 
further
explanation of why $\gamma=3/4$ is critical in our approach.  It also gives an outline of the contents of the entire paper.  

\medskip

\noindent{\bf Convention on Constants.} Constants whose value is unimportant and may change from line to line are denoted $c_1,c_2,\dots$, while constants whose values will be referred to later and appear initially in say, Lemma~i.j are denoted $c_{i.j}$ or $C_{i.j}$. 
\medskip
\paragraph{Acknowledgements.} The second author thanks the Technion for hosting him during a visit
where some of this research was carried out.  This project was initiated during 
the visit of the first  author to the UBC where he participated in a Workshop on SPDE's sponsored by PIMS and thanks go to the Pacific Institute for the Mathematical Sciences for its support. 

\section{Proof of Theorems \ref{theorem:nonLipschitzexistence} and \ref{theorem:unique}} \label{secmainres}

\paragraph{Proof of Theorem \ref{theorem:nonLipschitzexistence}} 
This is standard so we only give a sketch and set $b\equiv0$ for simplicity.  By taking weak limits as $T\to \infty$ we may assume $R_0(T)=R_0$ is independent of $T$.  Choose a symmetric $\psi_n\in C_c^\infty$ so that $0\le \psi_n\le 1$, $\Vert\psi'_n\Vert_\infty\le 1$, $\psi_n(x)=1$ if $|x|\le n$ and $\psi_n(x)=0$ if $|x|\ge n+2$.  Let 
\[\sigma_n(t,x,X)=\int\sigma(t,x,X')p_{2^{-n}}(X'-X)dX'\psi_n(X).\]
It is easy to then check the following:
\begin{equation}\label {unifgrowth}
|\sigma_n(t,x,X)|\le 2c_{\ref{growthcond}}(1+|X|),
\end{equation}
\begin{equation}\label{lipcond}
|\sigma_n(t,x,X')-\sigma_n(t,x,X)|\le c_n|X'-X|,
\end{equation}
and
\begin{align}\label{convrate}
|\sigma_n(t,x,X)-\sigma(t,x,X)|&\le c_{\ref{convrate}}\Bigl[e^{R_1|x|}(|X|^{R_2}+1)2^{-n\gamma/2}+(1+|X|)(1-\psi_n(|X|)\Bigr]\\
\nonumber&\to0\ \hbox{ uniformly on compacts as }n\to\infty.
\end{align}
Use \eqref{unifgrowth}, \eqref{lipcond} and Theorem~2.2 of \cite{bib:shiga94} to see there are solutions
$X^n$ to $\eqref{wspde}_n$ (all with respect to $W$)--here $\eqref{wspde}_n$ is \eqref{wspde} but with $\sigma_n$ in place of $\sigma$.   Now argue as in Section~6 of \cite{bib:shiga94} (see the proof of Theorem~2.2) or in the derivation of Theorem~1.2 of \cite{bib:mps} (the present white noise setting simplifies those arguments) to see that $\{X^n\}$ is tight in $C(\Rp,C_{tem})$.  More specifically, using the growth condition \eqref{growthcond}, it is straightforward to carry over the proof of Proposition~1.8(a) of \cite{bib:mps} (see Lemmas A.3 and A.5 of that paper) and show
\begin{equation} \label{LpXbnd}
\hbox{for all }T,\lambda,p>0,\quad\sup_n E(\sup_{0\le t\le T}\sup_{x\in\Re}|X^n(t,x)|^pe^{-\lambda|x|})<\infty.
\end{equation}
The above bounds in turn give uniform bounds on the $p$th moments of the space-time increments of $X^n$ (see Lemma~A.4 of \cite{bib:mps}) and hence tightness.
Indeed, the orthogonality of white noise makes all these calculation somewhat easier.  By Skorohod's theorem we may assume $X^{n_k}$ converges a.s. to $X$ in $C(\Rp,C_{tem})$ on some probability space. It is now easy to use \eqref{convrate} to see that (perhaps on a larger space), $X$ solves \eqref{wspde}. \gdm

Next consider Theorem~\ref{theorem:unique} and assume its hypotheses throughout.  By Remark~\ref{holderhyp}(a) decreasing $\gamma$ only weakens the hypotheses and so we may, and shall,  assume that 
\begin{equation}\label{gammabound}
3/4<\gamma<1.
\end{equation}
Let $X^{1}$ and $X^{2}$ be two solutions of (\ref{wspde}) on $(\Omega,\cF,\cF_t,\Pe)$ with sample
paths in $C(\IR_+,C_{tem})$ a.s.,
 with the same
initial condition, $X^{1}(0)=X^{2}(0)=X_0\in C_{tem}$, and of course the same noise $W$.  
For adapted processes with sample paths in $C(\Rp,C_{tem})$, \eqref{wspde} is equivalent to the distributional form of~(\ref{spde}) (see Theorem~2.1 of \cite{bib:shiga94}).  That is, for $i=1,2$ and $\Phi \in C_{c}^{\infty}(\IR):$
\begin{align}
\label{weakheatb}
\int_{\IRsm} X^i(t,x) \Phi(x)dx
= &\int_{\IRsm} X^i_{0}(x)  \Phi(x)dx
+\int_{0}^{t}\int_{\IRsm} X^i(s,x) \frac{1}{2} \Delta\Phi(x) dx ds
\\
\nonumber
&+\int_{0}^{t}\int_{\IRsm}\sigma(s,x,X^i(s,x))\Phi(x) W(ds, dx)\\
\nn&+\int_{0}^{t}\int_{\IRsm}b(s,x,X^i(s,x))\Phi(x)  dxds
\quad\forall t\ge 0 \quad a.s.
\end{align}
Let
\begin{equation}\label{TKdef}
T_K=\inf\{s\ge 0:\sup_y(|X^1(s,y)|\vee|X^2(s,y)|)e^{-|y|}>K\}\wedge K.
\end{equation}  

We first show that \eqref{holder} may be strengthened to 
\begin{align}
\label{holder'}
&\hbox{for some $1>\gamma>3/4$ there are $R_0,R_1\ge1$ so that for all }t\ge 0\\
&\nonumber\hbox{ and all }
(x,X,X')\in\Re^3,\ |\sigma(t,x,X)-\sigma(t,x,X')|\le R_0e^{R_1|x|}|X-X'|^\gamma.
\end{align}
Assume that Theorem~\ref{theorem:unique} holds under \eqref{holder'} 
and that $\sigma$ satisfies \eqref{holder}. 
 Define 
\[\sigma_K(t,x,X)=\sigma(t,x,(X\vee(-Ke^{|x|}))\wedge Ke^{|x|})1(t\le K).\]
Then 
\[|\sigma_K(t,x,X)-\sigma_K(t,x,X')|\le R_0(K)e^{R_1|x|}(1+2Ke^{|x|})^{R_2}|X-X'|^\gamma,\]
and so \eqref{holder'} holds with $R_1+R_2$ in place of $R_1$ (the restriction that $R_i\ge1$ is for convenience and is no restriction).  Providing that for $\lambda=1$, $\Vert X_0\Vert_\lambda<K$, we have
\begin{equation}\label{sigmas}
\sigma(t,x,X^i(t,x))=\sigma_K(t,x,X^i(t,x))\hbox{ for all }x \hbox{ and }t\le T_K.
\end{equation}
Therefore $\sigma_K$ satisfies\eqref{holder'} and of course \eqref{growthcond}.  So we may apply Theorem~\ref{strong} with $\sigma_K$ in place of $\sigma$.  Using the law $P_{K,X_0}$ of $(X,W)$ on $C(\Rp,C_{tem})^2$ (Borel in $X_0$) it is easy now to continue the solutions $X^i$ to $\eqref{weakheatb}_K$ (the $K$ reminds us we are dealing with $\sigma_K$) beyond $T_K$ and construct solutions $\tilde X^i$, $i=1,2$ to $\eqref{weakheatb}_K$ starting at $X_0$ such that $(\tilde X^1(\cdot\wedge T_K), \tilde X^2(\cdot\wedge T_K))$ is equal in law to $(X^1(\cdot\wedge T_K), X^2(\cdot\wedge T_K))$.  By pathwise uniqueness in $\eqref{weakheatb}_K$ we get $\tilde X^1=\tilde X^2$ and so $X^1(\cdot\wedge T_K)=X^2(\cdot\wedge T_K)$. Letting $K\to\infty$ gives $X^1=X^2$, as required. 

We now follow the approach in Section 2 of \cite{bib:mps} and reduce the theorem to showing the analogue of the ``local time term" in the Yamada-Watanabe proof is zero. Let 
$$a_n=\exp\{-n(n+1)/2\}$$ so that 
\begin{equation}\label{an}
a_{n+1}=a_ne^{-n-1}=a_na_n^{2/n}.
\end{equation}
Define functions $\psi_{n} \in C^{\infty}_{c}(\IR)$ such that
$supp(\psi_{n}) \subset  (a_{n}, a_{n-1})$, and
\begin{equation}
\label{psicond}
0 \leq \psi_{n}(x) \leq  {2\over nx}
\quad \mbox{ for all $x \in \IR$ as well as } \quad \int_{a_{n}}^{a_{n-1}} \psi_{n}(x) dx =1.
\end{equation}
Finally, set
\begin{equation}
\label{def:phi}
\phi_{n}(x) = \int_{0}^{|x|} \int_{0}^{y} \psi_{n}(z) dz dy.
\end{equation}
From this it is easy to see that $\phi_{n}(x) \uparrow |x|$ uniformly in $x$.
Note that each $\psi_{n}$, and thus also each $\phi_{n}$, is
identically zero in a neighborhood of zero. This implies that
$\phi_{n} \in C^{\infty}(\IR)$ despite the absolute value in its
definition. We have
\begin{align}
\label{phidiff1}
\phi_{n}'(x) &= \sgn(x) \int_{0}^{|x|}  \psi_{n}(y) dy,\\
\label{phidiff2}
\phi_{n}''(x) &=  \psi_{n}(|x|).
\end{align}
Thus, $|\phi_{n}'(x)| \leq 1$, and
$\int \phi_{n}''(x) h(x) dx \rightarrow h(0)$ for any function $h$ which
is continuous at zero.

Define 
\[{u} \equiv X^{1} - X^{2}.\] Let
$\Phi
\in C_{c}^{\infty}(\IR)$ satisfy $0\le \Phi\le 1$, $supp(\Phi) \subset (-1,1)$
and $\int_{\Resm}\Phi(x) dx=1$, and set
$\Phi_{x}^{m}( y) = m\Phi(m(x - y))$.
\noindent
Let $\la \cdot, \cdot \ra$ denote the scalar product on $\Le^{2}(\IR).$
By applying It\^{o}'s Formula to the semimartingales
$\la X^i_{t}, \Phi^{m}_{x} \ra$ in (\ref{weakheatb})
it follows that
\begin{align*}
\phi_{n}&(\la {u}_{t}, \Phi_{x}^{m} \ra)\\
=& \int_{0}^{t} \int_{\IRsm} \phi_{n}'(\la {u}_{s}, \Phi_{x}^{m} \ra)
\left( \sigma(s,y,X^{1}(s,y)) - \sigma(s,y,X^{2}(s,y)) \right)
\Phi_{x}^{m}(y)  W(ds,dy)\\
&+\int_{0}^{t} \phi_{n}'(\la{u}_{s}, \Phi_{x}^{m} \ra)
\la {u}_{s}, \frac{1}{2} \Delta \Phi_{x}^{m} \ra ds\\
&+ \frac{1}{2}
\int_{0}^{t} \int_{\IRsm}
\psi_{n}(|\la{u}_{s}, \Phi_{x}^{m} \ra|)
\left( \sigma(s,y,X^{1}(s,y)) - \sigma(s,y,X^{2}(s,y)) \right)^2\\
& \phantom{AAAAAAAAAAAAAAA}
\times\Phi_{x}^{m}(y)^2 dy ds\\
&+ \int_{0}^{t} \int_{\IRsm} \phi_{n}'(\la {u}_{s}, \Phi_{x}^{m} \ra)
\left( b(s,y,X^{1}(s,y)) - b(s,y,X^{2}(s,y)) \right)
\Phi_{x}^{m}(y)  dyds.\\
\end{align*}
We integrate this function of $x$ against another non-negative test
function
$\Psi \in C^{\infty}_{c}([0,t_0]\times\IR)$ ($t_0\in(0,\infty)$).  Choose $K_1\in\N$ large so that for $\lambda=1$, 
\begin{equation}
\label{Gamma}
\Vert X_0\Vert_\lambda<K_1\hbox{ and }
\Gamma\equiv\{x:\Psi_s(x)>0\ \exists s\le t_0\}
\subset (-K_1,K_1).
\end{equation}
 We then obtain by the
classical and stochastic versions of Fubini's Theorem (see Theorem~2.6 of \cite{bib:wal86} for the latter), and arguing as
in the proof of Proposition II.5.7 of \cite{bib:per02} to handle the time
dependence in $\Psi$, that for any $t\in [0,t_0]$,
\begin{align}
\label{eq:Iparts}
\la &\phi_{n}(\la {u}_{t}, \Phi_{.}^{m} \ra), \Psi_t\ra\\
\nonumber
&= \int_{0}^{t}\int_{\IRsm}
\la\phi_{n}'(\la {u}_{s}, \Phi_{\cdot}^{m} \ra) \Phi_{\cdot}^{m}(y)
 , \Psi_s \ra \left( \sigma(s,y,X^{1}(s,y)) - \sigma(s,y,X^{2}(s,y)) \right)
W(ds,dy)\\
\nonumber
&\quad+ \int_{0}^{t} \la \phi_{n}'(\la {u}_{s}, \Phi_{.}^{m} \ra)
\la {u}_{s}, \frac{1}{2} \Delta \Phi_{.}^{m} \ra, \Psi_s \ra ds\\
\nonumber
&\quad+ \frac{1}{2}
\int_{0}^{t} \int_{\IRsm^{2}}
\psi_{n}(|\la {u}_{s}, \Phi_{x}^{m} \ra|)
\left( \sigma(s,y,X^{1}(s,y)) - \sigma(s,y,X^{2}(s,y)) \right)^2\\
\nonumber
&\phantom{AAAAAAAA}
\times\Phi_{x}^{m}(y)^2
dy \Psi_s(x) dx ds+\int_0^t\la\phi_n(\la
u_s,\Phi^m_\cdot\ra) ,\dot\Psi_s\ra\,ds\\
\nn&+ \int_{0}^{t}\int_{\IRsm}
\la\phi_{n}'(\la {u}_{s}, \Phi_{\cdot}^{m} \ra) \Phi_{\cdot}^{m}(y)
 , \Psi_s \ra \left( b(s,y,X^{1}(s,y)) - b(s,y,X^{2}(s,y)) \right)
dyds\\
\nonumber
&\equiv I_{1}^{m,n}(t) + I_{2}^{m,n}(t) + I_{3}^{m,n}(t)+I_{4}^{m,n}(t)+I_5^{m,n}(t).\\
\end{align}
The expectation condition in Walsh's Theorem~2.6 may be realized by localization, using the stopping times $\{T_K\}$.

Set $m_n=a_{n-1}^{-1/2}=\exp\{(n-1)n/4\}$ for $n\in\N$. 
In the integral defining $I_3^{m_{n+1},n+1}$ we may assume $|x|\le K_1$ by \eqref{Gamma} and so $|y|\le K_1+1$.  Let $K\in\N^{\ge K_1}$.  If  $s\le T_{K}$, then for such a $y$, 
\[|X^i(s,y)|\le K e^{|y|}\le Ke^{(K_1+1)}\hbox{ for }i=1,2.\]  
Therefore \eqref{locholder}, \eqref{psicond} and \eqref{an} show that if $K'= Ke^{(K_1+1)}(\ge K_1+1)$, then for all $t\in [0,t_0]$, 
\begin{eqnarray}
\nn&&I_3^{m_{n+1},n+1}(t\wedge T_{K})\\
\nn&&\le {1\over 2}\int_0^{t\wedge T_{K}}\int\int2(n+1)^{-1}|\la u_s,\Phi_x^{m_{n+1}}\ra|^{-1}1(a_{n+1}<|\la u_s,\Phi_x^{m_{n+1}}\ra|<a_n)\\
\nn&&\phantom{\le {1\over 2}\int_0^{t\wedge T_{K}}\int\int 2(n+1)^{-1}|\la u_s,\Phi_x^{m_{n+1}}\ra|^{-1}}\times L_{K'}^2|u(s,y)|^{2\gamma}m_{n+1}\Phi_x^{m_{n+1}}(y)\Psi_s(x)dydxds\\
\nn&&\le L^2_{K'}a_{n+1}^{-1}a_n^{-1/2}\int_0^{t\wedge T_{K}}\int\int1(a_{n+1}<|\la u_s,\Phi_x^{m_{n+1}}\ra|<a_n)|u(s,y)|^{2\gamma}\Phi_x^{m_{n+1}}(y)\Psi_s(x)dydxds\\
\label{I3bnd1}&&\le L^2_{K'}a_n^{-3/2-2/n}\int_0^{t\wedge T_{K}}\int\int1(a_{n+1}<|\la u_s,\Phi_x^{m_{n+1}}\ra|<a_n)|u(s,y)|^{2\gamma}\Phi_x^{m_{n+1}}(y)\Psi_s(x)dydxds.
\end{eqnarray}

We define
\begin{equation}\label{Indef}
I^n(t)=a_n^{-3/2-2/n}\int_0^{t}\int\int1(|\la u_s,\Phi_x^{m_{n+1}}\ra|<a_n)|u(s,y)|^{2\gamma}\Phi_x^{m_{n+1}}(y)\Psi_s(x)dydxds.
\end{equation}

\begin{proposition}\label{prop:Inbound}
Suppose $\{U_{M,n,K}:M,n,K\in\N, K\ge K_1\}$ are $\cF_t$-stopping times such that for each $K\in\N^{\ge K_1}$, 
\begin{align*}
&\hbox{$(H_1)\quad U_{M,n,K}\le T_K$,  $U_{M,n,K}\uparrow T_K$ as $M\to\infty$ for each $n$,  and
}\phantom{AAAAAAAAAAAAAAAAAAAAAAAAA}\\
&\phantom{AAAAAAAAAAAA}\hbox{ $\lim_{M\to\infty}\sup_nP(U_{M,n,K}<T_K)=0$,}\\
&\hbox{and}\\
&\hbox{$(H_2)\quad$ For all $M\in\N$, $\lim_{n\to\infty}E(I^n(t_0\wedge U_{M,n,K}))=0.$}
\end{align*}
Then the conclusion of Theorem~\ref{theorem:unique} holds.
\end{proposition}
\paragraph{Proof.}  We adapt the reasoning in Lemma 2.2 of \cite{bib:mps} for the coloured noise setting to our white noise driven equation.   As in \eqref{LpXbnd} we have 
\begin{equation} \label{Lpbnd}
\hbox{for all }T,\lambda,p>0,\quad E(\sup_{0\le t\le T}\sup_{x\in\Re}|u(t,x)|^pe^{-\lambda|x|})<\infty.
\end{equation}
Let
\[Z_n(t)=\int\phi_n(\langle u_t,\Phi_x^{m_n}\rangle)\Psi_t(x)\,dx.\]
Fix $K\in\N^{\ge K_1}$ and $0\le t\le t_0$.  Note that since $0\le \phi_n(z)\le |z|$ and $\Psi\ge 0$, 
\begin{align}
\nn0\le Z_n(t\wedge T_K)&\le \int \int|u(t\wedge T_K,y)|\Phi_x^{m_n}(y)\Psi(t\wedge T_K,x)\,dydx\\
\nn&\le 2K\int\int e^{|y|}\Phi_x^{m_n}(y)\Psi(t\wedge T_K,x)dy 1(|x|\le K_1)dx\\
\label{Zbnd1}&\le 2Ke^{K_1+1}c_1(\Psi).
\end{align}
With \eqref{Lpbnd} in hand, the proof of Lemma~2.2(a) of \cite{bib:mps} is easily adapted (again it is in fact easier) to show
\begin{equation}\label{I1bnd}\{I_1^{m_n,n}(s):s\le t_0\}\hbox{ is an $\Le^2$-bounded sequence of $\Le^2$ martingales.}
\end{equation}
The proof of Lemma~2.2(b) of \cite{bib:mps} applies directly to show that
\[I_2^{m_n,n}=I_2^{m_n,n,1}+I_2^{m_n,n,2},\]
where for any stopping time $T$,
\begin{equation}\label{I21bnd}
I_2^{m_n,n,1}(t\wedge T)\to\int_0^{t\wedge T}\int |u(s,x)|{1\over 2}\Delta\Psi_s(x)dxds\hbox{ in }\Le^1\hbox{ as }n\to\infty
\end{equation}
(again the key bound here is \eqref{Lpbnd}), and we have the one-sided bound
\begin{equation}\label{I22bnd}
I_2^{m_n,n,2}(s)\le {a_n\over n}C(\Psi)\hbox{ for all }s\le t_0\hbox{ and }n.
\end{equation}
(In the notation of \cite{bib:mps}, $I_2^{m_n,n,1}(s)=\int_0^sI_{2,3}^{m_n,n}(r)dr$ and $I_2^{m_n,n,2}(s)=\int_0^sI^{m_n,n}_{2,1}(r)+I^{m_n,n}_{2,2}(r)dr$.)  The proof of Lemma~2.2(c) of \cite{bib:mps} also applies directly to show that for any stopping time $T$,
\begin{equation}\label{I4bnd}
I_4^{m_n,n}(t\wedge T)\to\int_0^{t\wedge T}\int |u(s,x)|\dot\Psi_s(x)dxds\hbox{ in }\Le^1\hbox{ as }n\to\infty.
\end{equation}
Since $|\phi'_n|\le 1$, \eqref{bLip} implies that for a stopping time $T$,
\begin{equation}\label{I5bound}
I_5^{m_n,n}(t\wedge T)\le B\int^{t\wedge T}_0\int\int |u(s,y)|\Phi_x^{m_n}(y)\Psi_s(x)dydxds\equiv B\tilde I_5^n(t\wedge T).
\end{equation}
It follows easily from \eqref{Lpbnd} that $\{\tilde I_5^n(t_0):n\in\N\}$ is $\Le^2$-bounded and,
as $n\to \infty$, 
\begin{equation}\label{I5conv}
\tilde I_5^n(t\wedge T)\to \int _0^{t\wedge T}\int |u(s,x)|\Psi_s(x)dxds\hbox{ a.s. and hence in }\Le^1\hbox{ by the above.}
\end{equation}
Let $\ve>0$.  Then $(H_1)$, \eqref{Zbnd1},\eqref{I1bnd},\eqref{I21bnd} and \eqref{I4bnd} show that, by a standard result for uniformly integrable random variables, there is an $M_0$ so that 
\begin{align}\label {errbnd}
&\sup_n E\Bigl(\Bigl(|Z_n(t\wedge T_K)|+|I_1^{m_n,n}(t\wedge T_K)|+|I_2^{m_n,n,1}(t\wedge T_K)|+|I_4^{m_n,n}(t\wedge T_K)|\Bigr)\\
\nn&\phantom{\sup_n E\Bigl((|Z_n(t\wedge T_K)|+|I_1^{m_n,n,1}(t\wedge T_K)|}
\times1(U_{M,n,K}<T_K)\Bigr)<\ve \hbox{ for all }M\ge M_0.
\end{align}
From \eqref{eq:Iparts}, \eqref{I5bound}, the non-negativity of $I_3^{m_n,n}$ and $\tilde I^n_5$, and Fatou's Lemma, we have for $M\ge M_0$,
\begin{align*}
&E\Bigl(\int|u(t\wedge T_K,x)|\Psi_{t\wedge T_K}(x)dx\Bigr)\\
&\le \liminf_{n\to\infty}E(Z_n(t\wedge T_K)1(U_{M,n,K}=T_K))+E(Z_n(t\wedge T_K)1(U_{M,n,K}<T_K))\\
&\le \liminf_{n\to\infty}E(I^{m_n,n}_1(t\wedge T_K))+E(I^{m_n,n,1}_2(t\wedge T_K))+E(I^{m_n,n,2}_2(t\wedge T_K)1(U_{M,n,K})=T_K))\\
&\quad+E(I_3^{m_n,n}(t\wedge U_{M,n,K}))+E(I^{m_n,n}_4(t\wedge T_K))+BE(\tilde I_5^n(t\wedge T_K))\\
&\quad-E(I_1^{m_n,n}(t\wedge T_K)1(U_{M,n,K}<T_K))
-E(I_2^{m_n,n,1}(t\wedge T_K)1(U_{M,n,K}<T_K))\\
&\quad-E(I_4^{m_n,n}(t\wedge T_K)1(U_{M,n,K}<T_K))
+E(Z_n(t\wedge T_K)1(U_{M,n,K}<T_K))\\
&\le \liminf_{n\to\infty}E\Bigl(\int_0^{t\wedge T_K}\int |u(s,x)|{1\over 2}\Delta \Psi_s(x)dxds\Bigr)+{a_n\over n}C(\Psi)\\
&\quad +E\Bigl(\int_0^{t\wedge T_K}\int |u(s,x)|\dot\Psi_s(x)dxds\Bigr)+BE\Bigl(\int _0^{t\wedge T_K}\int |u(s,x)|\Psi_s(x)dxds\Bigr)+\ve,
\end{align*}
by \eqref{I1bnd}, \eqref{I21bnd}, \eqref{I22bnd}, $(H_2)$ (together with the bound \eqref{I3bnd1}),\eqref{I4bnd} \eqref{I5conv},  and \eqref{errbnd}, respectively.  Let $\ve\downarrow 0$ to see that
\begin{align*}
E\Bigl(\int |u(t\wedge T_K,x)|\Psi_{t\wedge T_K}(x)dx\Bigr)&\le E\Bigl(\int_0^{t\wedge T_K}\int |u(s,x)|({1\over 2}\Delta\Psi_s(x)+\dot\Psi_s(x)+B\Psi_s(x))dxds\Bigr).
\end{align*}
Let $K\to \infty$ and use Dominated Convergence (recall \eqref{Lpbnd}) on each side to conclude that
\[ \int E(|u(t,x)|)\Psi_t(x)dx\le \int_0^t \int E(|u(s,x)|)({1\over 2}\Delta\Psi_s(x)+\dot\Psi_s(x)+B\Psi_s(x))dxds,\quad 0\le t\le t_0.\]
This gives (34) of \cite{bib:mps} with an additional drift term $B\Psi_s(x)$.  One proceeds exactly as in Section 3 of that reference, using the semigroup $e^{Bt}P_t$ in place of $P_t$, to see that since $E(|u(t,x)|)$ is a finite (by \eqref{Lpbnd})
non-negative subsolution of the heat equation with initial data zero, therefore $E(|u(t,x)|)=0$ and 
so $X^1=X^2$ by continuity of paths.  \gdm

The construction of $\{U_{M,n,K}\}$ and verification of $(H_1)$ and $(H_2)$ will be the objective of the rest of this work.  

\medskip

\noindent{\bf Notation.} For $t,t'\ge 0$ and $x,x'\in\Re$ let $d((t,x),(t',x'))=\sqrt{|t'-t|}+|x'-x|$.  

\medskip

Note that the indicator function in the definition of $I^n$ implies there is an $\hat x_0\in(x-\sqrt a_n, x+\sqrt a_n)$ such that $|u(s,\hat x_0)|\le a_n$.  If we could take $\hat x_0=y$ we could bound $I_n(t)$ by $C(t)a_n^{-3/2-2/n+2\gamma}$, and $(H_1)$ and $(H_2)$ would follow immediately with $U_{M,n,K}=T_K$.  (The criticality of $3/4$ in this argument is illusory as it follows from our choice of $m_n$.)
The hypotheses of Proposition~\ref{prop:Inbound} now turn on getting good bounds on $|u(s,y)-u(s,\hat x_0)|$.  The standard $1/2-\ve$-H\"older modulus\thinspace\footnote{Although this is well-known ``folklore'' result we were not able to find the exact reference. One can easily check that the estimates in the proof of 
 Corollary~3.4 in~\cite{bib:wal86} give $1/2-\ve$-H\"older spatial modulus; 
 similarly the result of~\cite{SS02} can be immediately extended to cover the white noise case; in both 
 works   the Lipschitz assumptions on noise coefficients can be relaxed to linear growth assumptions  and 
 the proofs still go through. }
not surprisingly, gives nothing. In \cite{bib:mps} this was refined to a $1-\ve$-H\"older modulus {\it near points where $u$ is small} as we now describe.  Let
\begin{align*} Z(N,K)(\omega)=\{(t,x)\in [0,T_K]\times [-K,K]:& \hbox{ there is a }(\hat t_0,\hat x_0)\in [0,T_K]\times \Re \hbox{ such that }\\
&\ d((\hat t_0,\hat x_0),(t,x))\le 2^{-N},\hbox{ and }|u(\hat t_0,\hat x_0)|\le 2^{-N}\}.
\end{align*}

 Return now to the
SPDE driven by coloured noise~(\ref{spde_c}) from 
 Section 1.  Let $\tilde u$ be the difference of two solutions of~(\ref{spde_c}) 
and $\alpha\in (0,2\wedge d)$ be the covariance kernel exponent as in \eqref{kform}. 
Let $\tilde Z(N,K)$  be defined as $Z(N,K)$ with $\tilde u$ instead of $u$. 
The following improved modulus of continuity was proved in \cite{bib:mps} (see Theorem 4.1 and the first two paragraphs of the proof of Corollary 4.2 in that reference).

\begin{theorem}\label{theorem:collipmod} For each $K\in \N$ and $0<\xi<{1-{\alpha\over 2}\over 1-\gamma}\wedge 1$ there is an $N_0=N_0(\xi,K,\omega)\in\N$ a.s. such that for all natural numbers $N\ge N_0$ and all 
$(t,x)\in \tilde Z(N,K)$, 
\[ d((t',x'),(t,x))\le 2^{-N}\hbox{ and }t'\le T_K \hbox{ implies }|\tilde u(t',x')-\tilde u(t,x)|\le 2^{-N\xi}.\]
\end{theorem}

In the white noise setting the result holds with $\alpha=1$.  Recall $u$ is the difference of two solutions   to \eqref{spde}.

\begin{theorem}\label{theorem:lipmod} Assume $\gamma\ge1/2$.  For each $K\in \N$ and $\xi\in(0,1)$ there is an $N_0=N_0(\xi,K,\omega)\in\N$ a.s. such that for all natural numbers $N\ge N_0$ and all $(t,x)\in Z(N,K)$, 
\[ d((t',x'),(t,x))\le 2^{-N}\hbox{ and }t'\le T_K \hbox{ implies }|u(t',x')-u(t,x)|\le 2^{-N\xi}.\]
\end{theorem}
\paragraph{Proof.} The proof of Theorem~\ref{theorem:collipmod} applies with $\alpha=1$ (the dependance of $\sigma$ on $t,x$ alters nothing in the proof).  In fact it is now
considerably simpler because of the orthogonality of white noise increments.    The required tools are Lemma~\ref{Jbnd} and Lemma~\ref{ptbnds} below.  With this choice of $\alpha$ and $\gamma\ge 1/2$, the upper bound on $\xi$ in Theorem~\ref{theorem:collipmod} becomes $1$.  In \cite{bib:mps} there was no drift term, but the calculations for the Lipschitz drift term are simpler still. Here one uses H\"older's inequality to utilize the $\Le^2$ bounds in Lemma~\ref{ptbnds}.
\gdm

The proof of Theorem~\ref {theorem:unique}
is long and involved so before descending into the technical details of the derivation of $(H_1)$ and $(H_2)$, we now give a heuristic description of the method with $b\equiv 0$ throughout, and also try to explain why $\gamma=3/4$ is critical in our approach.  The choice of $m_n=a_{n-1}^{-1/2}$ appears arbitrarily in the above so let us for the moment set $m_n=a_{n-1}^{-\alpha_0}\approx a_n^{-\alpha_0}$ for some $\alpha_0>0$.  
  $(H_1)$ and $(H_2)$ are delicate ways of ensuring $I^n(t)$ approaches zero as $n\to\infty$ and so our goal is to show that 
\begin{align}\label{Intozero}
I^n(t)&\approx a_n^{-1-\alpha_0}\int_0^t\int\int 1(|\langle u_s,\Phi_x^{m_{n+1}}\rangle |<a_n)|u(s,y)|^{2\gamma}\Phi_x^{m_{n+1}}(y) \Psi_s(x)dy dx ds\\
\nonumber&\to 0\hbox{ as n }\to \infty.
\end{align}
We have taken (and will take) some small liberties  with the ``local time term" $I^n(t)$ (with this new choice of $m_n$) in the first line.  In the integrand in \eqref{Intozero} the variable $y$ must be within $2a_n^{\alpha_0}$ of a point $\hat x_0$ where $|u(s,\hat x_0)|<a_n$.  If we simply replace $y$ with $\hat x_0$, $I^n(t)$ is at most
\begin{align*}
a_n^{-1-\alpha_0}\int_0^t\int\int a_n^{2\gamma}\Phi_x^{m_{n+1}}(y)\Psi_s(x)dydxds
\le Cta_n^{2\gamma-1-\alpha_0}\to 0,\\
\end{align*}
 if $\gamma>1/2$ and $\alpha_0$ is small enough. This is a bit too crude but shows it will be crucial to get good estimates on $u(s,\cdot)$ near points where it is small (and also shows we are already forced to assume $\gamma>1/2$).  Theorem~\ref{theorem:lipmod} implies that 
 \begin{equation}\label{uliponZ}\gamma\ge1/2\hbox{ implies } u(t,\cdot) \hbox{ is $\xi$-H\"older continuous near its zero set for $\xi<1$},
 \end{equation}
 and so allows us to bound $|u(s,y)-u(s,\hat x_0)|$.  Use this in \eqref{Intozero} and take $0<\alpha_0\le 1$ to bound $I^n(t)$ by
 \begin{align*}
 a_n^{-1-\alpha_0}\int_0^t \int\int Ca_n^{\alpha_0\xi 2\gamma}\Phi_x^{m_{n+1}}(y)\Psi_s(x)dydxds\\
 \le Ct a_n^{-1+\alpha_0(\xi2\gamma-1)}
 \to 0\hbox{ as }n\to\infty,
 \end{align*}
if $\gamma>1$ and we choose $\alpha_0,\xi$ close to one.  Of course $\gamma>1$ is not a viable choice but this shows we are now getting close, and in fact in the coloured noise setting of \cite{bib:mps} the above argument sufficed for the results there, although there was some work to be done to implement this idea carefully.  
 
To increase our control on $u(s,\cdot)$ near its zero set we will improve \eqref{uliponZ} to
\begin{equation}\label{uC2-onZ2} \gamma>3/4\hbox{ implies $u'(s,\cdot)$ is $\xi$-H\"older on }\{x: u(s,x)\approx u'(s,x)\approx 0\} \hbox{ for }\xi <1,
\end{equation}
where $u'$ denotes the spatial derivative.  Corollary~\ref{u1'modulus} below with $m=\bar m+1$ is the closest result which comes to a formal statement of the above, although the condition on $\gamma$ is implicit.  

We first make the case that for $\gamma<3/4$, we cannot expect the following slight strengthening of \eqref{uC2-onZ2}:
\begin{equation}\label{uC2onZ2}
u(s,\cdot)\hbox{ is 
 $C^2$ on }\{x:u(s,x)\approx u'(s,x)\approx 0\}.
\end{equation}
A formal differentiation of \eqref{wspde} (recall $b\equiv 0$ and $u$ is the difference of the $X^i$'s) gives for $u(t,x)\approx u'(t,x)\approx 0$,
\begin{equation}\label{d2form} u''(t,x)=\int _0^t\int p''_{t-s}(y-x)[\sigma(s,y,X^1(s,y))-\sigma(s,y,X^2(s,y))]W(ds,dy).
\end{equation}
If $\sigma$ is a Weierstrass-type function that realizes its H\"older modulus at typical points we have 
\[|\sigma(s,y,X^1(s,y))-\sigma(s,y,X^2(s,y))|\approx L|u(s,y)|^\gamma,\]
and for $s<t$ and very close to $t$, we have by a Taylor series expansion in space,
\[ |u(s,y)|\approx |u''(s,x)|{(y-x)^2\over 2}.\]
Use these approximations in the finite square function associated with the right hand side of \eqref{d2form} and conclude
that
\begin{align*}
\infty&>\int _{t-\delta}^t\int p''_{t-s}(y-x)^2 \Bigl[{|u''(s,x)|(y-x)^2\over 2}\Bigr]^{2\gamma}dyds\\
&\approx c |u''(t,x)|^{2\gamma}\int _{t-\delta}^t \int p_{t-s}(z)^2[z^2(t-s)^{-2}-(t-s)^{-1}]^2 z^{4\gamma}dzds\\
&\approx c |u''(t,x)|^{2\gamma} \int _{t-\delta}^t (t-s)^{2\gamma-{5\over 2}}ds,
\end{align*}
which implies $\gamma>3/4$.  

We next show how \eqref{uC2-onZ2} will lead to \eqref{Intozero}.  Taking further liberties with 
$I^n(t)$ and recalling $m_n\approx a_n^{\alpha_0}$, we get
\begin{align}\label{betadecomp}I^n(t)&\approx a_n^{-1-\alpha_0}\sum_\beta \int _0^t \int\int 1(|u(s,x|\le a_n, u'(s,x)\approx \pm a_n^\beta)|u(s,y)|^{2\gamma}\Phi_x^{m_{n+1}}(y)\Psi_s(x)dydxds\\
\nonumber&\equiv\sum_\beta I^n_\beta(t),
\end{align}
where $\sum_\beta$ indicates we are summing over a finite grid $\beta_i\in[0,\bar\beta]$ ($\bar\beta$ to be determined below) and are bounding $u'(s,x)$ in the appropriate grid interval and its mirror image in the origin.  As the sum is finite we may fix
$\beta\in(0,\bar\beta]$ and consider only $u'(s,x)\approx a_n^\beta$.   The value $\beta=0$ is a bit special but should be clear from the argument below.   A Taylor series expansion and \eqref{uC2-onZ2} with $\xi\approx 1$ show that for $y$ as in the integrand of $I^n_\beta(t)$,
\begin{align*}|u(s,y)|&\le |u(s,x)|+(|u'(s,x)|+L|y-x|^\xi)|y-x|\\
&\le a_n+a_n^{\beta+\alpha_0}+La_n^{\alpha_0(\xi+1)}\\
&\le Ca_n^{({\xi\over 2}\wedge\beta)+{1\over 2}},
\end{align*}
where a comparison of the first and last terms in the second line leads naturally to $\alpha_0=1/2$.
Substitute this into the integrand of $I^n_\beta$, integrate out $y$,  and conclude
\begin{equation}\label{Inbetabnd}
I_\beta^n(t)\le Ca_n^{{-3\over 2}+\gamma+\gamma(2\beta\wedge \xi)}\int _0^t\int 1(|u(s,x)|\le a_n, u'(s,x)\approx a_n^\beta)\Psi_s(x)dxds.
\end{equation}

For $\beta=\bar\beta$ the precise meaning of $u'(s,x)\approx a_n^{\bar\beta}$ is $0\le u'(s,x)\le a_n^{\bar\beta}$ and we have from \eqref{Inbetabnd},
\begin{equation}\label{In1bnd}
I^n_{\bar\beta}(t)\le C ta_n^{{-3\over 2}+\gamma+\gamma(2\bar\beta\wedge\xi)}.
\end{equation}
Consider $0<\beta<\bar\beta$.  
Recall that $\{x:\Psi_s(x)>0\hbox{ for some }s\le t_0\}\subset[-K_1,K_1]$, let
\[S_n(s)=\{x\in[-K_1,K_1]:|u(s,x)|\le a_n, u'(s,x)\ge a_n^\beta\}\]
and $|S_n(s)|$ denote the Lebesgue measure of $S_n(s)$.    From \eqref{uC2-onZ2} we see that if $x\in S_n(s)$, then $u'(s,y)\ge {a_n^\beta\over 2}$ if $|y-x|\le L^{-1}a_n^{\beta/\xi}$, and so by the Fundamental Theorem of Calculus,
\[u(s,y)>a_n\hbox{ if }4a_n^{1-\beta}< |y-x|\le L^{-1}a_n^{\beta/\xi}.\]
A simple covering argument now shows that $|S_n(s)|\le c(L,K_1)a_n^{1-\beta}a_n^{-\beta/\xi}$ and  \eqref{Inbetabnd} implies 
\begin{align}\nonumber
I_\beta^n(t)&\le Cta_n^{{-3\over 2}+\gamma+\gamma(2\beta\wedge\xi)+1-\beta-{\beta\over \xi}}\\
\label{In2bnd}&\le Ct a_n^{\gamma(1+(2\bar\beta)\wedge\xi)-{1\over 2}-\bar\beta(1+{1\over\xi})}.
\end{align}
So from \eqref{In1bnd} and \eqref{In2bnd} we see that $\lim_{n\to\infty}I^n_\beta(t)=0$ will follow for all $\beta\le \bar\beta$ if 
\begin{equation*}
\gamma(1+(2\bar\beta\wedge 1))>{3\over 2}\quad\hbox{ and }\quad\gamma(1+(2\bar\beta\wedge 1))>{1\over 2}+2\bar\beta,
\end{equation*}
that is, $\gamma>(1+(2\bar\beta\wedge 1))^{-1}({3\over 2}\vee({1\over 2}+2\bar\beta))$.  The right-hand side is minimized when $\bar\beta={1\over 2}$, and leads to $\gamma>{3\over 4}$, as required, and also establishes the range $0\le \beta\le {1\over 2}$, which will be used below.

The above heuristics show that $\gamma>3/4$ and the regularity of $u$ given in \eqref{uC2-onZ2} (or \eqref{uC2onZ2}) is optimal for our approach.  If we try weakening the regularity condition on $u$, the above discussion shows we would have to increase $3/4$ to show that $I^n(t)\to 0$.  The earlier discussion shows that a strengthening of the regularity on $u$ would require increasing $3/4$ as well.  

A major obstruction to \eqref{uC2-onZ2} is the fact that we cannot expect $u'(s,x)$ to exist as soon as $u(s,x)\neq 0$ (and don't even know this is the case for $u(s,x)=0$).  So instead, if $D(r,y)=\sigma(r,y,X^1(r,y))-\sigma(r,y,X^2(r,y))$, then we will use \eqref{wspde} to decompose $u$ as 
\begin{align}\label {uidecomp}
u(t,x)&=\int_0^{t-a_n}\int p_{t-r}(y-x)D(r,y)W(dr,dy)+\int_{t-a_n}^t \int p_{t-r}(y-x)D(r,y)W(dr,dy)\\
\nonumber&\equiv u_{1,a_n}(t,x)+u_{2,a_n}(t,x).
\end{align}
$u_{1,a_n}$ is smooth in the spatial variable and so the above arguments may be applied with $u'_{1,a_n}(t,x)$ playing the role of $u'(t,x)$, while $u_{2,a_n}$ and its increments should lead to small and manageable error terms.  
Proposition~\ref{u2modulus} gives the required bounds on the increments of $u_{2,a_n^\alpha}$, and (as noted above) Corollary~\ref{u1'modulus} is the analogue of \eqref{uC2-onZ2} for $u'_{1,a_n^\alpha}$ ($\alpha\in[0,1]$).  
(The reason for the extension to $a_n^\alpha$ is discussed below.)
The proofs of these results are incorporated into an inductive proof of a space-time bound $(P_m)$ for $u(t,x)$ when $(t,x)$ is close to a point $(\hat t_0,\hat x_0)$ where 
\begin{equation}\label{hatcond}|u(\hat t_0,\hat x_0)|\le a_n\hbox{ and }|u'_{1,a_n^\alpha}(\hat t_0,\hat x_0)|\le a_n^\beta.
\end{equation}
If 
\begin{equation}\label{ddef} d=\sqrt{|t-\hat t_0|}+|x-\hat x_0|, \end{equation} 
then, roughly speaking, $(P_m)$ bounds $|u(t,x)|$ by 
\begin{equation}\label{Pmintu}d^\xi[d^{\tilde \gamma_m-1}+a_n^\beta],
\end{equation}
where $\tilde \gamma_m$ increases in $m$ and equals $2$ for $m$ large, and $\xi<1$ as usual.  When $\tilde\gamma_m=2$ this   
does capture the kind of bound one expects from \eqref{uC2-onZ2}.   The reader may find a precise statement of $(P_m)$ prior to Proposition~\ref{Pminduction} (the statement of its validity).

The $m=0$ case will be an easy consequence of our improved local modulus of continuity, Theorem~\ref{theorem:lipmod}.  Note that  \eqref{holder'}) implies
\begin{equation}\label{Dlipbound}|D(r,y)|\le R_0e^{R_1|y|}|u(r,y)|^\gamma.\end{equation}
The inductive proof of $(P_m)$ proceeds by using \eqref{Dlipbound} and then \eqref{Pmintu} to bound the  square functions associated with the space-time increments of $u'_{1,a_n^\alpha}$ and $u_{2,a_n^\alpha}$ for points near $(\hat t_0,\hat x_0)$ as in \eqref{hatcond}(recall \eqref{uidecomp}).  These give good control of the integrands of these square functions near the points where they have singularities.  
This will then lead to Corollary~\ref{u1'modulus} and Proposition~\ref{u2modulus}, our ``$m$th order" bounds for the increments of $u'_{1,a_n^\alpha}$ and $u_{2,a_n^\alpha}$.  We then use the slightly generalized version of \eqref{uidecomp},
\begin{equation}\label{uidecompa}u(t,x)=u_{1,a_n^\alpha}(t,x)+u_{2,a_n^\alpha}(t,x)
\end{equation}
to derive $(P_{m+1})$.  At this point we will optimize over $\alpha$ since decreasing $\alpha$ increases the regularity of $u_{1,a_n^\alpha}$ but increases the size of the error term $u_{2,a_n^\alpha}$.   The optimal choice will be so that $a_n^\alpha\approx d$, where $d$ is as in \eqref{ddef}.  

There are at least two issues to address here.  First, how do you control $u'_{1,a_n^\alpha}(t,x)$ when all you know is $|u'_{1,a_n}(t,x)|\le a_n^\beta$? Second, how do you control the {\it{time}} increments of $u_{1,a_n^\alpha}$ when you only have good estimates on the spatial derivatives?  The first question is answered in Proposition~\ref{Fmodulus2} which will give surprisingly good bounds on $|u'_{1,a_n^\alpha}(\hat t_0,\hat x_0)-u'_{1,a_n}(\hat t_0,\hat x_0)|$.  The second question is answered in Proposition~\ref{u1modulus}, where the key step is to note (see \eqref{u1andecomp}) that for $t>t'$,
\[ |u_{1,a_n^\alpha}(t,x)-u_{1,a_n^\alpha}(t',x)|\approx |P_{t-t'}(u_{1,a_n^\alpha}(t',\cdot))(x) -u_{1,a_n^\alpha}(t',x)|,\]
where $P_t$ is the Brownian semigroup.  The fact that the Brownian semigroup, $P_tf$, inherits temporal regularity from spatial regularity of $f$ will give the required regularity in time. 

A critical step in the above argument was finding a form of $(P_m)$ which actually iterates to produce $(P_{m+1})$.  Note also that although the required bound on $I^n_\beta(t)$ (see \eqref{Inbetabnd}) only required good spatial estimates for $u(s,\cdot)$  near points $(s,x)=(\hat t_0,\hat x_0)$ as in \eqref{hatcond}, the iteration of estimates requires an expansion in both space and time.

Turning now to a brief description of the contents of the paper, we first set $b\equiv0$.  In Section~\ref{sec:verifiction_Hyp}  we reduce $(H_1)$ and $(H_2)$ to a result (Proposition~\ref{tildeJ}) on control of the spatial increments of $u_{2,a_n^\alpha}$ and size of $u'_{1,a_n^\alpha}$ on relatively long intervals near a spatial point where $|u(s,\hat x_0)|$ is small and $u'_{1,a_n^\alpha}(s,\hat x_0)\approx a_n^\beta$.    This includes the covering argument sketched above.  Section~\ref{sec:integralbounds} gives some integral bounds for heat kernels and their derivatives which will help bound the square functions of the increments of $u'_{1,a_n^\alpha}$ and $u_{2,a_n^\alpha}$.   The heart of the proof of Theorem~\ref{theorem:unique} is given in Section~\ref{sec4} where the inductive proof of $(P_m)$ is given.  As was sketched above, this argument includes good local expansions for $u'_{1,a_n^\alpha}$ and $u_{2,a_n^\alpha}$ near points where $|u|$ and $|u'_{1,a_n^\alpha}|$ are small, although the (easier) proof for  $u_{2,a_n^\alpha}$ is deferred until Section~\ref{sec6}.  These expansions, with $m$ large enough, are then used in Section~\ref{sec5} to prove Proposition~\ref{tildeJ} and so complete the proof of Theorem~\ref{theorem:unique} for $b\equiv 0$.   In Section~\ref{secdrifts} we describe the relatively simple additions that are needed to include a Lipschitz drift $b$ in the argument already presented.   

\section{Verification of the Hypotheses of Proposition~\ref{prop:Inbound}}\label{sec:verifiction_Hyp}
\setcounter{equation}{0}
\setcounter{theorem}{0}
We assume throughout this Section that $b\equiv0$--the relatively simple refinements required to include the drift are outlined in Section~\ref{secdrifts}.  Let $X^1,X^2$ be as in Section~\ref{secmainres}, $u=X^1-X^2$, and assume the hypotheses of Theorem~\ref{theorem:unique} as well as \eqref{holder'}.  If $D(s,y)=\sigma(s,y,X^1(s,y))-\sigma(s,y,X^2(s,y))$, then 
\begin{equation}\label{ueq}
u(t,x)=\int_0^t\int p_{t-s}(y-x)D(s,y)\,W(ds,dy)\hbox{ a.s.  for all }(t,x),
\end{equation}
and by \eqref{holder'},
\begin{equation}\label{Dbound}
|D(s,y)|\le R_0e^{R_1|y|}|u(s,y)|^\gamma.
\end{equation}
$\delta$ will always take values in $(0,1]$.  Let 
\begin{equation}u_{1,\delta}(t,x)=P_\delta(u_{(t-\delta)^+})(x)\hbox{ and }u_{2,\delta}(t,x)=u(t,x)-u_{1,\delta}(t,x).\end{equation}
Since $P_\delta:C_{tem}\to C_{tem}$ is uniformly continuous (by Lemma~6.2(ii) of \cite{bib:shiga94}), $u_{1,\delta}$ and $u_{2,\delta}$ both have sample paths in $C(\Re,C_{tem})$.  \eqref{ueq} implies that
\[u_{1,\delta}(t,x)=\int\Bigl[\int_0^{(t-\delta)^+} \int p_{(t-\delta)^+-s}(y-z)D(s,y) W(ds,dy)\Bigr] p_\delta(z-x)dz.\]
A stochastic Fubini argument (Theorem~2.6 of \cite{bib:wal86}) then gives
\begin{equation}\label{u1eq}
u_{1,\delta}(t,x)=\int_0^{(t-\delta)^+} \int p_{t-s}(y-x)D(s,y)W(ds,dy)\hbox{ a.s. for all }(t,x)\in\Rp\times\Re.
\end{equation}
(The above identity is trivial for $t\le \delta$ since $u(0,\cdot)\equiv 0$.)
The expectation condition in Walsh's Theorem~2.6 may be realized by localization with the stopping times $\{T_K\}$, working with $D(s\wedge T_K)$, and letting $K\to\infty$.  It follows that
\begin{equation}\label{u2eq}
u_{2,\delta}(t,x)=\int_{(t-\delta)^+}^t \int p_{t-s}(y-x)D(s,y)W(ds,dy)\hbox{ a.s. for all }(t,x)\in\Rp\times\Re.
\end{equation}
Hence $u_{j,\delta}$ ($j=1,2$) define jointly continuous versions of the right-hand sides of \eqref{u1eq} and \eqref{u2eq}.

\medskip

\noindent {\bf Notation.} If $s,t\ge0$ and $x\in\Re$, let $G_\delta(s,t,x)=P_{(t-s)^++\delta}(u_{(s-\delta)^+})(x)$ and \hfil\break $F_\delta(s,t,x)=-{d\over dx}G_\delta(s,t,x)\equiv -G'_\delta(s,t,x)$, if the derivative exists.

\begin{lemma}\label{G'} $G'_\delta(s,t,x)$ exists for all $(s,t,x)\in\Rp^2\times\Re$, is jointly continuous in $(s,t,x)$, and satisfies
\begin{equation}\label{G'form}
F_\delta(s,t,x)=
\int_0^{(s-\delta)^+}\int p'_{(t\vee s)-r}(y-x)
D(r,y)W(dr,dy)\hbox{ for all }s\hbox{ a.s. for all } (t,x).
\end{equation}
\end{lemma}
\paragraph{Proof.} Since $G_\delta(s,t,x)=\int p_{(t-s)^++\delta}(y-x)u((s-\delta)^+,y)dy$ and 
\begin{equation}\label{dctbound}
\sup_{s\le T, y} 
e^{-|y|}|u((s-\delta)^+,y)|<\infty\hbox{ for all }T>0\hbox{ a.s.,}
\end{equation}
a simple Dominated Convergence argument shows that 
\begin{equation}\label{G'formb}
G_\delta'(s,t,x)=-\int p'_{(t-s)^++\delta}(y-x)u((s-\delta)^+,y)dy\quad\hbox{for all }(s,t,x)\ a.s.
\end{equation}
Another application of \eqref{dctbound} and Dominated Convergence gives the a.s. joint continuity of the right-hand side of \eqref{G'formb}, and hence of $G_\delta'$.  

To prove \eqref{G'form} we may assume without loss of generality that $t\ge s>\delta$.  From \eqref{G'formb} and \eqref{ueq} we have w.p. 1,
\begin{equation*}
G'_\delta(s,t,x)=-\int p'_{t-s+\delta}(y-x)\Bigl[\int_0^{s-\delta}\int p_{s-\delta-r}(z-y)D(r,z)W(dr,dz)\Bigr]dy\ a.s.
\end{equation*}
Now use the stochastic Fubini theorem, as in the derivation of \eqref{u1eq} above, to see that
\begin{align*}
G'_\delta(s,t,x)&=-\int_0^{s-\delta}\int\Bigl[\int p'_{t-s+\delta}(y-x)p_{s-\delta-r}(z-y)dy\Bigr]D(r,z)W(dr,dz)\\
&=-\int_0^{s-\delta}\int p'_{t-r}(z-x)D(r,z)W(dr,dz) \quad a.s.
\end{align*}
In the last line we have used Dominated Convergence yet again to differentiate through the integral in the Chapman-Kolmorgorov equation.  As both sides of \eqref{G'form} are continuous in $s$ we may take the null set to be independent of $s$.
\gdm

\begin{remark} Since $G_\delta(t,t,x)=u_{1,\delta}(t,x)$, as a special case of the above we see that $u'_{1,\delta}(t,x)$ is a.s. jointly continuous and satisfies
\begin{equation}\label{u1'eq}
u'_{1,\delta}(t,x)=-\int_0^{(t-\delta)^+}\int p'_{t-s}(y-x)D(s,y) W(ds,dy)\quad \hbox{a.s. for all }(t,x).
\end{equation}
\end{remark}

\noindent{\bf Definition.} For $(t,x)\in\Rp\times\Re$, 
\begin{align*}\hat x_n(t,x)(\omega)&=\inf\{y\in[x-\sqrt{a_n},x+\sqrt{a_n}]: |u(t,y)|=\inf\{|u(t,z)|:|z-x|\le\sqrt{a_n}\}\}\\
&\in[x-\sqrt{a_n},x+\sqrt{a_n}].\end{align*}

\noindent It is easy to use the continuity of $u$ to check that $\hat x_n$ is well-defined and $\cB(\Rp\times\Re)\times \cF$-measurable.  

We fix a $K_0\in\N^{\ge K_1}$ and positive constants satisfying
\begin{equation}\label{veconditions}
0<\ve_1<{1\over 100}(\gamma-{3\over 4}),\ 0<\ve_0<{\ve_1\over 100}.
\end{equation}
We introduce a grid of $\beta$ values by setting
\[L=L(\ve_0,\ve_1)=\lfloor((1/ 2)-6\ve_1)/\ve_0\rfloor,\]
and
\begin{eqnarray}\label{betai}
\beta_i&=&i\ve_0\in [0,{1\over 2}-6\ve_1],\ \ \alpha_i=2(\beta_i+\ve_1)\in[0,1],\ \ i=0,\dots,L,\\
\nonumber
\beta_{L+1}&=&{1\over 2}-\ve_1\,.
\end{eqnarray}
Note that $\beta=\beta_i\,,i=0,\ldots,L+1,$ satisfies
\begin{equation}\label{betarange}
0\le \beta\le {1\over 2}-\ve_1
\end{equation}

If $s\ge 0$ set
\[J_{n,0}(s)=\Bigl\{x: |x|\le K_0, |\langle u_s,\Phi_x^{m_{n+1}}\rangle|\le a_n, u'_{1,a_n}(s,\hat x_n(s,x))\ge {a_n^{\ve_0} \over 4}\Bigr\},\]
\[J_{n,L}(s)=\Bigl\{x: |x|\le K_0, |\langle u_s,\Phi_x^{m_{n+1}}\rangle|\le a_n, u'_{1,a_n}(s,\hat x_n(s,x))\in [0,{a_n^{\beta_L} \over 4}]\Bigr\},\]
and for $i=1,\dots, L-1$ set
\[J_{n,i}(s)=\Bigl\{x: |x|\le K_0, |\langle u_s,\Phi_x^{m_{n+1}}\rangle|\le a_n, u'_{1,a_n}(s,\hat x_n(s,x))\in[{a_n^{\beta_{i+1}}\over 4},{a_n^{\beta_i}\over 4}]\Bigr\}.\]
If $t_0>0$ is as in 
$(H_2)$ and $i=0,\dots,L$, define
\[J_{n,i}=\Bigl\{(s,x):0\le s,\ x\in J_{n,i}(s)\Bigr\},\]
and if $0\le t\le t_0$, let
\[I^n_i(t)=a_n^{-{3\over 2}-{2\over n}}\int_0^t\int\int 1_{J_{n,i}(s)}(x)|u(s,y)|^{2\gamma}\Phi_x^{m_{n+1}}(y)\Psi_s(x)dydxds.\]
Let 
\begin{align*}
&I^n_+(t)
=a_n^{-{3\over 2}-{2\over n}}\int_0^t\int\int 1(u'_{1,a_n}(s,\hat x_n(s,x))\ge 0)1(|\langle u_s,\Phi_x^{m_{n+1}}\rangle|<a_n)\\
&\phantom{I^N_+(t)=a_n^{-{3\over 2}-{2\over n}}\int_0^t\int\int 1(u'_{1,a_n}}\times|u(s,y)|^{2\gamma}\Phi_x^{m_{n+1}}(y)\Psi_s(x)dydxds.
\end{align*}
Then to prove Proposition~\ref{prop:Inbound} it suffices to construct the stopping times 
$\{U_{M,n}\equiv U_{M,n,K_0}:M,n\in\N\}$ satisfying $(H_1)$ such that 
\begin{eqnarray}
\label{Inplus}
\hbox{for each }M\in\N, \ \lim_{n\to\infty}E(I_+^n(t_0\wedge U_{M,n}))=0.
\end{eqnarray}
Note that~(\ref{Inplus}) implies $(H_2)$ by symmetry (interchange $X_1$ and $X_2$). 

Our definitions imply
\[I_+^n(t)\le \sum_{i=0}^LI^n_i(t)\hbox{ for all }t\le t_0,\]
and so to prove 
$(H_2)$ it suffices to show that for $i=0,\dots,L$,
\begin{equation*}
(H_{2,i})\phantom{AAAAAAAAA}\hbox{ for all }M\in\N,\quad \lim_{n\to\infty}E(I^n_i(t_0\wedge U_{M,n}))=0.\phantom{AAAAAAAAAAAAAAAAAAA}
\end{equation*}

\noindent{\bf Notation.} $\oln(\beta)=a_n^{\beta+5\ve_1}$. 

Now introduce the related sets: 
\begin{align*}
\tilde J_{n,0}(s)=&\Bigl\{x\in[-K_0,K_0]: |\langle u_s,\Phi_x^{m_{n+1}}\rangle|\le a_n, u'_{1,a_n^{\alpha_0}}(s,x')\ge a_n^{\beta_1}/16\\
&\phantom{\Bigl\{x\in[-K_0,K_0]: |\langle u_s,\Phi_x^{m_{n+1}}\rangle|\le a_n, u'}\hbox{ for all }
x'\in[x-5\oln(\beta_0),x+5\oln(\beta_0)],\\
& |u_{2,a_n^{\alpha_0}}(s,x')-u_{2,a_n^{\alpha_0}}(s,x'')|\le 2^{-75}a_n^{\beta_{1}}(|x'-x''|\vee a_n^{\gamma-2\beta_0(1-\gamma)-\ve_1})\\
&\hbox{ for all }x'\in[x-4\sqrt{a_n},x+4\sqrt{a_n}],x''\in[x'-\oln(\beta_0),x'+\oln(\beta_0)],\\
&\hbox{ and }|u(s,x')|\le 3a_n^{(1-\ve_0)/2}\hbox{ for all }x'\in[x-\sqrt{a_n},x+\sqrt{a_n}]\Bigr\},\\
\tilde J_{n,L}(s)=&\Bigl\{x\in[-K_0,K_0]:  |\langle u_s,\Phi_x^{m_{n+1}}\rangle|\le a_n, |u'_{1,a_n^{\alpha_L}}(s,x')|\le a_n^{\beta_L}\\
&\phantom{\Bigl\{x\in[-K_0,K_0]:  |\langle u_s,\Phi_x^{m_{n+1}}\rangle|\le a_n, |u'}\hbox{ for all }
x'\in[x-5\oln(\beta_L),x+5\oln(\beta_L)],\\
&\hbox{ and }
 |u_{2,a_n^{\alpha_L}}(s,x')-u_{2,a_n^{\alpha_L}}(s,x'')|\le 2^{-75}a_n^{\beta_{L+1}}(|x'-x''|\vee a_n^{\gamma-2\beta_L(1-\gamma)-\ve_1})\\
&\hbox{ for all }x'\in[x-4\sqrt{a_n},x+4\sqrt{a_n}],x''\in[x'-\oln(\beta_L),x'+\oln(\beta_L)]\Bigr\},
\end{align*}
and for $i\in\{1,\dots,L-1\}$,
\begin{align*}
\tilde J_{n,i}(s)=&\Bigl\{x\in[-K_0,K_0]:  |\langle u_s,\Phi_x^{m_{n+1}}\rangle|\le a_n, u'_{1,a_n^{\alpha_i}}(s,x')\in[a_n^{\beta_{i+1}}/16, a_n^{\beta_i}]\\
&\phantom{\Bigl\{x\in[-K_0,K_0]:  |\langle u_s,\Phi_x^{m_{n+1}}\rangle|\le a_n, u'}\hbox{ for all }
x'\in[x-5\oln(\beta_i),x+5\oln(\beta_i)],\\
&\hbox{ and }
 |u_{2,a_n^{\alpha_i}}(s,x')-u_{2,a_n^{\alpha_i}}(s,x'')|\le 2^{-75}a_n^{\beta_{i+1}}(|x'-x''|\vee a_n^{\gamma-2\beta_i(1-\gamma)-\ve_1})\\
&\hbox{ for all }x'\in[x-4\sqrt{a_n},x+4\sqrt{a_n}],x''\in[x'-\oln(\beta_i),x'+\oln(\beta_i)]\Bigr\},
\end{align*}
Finally for $0\le i\le L$, set 
\[\tilde J_{n,i}=\{(s,x): s\ge 0, x\in\tilde J_{n,i}(s)\}.\]
\noindent{\bf Notation.} $n_M(\ve_1)=\inf\{n\in\N:a_n^{\ve_1}\le 2^{-M-4}\}$, $n_0(\ve_0,\ve_1)=\sup\{n\in\N:\sqrt{a_n}<2^{-a_n^{-\ve_0\ve_1/4}}\}$, where $\sup\emptyset=1$.  

The following proposition  will be proved in Section~\ref{sec5}.

\begin{proposition}\label{tildeJ} $\tilde J_{n,i}(s)$ is a compact set for all $s\ge 0$.
There exist  stopping times\\
$\left\{U_{M,n}\equiv U_{M,n,K_0}:M,n\in\N\right\}$ satisfying $(H_1)$ from Proposition~\ref{prop:Inbound} such that
for $i\in\{0,\dots,L\}$,  $\tilde J_{n,i}(s)$ contains $J_{n,i}(s)$ for all $0\le s<U_{M,n}$ and 
\begin{equation}\label{ncond}
n>n_M(\ve_1)\vee n_0(\ve_0,\ve_1).
\end{equation} 
\end{proposition}

Since $K_0\in \N^{\geq K_1}$ was arbitrary this 
proposition  implies that there exist stopping times 
satisfying $(H_1)$ such that the inclusion $\tilde J_{n,i}(s)\supset J_{n,i}(s)$ holds up to these stopping times for $n$ sufficiently large.  
This inclusion means that given a value of the derivative of $u_{1,a_n}$  at {\it some} point in a {\it small} neighborhood $[x-\sqrt{a_n}, x+\sqrt{a_n}]$ of $x$ where $|u|$ is small, one can guarantee that the derivative of $u_{1,a_n^{\alpha}}$ (for a certain $\alpha$) is of the same order at {\it any} 
point in a much {\it larger} neighborhood of $x$ (note that 
$ \oln(\beta)\gg \sqrt{a_n}$ for $\beta \leq {1\over 2}- 6\ve_1$). Moreover we can also 
control $u_{2,a_n^{\alpha}}$ on those long intervals. 
Our goal now is to show that this implies  $(H_{2,i})$ for $i=0,\ldots, L$. The next three   lemmas provide necessary tools for this. 

Throughout the rest of the section we may, and shall, assume that the parameters $M,n\in\N$ satisfy~(\ref{ncond}), although the importance of $n_0$ in (\ref{ncond}) will not be clear until Section~\ref{sec5}.

 \medskip

\begin{lemma}\label{ucontrol} Assume $i\in\{0,\dots,L\}$, $x\in\tilde J_{n,i}(s)$ and $|x'-x|\le 4\sqrt{a_n}$.

\noindent(a) If $i>0$, then $|u(s,x'')-u(s,x')|\le 2a_n^{\beta_i}(|x''-x'|\vee a_n^{\gamma-2\beta_i(1-\gamma)-\ve_1})$ for all $|x''-x'|\le \oln(\beta_i)$.

\noindent(b) If $i<L$, and $a_n^{\gamma-2\beta_i(1-\gamma)-\ve_1}\le |x''-x'|\le \oln(\beta_i)$, then
\begin{equation*}u(s,x'')-u(s,x')\begin{cases}\ge 2^{-5}a_n^{\beta_{i+1}}(x''-x')& \text{ if }x''\ge x',\\
\le 2^{-5}a_n^{\beta_{i+1}}(x''-x')&\text{ if }x''\le x'.
\end{cases}
\end{equation*}
\end{lemma}
\paragraph{Proof.}  (a) For $n,i,s,x,x',x''$ as in (a), we have (since $\beta_i+5\ve_1<{1\over 2}$)
\begin{equation}\label{incrbound}
|x'-x|\vee|x''-x|\le 5\oln(\beta_i).
\end{equation}
 We can therefore apply the definition of $\tilde J_{n,i}$ and the Mean Value Theorem to conclude that
\begin{align*}
|u(s,x'')-u(s,x')|&\le |u_{1,a_n^{\alpha_i}}(s,x'')-u_{1,a_n^{\alpha_i}}(s,x')|+ |u_{2,a_n^{\alpha_i}}(s,x'')-u_{2,a_n^{\alpha_i}}(s,x')|\\
&\le a_n^{\beta_i}|x''-x'|+2^{-75}a_n^{\beta_{i+1}}(|x''-x'|\vee a_n^{\gamma-2\beta_i(1-\gamma)-\ve_1})\\
&\le 2a_n^{\beta_i}(|x''-x'|\vee a_n^{\gamma-2\beta_i(1-\gamma)-\ve_1}).
\end{align*}

\noindent(b) Consider $n,i,s,x,x',x''$ as in (b) with $x''\ge x'$.  We again have \eqref{incrbound} and can argue as in (a) to see that
\begin{align*}
u(s,x'')-u(s,x')&= u_{1,a_n^{\alpha_i}}(s,x'')-u_{1,a_n^{\alpha_i}}(s,x')+ u_{2,a_n^{\alpha_i}}(s,x'')-u_{2,a_n^{\alpha_i}}(s,x')\\
&\ge (a_n^{\beta_{i+1}}/16)(x''-x')-2^{-75}a_n^{\beta_{i+1}}(x''-x')\\
&\ge (a_n^{\beta_{i+1}}/32)(x''-x').
\end{align*}
The case $x''\le x'$ is similar.
\gdm

\noindent{\bf Notation.} $l_n(\beta_i)=(65a_n^{1-\beta_{i+1}})\vee a_n^{\gamma-2\beta_i(1-\gamma)-\ve_1}$, $F_n(s,x)=\int\Phi^{m_{n+1}}(y-x)u(s,y)dy$.

\begin{lemma}\label{Fnincrement} Assume $i\in\{0,\dots,L-1\}$ and $(s,x)\in \tilde J_{n,i}$.

\noindent(a) If $l_n(\beta_i)\le |x-\tilde x|\le \oln(\beta_i)$, then 
\begin{equation*}
F_n(s,\tilde x)-F_n(s,x)\begin{cases}\ge 2^{-5}a_n^{\beta_{i+1}}(\tilde x-x)&\text{ if } \tilde x\ge x,\\
\le 2^{-5}a_n^{\beta_{i+1}}(\tilde x-x)&\text{ if } \tilde x\le x.
\end{cases}
\end{equation*}

\noindent(b) $[x-\oln(\beta_i),x-l_n(\beta_i)]\cup[x+l_n(\beta_i),x+\oln(\beta_i)]\subset \tilde J_{n,i}(s)^c$. 
\end{lemma}
\paragraph{Proof.} (a) Assume $\tilde x\in[x+l_n(\beta_i),x+\oln(\beta_i)]$.  Then
\begin{equation}\label{Fexpression}
F_n(s,\tilde x)-F_n(s,x)=\int_{-\sqrt{a_n}}^{\sqrt{a_n}}\Phi^{m_{n+1}}(z)(u(s,\tilde x+z)-u(s,x+z))dz.
\end{equation}
For $|z|\le \sqrt{a_n}$, let $x''=\tilde x+z$ and $x'=x+z$.  Then $|x'-x|\le \sqrt{a_n}$ and 
\[x''-x'=\tilde x-x\in[l_n(\beta_i),\oln(\beta_i)].\]
Therefore Lemma~\ref{ucontrol}(b) and \eqref{Fexpression} imply 
\begin{align*}
F_n(s,\tilde x)-F_n(s,x)&\ge \int _{-\sqrt{a_n}}^{\sqrt{a_n}}\Phi^{m_{n+1}}(z)2^{-5}a_n^{\beta_{i+1}}(\tilde x-x)dz\\
&=2^{-5}a_n^{\beta_{i+1}}(\tilde x-x).
\end{align*}
The proof for $\tilde x<x$ is similar.
\medskip

\noindent(b) If $\tilde x\in [x-\oln(\beta_i),x-l_n(\beta_i)]\cup[x+l_n(\beta_i),x+\oln(\beta_i)]$, then
\begin{align*}
|F_n(s,\tilde x)|&\ge |F_n(s,\tilde x)-F_n(s,x)|-|F_n(s,x)|\\
&\ge 2^{-5}a_n^{\beta_{i+1}}l_n(\beta_i)-a_n\quad\text{(by (a) and $(s,x)\in\tilde J_{n,i}$)}\\
&\ge {33\over 32}a_n.
\end{align*}
Therefore $\tilde x\notin \tilde J_{n,i}(s)$.\gdm

To ensure that (b) is not vacuous we obtain some crude lower bounds on the interval given there.

\begin{lemma}\label{intervalbounds}
If $i\in\{0,\dots,L\}$, then 
\[l_n(\beta_i)<\sqrt{a_n}<{1\over 2}\oln(\beta_i).\]
\end{lemma}
\paragraph{Proof.} 
\begin{align*}
l_n(\beta_i)a_n^{-1/2}&=(65a_n^{{1\over 2}-\beta_{i+1}})\vee a_n^{\gamma-{1\over 2}-2\beta_i(1-\gamma)-\ve_1}\\
&\le (65a_n^{5\ve_1})\vee a_n^{2\gamma-{3\over 2}-\ve_1}\quad\hbox{(since $\beta_i<1/2$)}\\
&<1
\end{align*}
by \eqref{veconditions} and because $a_n^{5\ve_1}<2^{-20}$ by \eqref{ncond}. This gives the first inequality. For the second one, use $\beta_i\le {1\over 2}-6\ve_1$  and \eqref{ncond} to see that
\begin{equation*}\sqrt{a_n}\oln(\beta_i)^{-1}=a_n^{{1\over 2}-\beta_i-5\ve_1}\le a_n^{\ve_1}<1/2.
\end{equation*}
\gdm 

Let $|A|$ denote the Lebesgue measure of $A\subset\Re$.

\begin{lemma}\label{LebJ} For all $i\in\{0,\dots,L-1\}$ and $s\ge 0$,
\[|\tilde J_{n,i}(s)|\le 10K_0\oln(\beta_i)^{-1}l_n(\beta_i).\]
\end{lemma}
\paragraph{Proof.} Fix $s$, $i$ as above.  Let $\cI_{n,i}(x)=(x-l_n(\beta_i),x+l_n(\beta_i))\subset \tcI_{n,i}(x)=(x-\oln(\beta_i),x+\oln(\beta_i))$. This inclusion follows from Lemma~\ref{intervalbounds}.  The compactness of $\tilde J_{n,i}(s)$ (Proposition~\ref{tildeJ}) implies there are $x_1,\dots,x_Q\in \tilde J_{n,i}(s)$ so that $\tilde J_{n,i}(s)\subset \cup_{j=1}^Q \cI_{n,i}(x_j)$.  

Assume that for some $k\neq j$, $|x_k-x_j|\le \oln(\beta_i)/2$.  We claim that $\cI_{n,i}(x_j)\subset \tcI_{n,i}(x_k)$.  Indeed, if $y\in \cI_{n,i}(x_j)$ then
\[|y-x_k|\le |y-x_j|+|x_j-x_k|< l_n(\beta_i)+\oln(\beta_i)/2< \oln(\beta_i),\]
the last by Lemma~\ref{intervalbounds}, and the claim is proved.  Lemma~\ref{Fnincrement}(b) implies
 \[\tilde J_{n,i}(s)\cap (\tcI_{n,i}(x_k)-\cI_{n,i}(x_k))=\emptyset,\] 
 and so the above claim gives
\[\cI_{n,i}(x_j)\cap\tilde J_{n,i}(s)\subset \tcI_{n,i}(x_k)\cap \tilde J_{n,i}(s)=\cI_{n,i}(x_k)\cap\tilde J_{n,i}(s).\]
Therefore we may omit $\cI_{n,i}(x_j)$ and still have a cover of $\tilde J_{n,i}(s)$.  Doing this sequentially for $x_1,\dots,x_Q$, we may therefore assume that 
\begin{equation*}
|x_k-x_j|>\oln(\beta_i)/2\quad\hbox{ for all }k\neq j.
\end{equation*}
Since each $x_j\in \tilde J_{n,i}(s)\subset [-K_0,K_0]$, this implies $Q\le 2K_0(\oln(\beta_i)/2)^{-1}+1$, and therefore
\[|\tilde J_{n,i}(s)|\le (4K_0\oln(\beta_i)^{-1}+1)2l_n(\beta_i)\le 10K_0\oln(\beta_i)^{-1}l_n(\beta_i).\]\gdm

\medskip

\noindent{\bf Proof of $(H_2)$.} 
Fix $M\in\N$.  Recall from our discussion after the definition of 
$J_{n,i}$ sets that it suffices to show that for all $i\in\{0,\dots,L\}$,
\begin{equation*}
(H_{2,i})\phantom{AAAAAAAAAAAAAAA} \lim_{n\to\infty}E(I^n_i(t_0\wedge U_{M,n}))=0.\phantom{AAAAAAAAAAAAAAAAAAAAAAAAAAAA}
\end{equation*}
We will in fact show that if we strengthen \eqref{ncond} to
\begin{equation}\label{ncond2}
n>n_M(\ve_1)\vee n_0(\ve_0,\ve_1)\vee{2\over \ve_1},
\end{equation}
then
we have the stronger $\Le^\infty$ bound
\begin{equation}\label{Linfty}
I^n_i(t_0\wedge U_{M,n})\le c_1(\Psi)t_0K_0a_n^{\gamma-{3\over 4}},
\end{equation}
which clearly implies $(H_{2,i})$ since $\gamma>{3\over 4}$. Proposition~\ref{tildeJ}, $Supp(\Phi_x^{m_{n+1}})\subset[x-\sqrt{a_n},x+\sqrt{a_n}]$ and ${2\over n}<\ve_1$ (by \eqref{ncond2}) imply
\begin{align}\label{Inibound}
I^n_i(t_0\wedge U_{M,n})\le a_n^{-{3\over 2}-{\ve_1}}\int_0^{t_0}\int\int &1(s<U_{M,n})1_{\tilde J_{n,i}(s)}(x)|u(s,y)|^{2\gamma}\\
\nn&\phantom{}\times1(|y-x|\le \sqrt{a_n})\Phi_x^{m_{n+1}}(y)\Psi_s(x)dydxds.
\end{align}

Consider first \eqref{Linfty} for $i=0$.  For $x\in\tilde J_{n,0}(s)$ and $|y-x|\le \sqrt{a_n}$, we have $|u(s,y)|\le 3a_n^{(1-\ve_0)/2}$ and so from \eqref{Inibound},
\begin{align}
\nn I^n_{0}(t_0\wedge U_{M,n})&\le a_n^{-{3\over 2}-{\ve_1}}3^{2\gamma}a_n^{\gamma(1-\ve_0)}\Vert\Psi\Vert_\infty\int_0^{t_0}|\tilde J_{n,0}(s)|ds\\
\nn&\le a_n^{-{3\over 2}-{\ve_1}}3^{2\gamma}a_n^{\gamma(1-\ve_0)}\Vert\Psi\Vert_\infty t_010K_0a_n^{-5\ve_1}((65a_n^{1-\ve_0})\vee a_n^{\gamma-\ve_1})\quad\hbox{(by Lemma~\ref{LebJ})}\\
\nn&\le c_1(\Psi) t_0K_0 a_n^{2\gamma-{3\over 2}}a_n^{-\gamma\ve_0-7\ve_1}\\
\nn&\le c_1(\Psi)t_0K_0a_n^{\gamma-{3\over 4}},
\end{align}
as required, where \eqref{veconditions} is used in the last two lines.  

Consider now $i\in\{1,\dots,L\}$.  Assume $x\in\tilde J_{n,i}(s)$ and $|y-x|\le \sqrt{a_n}$. 
We have $|\langle u_s,\Phi^{m_{n+1}}_x\rangle|\le a_n$ and \[Supp(\Phi_x^{m_{n+1}})\subset [x-\sqrt{a_n},x+\sqrt{a_n}].\]  
Using the continuity of $u(s,\cdot)$, we conclude that
\[|u(s,\hat x_n(s,x))|\le a_n,\]
 and of course we have $|\hat x_n(s,x)-x|\le \sqrt{a_n}$.  Therefore
\begin{align}
\label{ineq1}|y-\hat x_n(s,x)|&\le 2\sqrt{a_n}\\
\nn&\le  \oln(\beta_i)\quad\hbox{(by Lemma~\ref{intervalbounds}).}\end{align}
Apply Lemma~\ref{ucontrol}(a) with $x''=y$ and $x'=\hat x_n(s,x)$,  to see that
\begin{align}
\nn |u(s,y)|&\le |u(s,\hat x_n(s,x))|+2a_n^{\beta_i}(|y-\hat x_n(s,x)|\vee a_n^{\gamma-2\beta_i(1-\gamma)-\ve_1})\\
\label{localubnd}&\le a_n+4a_n^{\beta_i+{1\over 2}}\le 5a_n^{\beta_i+{1\over 2}},
\end{align}
where \eqref{ineq1} and Lemma~\ref{intervalbounds} are used in the next to last inequality.  Use \eqref{localubnd}
in \eqref{Inibound} and conclude that
\begin{equation}\label{Inibound2}
I^n_i(t_0\wedge U_{M,n})\le a_n^{-{3\over 2}-\ve_1}5^{2\gamma}a_n^{2\gamma(\beta_i+{1\over 2})}\Vert\Psi\Vert_\infty\int_0^{t_0}|\tilde J_{n,i}(s)|ds,\ i=1,\dots,L.
\end{equation}

Assume now $1\le i\le L-1$.  Apply Lemma~\ref{LebJ} to the right-hand side of \eqref{Inibound2} to see that
\begin{align}
\nn I^n_i(t_0\wedge U_{M,n})&\le c_1(\Psi)t_0K_0a_n^{-{3\over 2}-\ve_1+2\gamma(\beta_i+{1\over 2})}(a_n^{1-\beta_{i+1}}\vee a_n^{\gamma-2\beta_i(1-\gamma)-\ve_1})a_n^{-5\ve_1-\beta_i}\\
\label{Inibound3}&=c_1(\Psi)t_0K_0[a_n^{\rho_{1,i}}\vee a_n^{\rho_{2,i}}].
\end{align}
A bit of arithmetic shows that
\begin{align*}
\rho_{1,i}&=\gamma-{1\over 2}-2\beta_i(1-\gamma)-\ve_0-6\ve_1\\
&>2\gamma-{3\over 2}-\ve_0-6\ve_1\quad\hbox{(use $\beta_i<1/2$)}\\
&>\gamma-{3\over 4},\\
\end{align*}
where \eqref{veconditions} is used in the last.  We also have
\begin{align*}
\rho_{2,i}&=-{3\over 2}+2\gamma+\beta_i(4\gamma-3)-7\ve_1\\
&>2\gamma-{3\over 2}-7\ve_1>\gamma-{3\over 4},
\end{align*}
again using \eqref{veconditions} in the last inequality.  Use these bounds on $\rho_{l,i}$, $l=1,2$ in \eqref{Inibound3} to prove \eqref{Linfty} for $i\le i\le L-1$.  

It remains to prove \eqref{Linfty} for $i=L$.  For this, use the trivial bound $|\tilde J_{n,i}(s)|\le 2K_0$ in \eqref{Inibound2} and obtain
\begin{align*}
I^n_L(t_0\wedge U_{M,n})&\le a_n^{-{3\over 2}-\ve_1}5^{2\gamma}a_n^{2\gamma(\beta_L+{1\over 2})}\Vert\Psi\Vert_\infty 2K_0t_0\\
& \le c_1(\Psi) K_0t_0 a_n^{-{3\over 2}-\ve_1+\gamma+2\gamma({1\over 2}-6\ve_1-\ve_0)}\\
&\le c_1(\Psi) K_0t_0 a_n^{2\gamma-{3\over 2}-15\ve_1}\\
&\le c_1(\Psi) K_0t_0 a_n^{\gamma-{3\over 4}},
\end{align*}
yet again using \eqref{veconditions} in the last.  This proves \eqref{Linfty} in the last case of $i=L$.  Having proved \eqref{Linfty} in all cases, we have finished the proof of 
$(H_2)$.\gdm

Proposition~\ref{prop:Inbound} therefore applies and establishes Theorem~\ref{theorem:unique} for $b\equiv0$, except for the proof of Proposition~\ref{tildeJ}.  This will be the objective of the next three sections.

\section{Some Integral Bounds for Heat Kernels}\label{sec:integralbounds}
\setcounter{equation}{0}
\setcounter{theorem}{0}
If $0<p\le 1$, $q\in\Re$ and $0\le \Delta_2\le \Delta_1\le t$, define
\[J_{p,q}(\Delta_1,\Delta_2,\Delta)=\int_{t-\Delta_1}^{t-\Delta_2}(t-s)^q\Bigl(1\wedge{\Delta\over t-s}\Bigr)^p\,ds.
\]
These integrals will arise frequently in our modulus of continuity estimates.
\begin{lemma}\label{Jbnd}
\noindent(a) If $q>p-1$, then
\begin{equation}\label{J1}
J_{p,q}(\Delta_1,\Delta_2,\Delta)\le {2\over q+1-p}(\Delta\wedge \Delta_1)^p\Delta_1^{q+1-p}.
\end{equation}
\noindent(b) If $-1<q<p-1$, then
\begin{align}
\nonumber J&_{p,q}(\Delta_1,\Delta_2,\Delta)\\
&\label{J2}\le((p-1-q)^{-1}+(q+1)^{-1})[(\Delta\wedge \Delta_1)^{q+1}1(\Delta_2\le \Delta)+(\Delta\wedge \Delta_1)^p \Delta_2^{q-p+1}1(\Delta_2>\Delta)]\\
&\label{J3}\le((p-1-q)^{-1}+(q+1)^{-1})\Delta^p(\Delta\vee\Delta_2)^{q-p+1}.
\end{align}
\noindent(c) If $q<-1$, then
\begin{equation}
J_{p,q}(\Delta_1,\Delta_2,\Delta)\le 2|q+1|^{-1}(\Delta\wedge \Delta_2)^p\Delta_2^{q+1-p}.
\end{equation}
\end{lemma}
\paragraph{Proof.} For all $p,q$ as above, 
\begin{align}\label{Jform}\nonumber J_{p,q}&=\int_{\Delta_2}^{\Delta_1}u^q\Bigl(1\wedge{\Delta\over u}\Bigr)^p\,du\\
&=1(\Delta_2< \Delta)\int_{\Delta_2}^{\Delta\wedge\Delta_1}u^q\,du+
1(\Delta_1>\Delta)\int_{\Delta_2\vee \Delta}^{\Delta_1}\Delta^pu^{q-p}\,du.
\end{align}

\noindent(a) From \eqref{Jform},
\begin{align*}
J_{p,q}&\le 1(\Delta_2<\Delta){(\Delta\wedge \Delta_1)^{q+1}\over q+1}+1(\Delta_1>\Delta)\Delta^p{\Delta_1^{q-p+1}\over q-p+1}\\
&\le ((q+1)^{-1}+(q-p+1)^{-1})(\Delta\wedge \Delta_1)^p \Delta_1^{q+1-p},
\end{align*}
which gives the required bound.

\noindent(b) Again \eqref{Jform} implies
\begin{align*}
J_{p,q}&\le 1(\Delta_2<\Delta){(\Delta\wedge \Delta_1)^{q+1}\over q+1}+1(\Delta_1>\Delta)\Delta^p{(\Delta_2\vee \Delta)^{q+1-p}\over p-1-q}\\
&\le 1(\Delta_2\le\Delta)(\Delta\wedge \Delta_1)^{q+1}((q+1)^{-1}+(p-1-q)^{-1})+1(\Delta_2> \Delta)(\Delta\wedge \Delta_1)^p\Delta_2^{q-p+1},
\end{align*}
which gives the first inequality.  The second inequality is elementary.

\noindent(c) By \eqref{Jform},
\begin{align*}
J_{p,q}&\le 1(\Delta_2<\Delta){\Delta_2^{q+1}\over |q+1|}+1(\Delta_1>\Delta)\Delta^p{(\Delta_2\vee\Delta)^{q+1-p}\over p-1-q}\\
&\le {1(\Delta_2<\Delta)\over |q+1|}(\Delta_2\wedge \Delta)^p\Delta_2^{q+1-p}+1(\Delta_2\le \Delta<\Delta_1){\Delta^{q+1}\over p-1-q}\\
&\quad +1(\Delta<\Delta_2){(\Delta\wedge \Delta_2)^p\Delta_2^{q+1-p}\over p-1-q}\\
&\le {2\over |q+1|}(\Delta\wedge \Delta_2)^p\Delta_2^{q+1-p},\\
\end{align*}
where we used $\Delta^{q+1}\le \Delta_2^{q+1}=(\Delta\wedge \Delta_2)^p\Delta_2^{q+1-p}$ if $\Delta_2\le \Delta$, and $|q+1|^{-1}\ge (p-1-q)^{-1}$ in the last line.  
\gdm

\medskip
We let $p'_t(x)={d\over dx}p_t(x)$.
\medskip

\begin{lemma}\label{pt'triv}
\[|p'_t(z)|\le c_{\ref{pt'triv}}t^{-1/2}p_{2t}(z).\]
\end{lemma}
\paragraph{Proof.} Trivial.
\gdm

\begin{lemma}\label{ptbnds}
(a) There is a $c_{\ref{ptbnds}}$ so that for any $s<t\le t', x,x'\in\Re$,
\begin{equation} \label{psq}\int(p_{t'-s}(y-x')-p_{t-s}(y-x))^2\,dy\le c_{\ref{ptbnds}}(t-s)^{-1/2}\Bigl[1\wedge{d((t,x),(t',x'))^2\over t-s}\Bigr].
\end{equation}
\noindent(b) For any $R>2$ there is a $c_{\ref{ptbnds}}(R)$ so that for any $0\le p, r\le R$, $\eta_0,\eta_1\in(1/R,1/2)$, $0\le s<t\le t'\le R$, $x,x'\in\Re$, 
\begin{align}\nonumber\int &e^{r|y-x|}|y-x|^p(p_{t-s}(y-x)-p_{t'-s}(y-x'))^21(|y-x|>(t'-s)^{1/2-\eta_0}\vee2|x'-x|)dy\\
&\le c_{\ref{ptbnds}}(R)(t-s)^{-1/2}\exp\{-\eta_1(t'-s)^{-2\eta_0}/33\}\Bigl[1\wedge{d((t,x),(t',x'))^2\over t-s}\Bigr]^{1-(\eta_1/2)}.
\end{align}
\end{lemma}
\paragraph{Proof.} (a) Let $f(u)=u^{-1/2}$. By Chapman-Kolmogorov, the integral in \eqref{psq} equals
\begin{align} \nonumber p&_{2(t'-s)}(0)+p_{2(t-s)}(0)-2p_{t'-s+t-s}(0)+2(p_{t'-s+t-s}(0)-p_{t'-s+t-s}(x-x'))\\
&\nonumber\le (2\pi)^{-1/2}\Bigl[|f(2(t'-s))+f(2(t-s))-2f(t'-s+t-s)|\\
&\nonumber\phantom{\le (2\pi)^{-1/2}\Bigl[}+(t'-s+t-s)^{-1/2}\Bigl(1-\exp\Bigl\{{-(x-x')^2\over 2(t'-s+t-s)}\Bigr\}\Bigr)\\
&\label{ptexp}\equiv (2\pi)^{-1/2}[T_1+T_2].
\end{align}
Clearly $T_2$ is at most $(t'-s)^{-1/2}\Bigl[1\wedge {|x-x'|^2\over t'-s}\Bigr]$.  If $0<u\le u'$, $0\le f(u)-f(u')\le u^{-1/2}\wedge[u^{-3/2}|u'-u|]$ (by the Mean Value Theorem) and so
\begin{equation}\label{t1bnd}
T_1\le \Bigl[\sqrt 2(t-s)^{-1/2}\wedge \Bigl({2\over (2(t-s))^{3/2}}|t'-t|\Bigr)\Bigr]
\end{equation}
Use the above bounds on $T_1$ and $T_2$ in \eqref{ptexp} to complete the proof of (a). 

\medskip

\noindent(b) This proof is very similar to that of Lemma~\ref{pt'bnds} (b) below and so is omitted.  There are some minor differences leading to the factor of $1/33$ (rather than the $1/64$ in Lemma~\ref{pt'bnds} (b))--e.g., the much simpler analogue of \eqref{p'upperbound} below has $p_u(w)$ on the right side.
\gdm

\medskip

\begin{lemma}\label{pt'bnds}
(a) There is a $c_{\ref{pt'bnds}}$ so that for any $s<t\le t', x,x'\in\Re$,
\begin{equation}\label{pt'ident0} \int(p'_{t'-s}(y-x')-p'_{t-s}(y-x))^2\,dy\le c_{\ref{pt'bnds}}(t-s)^{-3/2}\Bigl[1\wedge{d((t,x),(t',x'))^2\over t-s}\Bigr].
\end{equation}
\noindent(b) For any $R>2$ there is a $c_{\ref{pt'bnds}}(R)$ so that for any $0\le p, r\le R$, $\eta_0,\eta_1\in(1/R,1/2)$, $0\le s<t\le t'\le R$, $x,x'\in\Re$, 
\begin{align}\nonumber\int &e^{r|y-x|}|y-x|^p(p'_{t-s}(y-x)-p'_{t'-s}(y-x'))^21(|y-x|>(t'-s)^{1/2-\eta_0}\vee2|x'-x|)dy\\
&\label{pt'bndseq}\le c_{\ref{pt'bnds}}(R)(t-s)^{-3/2}\exp\{-\eta_1(t'-s)^{-2\eta_0}/64\}\Bigl[1\wedge{d((t,x),(t',x'))^2\over t-s}\Bigr]^{1-(\eta_1/2)}.
\end{align}
\end{lemma}
\paragraph{Proof.} (a) We first claim that
\begin{equation}\label{pt'ident1}
\int p'_t(w)p'_t(w-x)dw=\Bigl({t\over 2}-{x^2\over 4}\Bigr){p_{2t}(x)\over t^2}.
\end{equation}
To see this, first do a bit of algebra to get
\begin{equation}\label{ptident1}p_t(w)p_t(w-x)=p_{t/2}(w-(x/2))p_{2t}(x).
\end{equation}
Therefore the left-hand side of \eqref{pt'ident1} equals
\begin{align*}
\int {w(w-x)\over t^2}p_t(w)p_t(w-x)dw&=\int {w(w-x)\over t^2}p_{t/2}(w-(x/2))dwp_{2t}(x)\\
&=\int \Bigl(u^2-{x^2\over 4}\Bigr) p_{t/2}(u)du{p_{2t}(x)\over t^2}\quad (u=w-(x/2))\\
&=\Bigl({t\over 2}-{x^2\over 4}\Bigr){p_{2t}(x)\over t^2},
\end{align*}
giving the right-hand side of \eqref{pt'ident1}.

Next we claim that 
\begin{equation}\label{pt'ident2}
\int p'_{t'}(w-x)p'_t(w-x)dw=(t+t')^{-1}p_{t'+t}(0).
\end{equation}
Some algebra shows that
\begin{equation}\label{ptident2}
p_{t'}(w)p_t(w)=p_{t'+t}(0)p_{{t't\over t+t'}}(w).
\end{equation}
Therefore the left-hand side of \eqref{pt'ident2} equals
\begin{align*}\int {w^2\over t't}p_{t'}(w)p_{t}(w)dw
=\int {w^2\over t't}p_{{tt'\over t+t'}}(w)dwp_{t'+t}(0)
=(t+t')^{-1}p_{t'+t}(0),
\end{align*}
and we have \eqref{pt'ident2}.

The left-hand side of \eqref{pt'ident0} is bounded by
\begin{align}\label{Tdecomp1}
2&\Bigl[\int(p'_{t'-s}(y-x')-p'_{t'-s}(y-x))^2dy+\int(p'_{t'-s}(y-x)-p'_{t-s}(y-x))^2dy\Bigr]\\
\nonumber&\equiv 2(T_1+T_2).
\end{align}
Now expand $T_1$ and use \eqref{pt'ident1} to see that 
\begin{align*}
T_1&=2\int p'_{t'-s}(y-x)^2dy-2\int p'_{t'-s}(w)p'_{t'-s}(w+x-x')dw\\
&=(t'-s){p_{2(t'-s)}(0)\over (t'-s)^2}-\Bigl((t'-s)-{(x-x')^2\over 2}\Bigr){p_{2(t'-s)}(x'-x)\over (t'-s)^2}\\
&=(t'-s)^{-1}(p_{2(t'-s)}(0)-p_{2(t'-s)}(x'-x))+{(x'-x)^2\over 2(t'-s)^2}p_{2(t'-s)}(x'-x)\\
&\le (t'-s)^{-3/2}\Bigl[1\wedge {(x'-x)^2\over t'-s}\Bigr]+(\sup_z(z e^{-z})(t'-s)^{-3/2})\wedge((x'-x)^2(t'-s)^{-5/2})\\
&\le c_0  (t'-s)^{-3/2}\Bigl[1\wedge {(x'-x)^2\over t'-s}\Bigr].
\end{align*}
Finally let $g(u)=u^{-3/2}$, and expand $T_2$ and use \eqref{pt'ident2} to conclude
\begin{align*}
T_2&=(2(t'-s))^{-1}p_{2(t'-s)}(0)+(2(t-s))^{-1}p_{2(t-s)}(0)-2(t-s+t'-s)^{-1}p_{t-s+t'-s}(0)\\
&=(2\pi)^{-1/2}[g(2(t'-s))+g(2(t-s))-2g(t'-s+t-s)]\\
&\le (2\pi)^{-1/2}2\Bigl[(2(t-s))^{-3/2}\wedge(2(t-s))^{-5/2}|t'-t|\Bigr].
\end{align*}
The last inequality follows as in \eqref{t1bnd}.  Use the above bounds on $T_1$ and $T_2$ in \eqref{Tdecomp1}
to complete the proof of (a). 

\medskip
\noindent (b) Note that $|y-x|>(t'-s)^{1/2-\eta_0}\vee2|x'-x|$ implies that
\begin{equation}\label{incbound1}
|y-x'|\ge |y-x|-|x'-x|\ge |y-x|/2\ge {(t'-s)^{1/2-\eta_0}\over 2}
\end{equation}
and in particular from the second inequality,
\begin{equation}\label{incbound2}
|y-x|\le 2|y-x'|.\end{equation}
Assume $p,r,\eta_i,s,t,t'$ as in (b) and let $d=d((t,x),(t',x'))$.  By H\"older's inequality and then (a),  the integral on the left-hand side of \eqref{pt'bndseq} is at most
\begin{align}
\nonumber &\Bigl[\int(p'_{t-s}(y-x)-p'_{t'-s}(y-x'))^2dy\Bigr]^{1-\eta_1/2}\\
\nonumber&\phantom{\Bigl[}\times\Bigl[\int e^{{2r\over \eta_1}|y-x|}|y-x|^{2p/\eta_1}(p'_{t-s}(y-x)-p'_{t'-s}(y-x'))^21(|y-x|>(t'-s)^{1/2-\eta_0}\vee 2|x-x'|)dy\Bigr]^{\eta_1/2}\\
\nonumber&\le c_1 (t-s)^{-{3\over 2}+{3\eta_1\over 4}}\Bigl[1\wedge{d^2\over t-s}\Bigr]^{1-{\eta_1\over 2}}\Bigl[\int e^{{2r\over \eta_1}|y-x|}|y-x|^{2p/\eta_1}(p'_{t-s}(y-x)^2+p'_{t'-s}(y-x')^2)\\
\nonumber&\phantom{\le c_1 (t-s)^{-{3\over 2}+{3\eta_1\over 4}}\Bigl[1\wedge{d^2\over t-s}\Bigr]^{1-{\eta_1\over 2}}\Bigl[}\times1(|y-x|>(t'-s)^{1/2-\eta_0}\vee 2|x-x'|)dy\Bigr]^{\eta_1/2}\\
\nonumber&\le c_2(R) (t-s)^{-{3\over 2}+{3\eta_1\over 4}}\Bigl[1\wedge{d^2\over t-s}\Bigr]^{1-{\eta_1\over 2}}\Bigl[\int e^{{2r\over \eta_1}|w|}|w|^{{2p\over \eta_1}}p'_{t-s}(w)^21(|w|>(t-s)^{1/2-\eta_0})dw\\
\label{ubound1}&\phantom{\le c_2(R) (t-s)^{-{3\over 2}+{3\eta_1\over 4}}\Bigl[1\wedge{d^2\over t-s}\Bigr]^{1-{\eta_1\over 2}}\Bigl[}+\int e^{{4r\over \eta_1}|w|}|w|^{{2p\over \eta_1}}p'_{t'-s}(w)^21(|w|>(t'-s)^{1/2-\eta_0}/2)dw\Bigr]^{\eta_1/2},
\end{align}
where in the last we used \eqref{incbound1} and \eqref{incbound2}.  

If $|w|>{1\over 2}u^{1/2-\eta_0}$, then by Lemma~\ref{pt'triv}, 
\begin{align}\nonumber p'_u(w)^2\le c_{\ref{pt'triv}}^2u^{-1}p_{2u}(w)^2&\le c_{\ref{pt'triv}}^2u^{-3/2}e^{u^{-2\eta_0}/16}p_{2u}(w)\\
\label{p'upperbound}&\le c_3(R)p_{2u}(w).
\end{align}
Use this to show that
\begin{align*}
&\Bigl[\int e^{{4r\over \eta_1}|w|}|w|^{{2p\over \eta_1}}p'_{t'-s}(w)^21(|w|>(t'-s)^{1/2-\eta_0}/2)dw\Bigr]^{\eta_1/2}\\
&\le c_4(R)\Bigl[\int e^{{4r\over \eta_1}|w|}|w|^{{2p\over \eta_1}}p_{2(t'-s)}(w)1(|w|>(t'-s)^{1/2-\eta_0}/2)dw\Bigr]^{\eta_1/2}\\
&\le c_5(R)E(e^{{8r\over\eta_1}\sqrt{2(t'-s)}|B_1|}|B_1|^{{4p\over\eta_1}})^{\eta_1/4}P(|B_1|>{1\over 2\sqrt 2}(t'-s)^{-\eta_0})^{\eta_1/4}\quad\hbox{(by H\"older)}\\
&\le c_6(R)\exp\{-\eta_1(t'-s)^{-2\eta_0}/64\}.
\end{align*}
Use the same bound with $t$ in place of $t'$ to see that the right-hand side of \eqref{ubound1}, and hence of \eqref{pt'bndseq}, is at most
\begin{equation*}
c_7(R)(t-s)^{3\eta_1/4}(t-s)^{-3/2}\Bigl[1\wedge{d^2\over t-s}\Bigr]^{1-{\eta_1\over 2}}\exp\{-(t'-s)^{-2\eta_0}\eta_1/64\}.
\end{equation*}
The result follows because $(t-s)^{3\eta_1/4}\le R$.
\gdm

\section{Local Bounds on the Difference of Two Solutions}\label{sec4}
\setcounter{equation}{0}
\setcounter{theorem}{0}
This section is devoted to establishing the local bounds on the difference of two solutions to~(\ref{wspde}).
These bounds are crucial for the construction of the stopping times in Proposition~\ref{tildeJ}, which is then carried out in Section~\ref{sec5}. We continue to assume throughout this Section that $b\equiv0$.  Recall that  $X^1,X^2$ are two solutions as in Section~\ref{secmainres}, $u=X^1-X^2$, and we assume the hypotheses of Theorem~\ref{theorem:unique} as well as \eqref{holder'}.

We refine the earlier set $Z(N,K)$ and define, for $K,N,n\in\N$ and $\beta \in (0,1/2]$,
\begin{align*}
Z(N,n,K,\beta)(\omega)=&\{(t,x)\in[0,T_K]\times[-K,K]: \hbox{ there is a }(\hat t_0,\hat x_0)\in[0,T_K]\times\Re\hbox{ such that }\\
&d((\hat t_0,\hat x_0),(t,x))\le 2^{-N},\ |u(\hat t_0,\hat x_0)|\le a_n\wedge(\sqrt{ a_n}2^{-N}),\hbox{ and }
|u'_{1,{a_n}}(\hat t_0,\hat x_0)|\le 
 a_n^\beta\},
\end{align*}
and for $\beta=0$ define
$Z(N,n,K,0)(\omega)=Z(N,n,K)(\omega)$ as above, but with the condition on $|u'_{1,{a_n}}(\hat t_0,\hat x_0)|$ omitted.

Recalling $\gamma<1$, we let
\begin{equation}\label{gammamdef} 
\gamma_m={(\gamma-1/2)(1-\gamma^m)\over 1-\gamma}+1, \ \ \tilde\gamma_m=\gamma_m\wedge 2
\end{equation} so that we have the recursion relation
\begin{equation}\label{gammarecursion}
\gamma_{m+1}=\gamma\gamma_m+1/2,\quad \gamma_0=1.
\end{equation}
Clearly $\gamma_m$ increases to $\gamma_\infty=(\gamma-1/2)(1-\gamma)^{-1}+1=(2(1-\gamma))^{-1}>2$ and so we may define a finite natural number, $\overline m>1$, by
\begin{equation}\label{barmdef}
\om=\min\{m:\gamma_{m+1}>2\}=\min\{m:\gamma\gamma_m>3/2\}.
\end{equation}

\noindent{\bf Definition.} A collection of $[0,\infty]$-valued random variables, $\{N(\alpha):\alpha\in A\}$, is stochastically bounded uniformly in $\alpha$ iff 
\[\lim_{M\to\infty}\sup_{\alpha\in A}P(N(\alpha)\ge M)=0.\]

Finally we introduce the condition whose proof will be the goal of this section.   Recall that $K_1$ is as in \eqref{Gamma}.  
\paragraph{Property  $(P_m)$.} For $m\in\Z_+$ we will let $(P_m)$ denote the following property:
\begin{align*} &\hbox{For any }n\in \N, \xi,\ve_0\in(0,1), K\in\N^{\ge K_1} \hbox{ and }\beta\in[0,1/2], \hbox{ there is an }N_1(\omega)=N_1(m,n,\xi,\ve_0,K,\beta)\\
&\hbox{in $\N$ a.s. such that for all }
N\ge N_1,\hbox{ if }(t,x)\in Z(N,n,K,\beta),\ t'\le T_K,\hbox{ and } d((t,x),(t',x'))\le 2^{-N},\hbox{ then }\\
&|u(t',x')|\le a_n^{-\ve_0}2^{-N\xi}\Bigl[(\sqrt{a_n}\vee 2^{-N})^{\tilde\gamma_m-1}+a_n^\beta1(m>0)\Bigr].\\
&\hbox{Moreover }N_1\hbox{ is stochastically bounded uniformly in }(n,\beta).
\end{align*}

\medskip
Here is the main result of this section:

\begin{proposition}\label{Pminduction}
For any $m\le \om+1$, $(P_m)$ holds.
\end{proposition}
\paragraph{Proof.} $(P_0)$ is an easy consequence of Theorem~\ref{theorem:lipmod}, as we now show--and we may even take $\ve_0=0$.  
Let $\xi\in(0,1)$ and apply Theorem~\ref{theorem:lipmod} with $\xi'=(\xi+1)/2$ in place of $\xi$.  If $(t,x)\in Z(N,K) (\supset Z(N,n,K,\beta))$ and $(\hat t_0,\hat x_0)$ is as in the definition of $Z(N,K)$, then $(\hat t_0,\hat x_0)\in Z(N,K+1)$ (we need $K+1$ since $|\hat x_0|\le K+1$).  Theorem~\ref{theorem:lipmod} implies that if $N\ge N_0(\xi',K+1)\vee 
4(1-\xi)^{-1}\equiv N_1(0,\xi,K)$ (it doesn't depend on $(n,\beta)$ and there is no $\ve_0$), then
\begin{align*}
|u(t,x)|\le 2^{-N\xi'}+|u(\hat t_0,\hat x_0)|\le 2^{1-N\xi'}.
\end{align*}
If $(t',x')$ is as in $(P_0)$, the above and Theorem~\ref{theorem:lipmod} imply
\begin{eqnarray*}
 |u(t',x')|\le |u(t,x)|+|u(t',x')-u(t,x)|\le 2^{1-N\xi'}+2^{-N\xi'}
 \le2^{2-N\xi'}\le 2^{-N\xi},
\end{eqnarray*}
where the last inequality holds because $N\ge 4(1-\xi)^{-1}$.
  $(P_0)$ follows.

The induction step will require some additional continuity results which also will be used directly in the next section.  We start by noting that $(P_m)$ easily gives some global bounds on $|u|$. 

\begin{lemma}\label{uglobalbound} Let $0\le m\le \om+1$ and assume $(P_m)$.  For any $n,\xi,\ve_0,K$ and $\beta$ as in $(P_m)$, if $\overline d_N=2^{-N}\vee d((s,y),(t,x))$ and $\sqrt{C_{\ref{uglobalbound}}(\omega)}=(4a_n^{-\ve_0}+2^{2N_1(\omega)}2Ke^K)$, then for any $N\in\N$,  
\begin{equation}\label{Zcond}
\hbox{on }
\{\omega:N\ge N_1(m,n,\xi,\ve_0,K,\beta), (t,x)\in Z(N,n,K,\beta)\},
\end{equation}
 we have
\begin{align}
 |u(s,y)|
\le& \sqrt{C_{\ref{uglobalbound}}}e^{|y-x|}\overline d_N^\xi \\
&\nonumber\times\Bigl[(\sqrt{a_n}\vee\overline d_N)^{\tilde\gamma_m-1}+1(m>0)a_n^\beta
\Bigr]\hbox{ for all }s\le T_K\hbox{ and }y\in\Re.
\end{align}
\end{lemma}
\paragraph{Proof.} Assume $N,\omega,t,x$ are as in \eqref{Zcond}.

\noindent Case 1.  $d\equiv d((s,y),(t,x))\le 2^{-N_1}$.

\noindent If $d>2^{-N}$ choose $N_1\le N'<N$ so that $2^{-N'-1}<d\le 2^{-N'}$, and if $d\le 2^{-N}$ set $N'=N$.  Then $(t,x)\in Z(N',n,K,\beta)$, $d\le 2^{-N'}\le 2^{-N}\vee 2d\le 2\overline d_N$ and so by $(P_m)$ for $s\le T_K$,
\begin{align*}
|u(s,y)|&\le a_n^{-\ve_0}2^{-N'\xi}\Bigl[(\sqrt{a_n}\vee 2^{-N'})^{\tilde\gamma_m-1}+1(m>0)a_n^\beta\Bigr]\\
&\le 4a_n^{-\ve_0}(\overline d_N)^\xi\Bigl[(\sqrt{a_n}\vee \overline d_N)^{\tilde\gamma_m-1}+1(m>0)a_n^\beta\Bigr].
\end{align*}

\noindent Case 2. $d>2^{-N_1}$.

\noindent As $K\ge K_1$, for $s\le T_K$,
\begin{align*}
|u(s,y)|\le 2Ke^{|y|}&\le 2Ke^{|y|}(d2^{N_1})^{\xi+\tilde\gamma_m-1}\\
&\le 2K e^Ke^{|y-x|}2^{2N_1}(\overline d_N)^{\xi+\tilde\gamma_m-1}.
\end{align*}

The Lemma follows from the above two bounds.
\gdm

\begin{remark} \label{m0case}If $m=0$ we may set $\ve_0=0$ in the above and $N_1$ will not depend on $(n,\ve_0,\beta)$ by the above proof of $(P_0)$.
\end{remark}

To carry out the induction we first use $(P_m)$ to obtain a local modulus of continuity for $F_\delta$.  
From Lemma~\ref{G'}, we have for $s\le t\le t'$ and $s'\le t'$
\begin{align}\label{Fdecomp}
\nonumber |F_\delta(s,t,x)-F_\delta(s',t',x')|\le &|F_\delta(s,t',x')-F_\delta(s',t',x')|+|F_\delta(s,t',x')-F_\delta(s,t,x)|\\
=&\Bigl|\int_{(s-\delta)^+}^{(s'-\delta)^+}\int p'_{t'-r}(y-x')D(r,y)W(dr,dy)\Bigr|\\
\nonumber&+\Bigl|\int_0^{(s-\delta)^+}(p'_{t'-r}(y-x')-p'_{t-r}(y-x))D(r,y)W(dr,dy)\Bigr|.
\end{align}
This decomposition and \eqref{Dbound} suggest we introduce the following square functions for $\eta_0\in(0,1/2)$ and $\delta\in(0,1]$, and $s\le t\le t', s'\le t'$:
\begin{align*}
&Q_{T,\delta}(s,s',t',x')=\int_{(s\wedge s'-\delta)^+}^{(s\vee s'-\delta)^+}\int p'_{t'-r}(y-x')^2e^{2R_1|y|}|u(r,y)|^{2\gamma}dydr,\\
&Q_{S,1,\delta,\eta_0}(s,t,x,t',x')=\int _0^{(s-\delta)^+}\int 1(|y-x|>(t'-r)^{1/2-\eta_0}\vee 2|x-x'|)\\
&\phantom{Q_{S,1,\eta_0}(s,t,x,t',x')=\int _0^{(s-\delta)^+}\int}\times(p'_{t'-r}(y-x')-p'_{t-r}(y-x))^2e^{2 R_1|y|}|u(r,y)|^{2\gamma}dydr,\\
&Q_{S,2,\delta,\eta_0}(s,t,x,t',x')=\int _0^{(s-\delta)^+}\int 1(|y-x|\le (t'-r)^{1/2-\eta_0}\vee 2|x-x'|)\\
&\phantom{Q_{S,1,\delta,\eta_0}(s,t,x,t',x')=\int _0^{(s-\delta)^+}\int}\times(p'_{t'-r}(y-x')-p'_{t-r}(y-x))^2e^{2R_1|y|}|u(r,y)|^{2\gamma}dydr.
\end{align*}

\begin{lemma}\label{QS1bound}
For all $K\in\N^{\ge K_1}, R>2$ there is a $c_{\ref{QS1bound}}(K,R)$ and an $N_{\ref{QS1bound}}=N_2(K,\omega)\in \N$ a.s. so that for all $\eta_0,\eta_1\in(1/R,1/2), \delta\in(0,1],\beta\in[0,1/2]$ and $N,n\in\N$, for any $(t,x)\in \Rp\times\Re$, on
\begin{equation}\label {omegacond2}
\{\omega:(t,x)\in Z(N,n,K,\beta), N\ge N_{\ref{QS1bound}}\},
\end{equation}
\begin{equation}
Q_{S,1,\delta,\eta_0}(s,t,x,t',x')\le c_{\ref{QS1bound}}2^{4N_{\ref{QS1bound}}(\omega)}\Bigl[d^{2-\eta_1}+(d\wedge\sqrt\delta)^{2-\eta_1}\delta^{-3/2}(d\wedge1)^{4\gamma}\Bigr]\hbox{ for all }s\le t\le t'\hbox{ and }x'\in\Re.
\end{equation}
Here $d=d((t',x'),(t,x))$. 
\end{lemma}
\paragraph{Proof} We let $N_{\ref{QS1bound}}(K,\omega)=N_1(0,3/4,K)$, that is we recall from Remark~\ref{m0case} that for $m=0$, $N_1$ depends only on $\xi$ and $K$ and we take $\xi=3/4$. We may assume $\delta<s$ as the left-hand side is $0$ otherwise. Then for $\omega$ as in \eqref{omegacond2} and $s\le t\le t'$, Lemma~\ref{uglobalbound} with $m=0$ implies 
\begin{align*}
&Q_{S,1,\delta,\eta_0}(s,t,x,t',x')\\
&\le C_{\ref{uglobalbound}}(\omega)\int  _0^{(s-\delta)^+}\int 1(|y-x|>(t'-r)^{1/2-\eta_0}\vee 2|x-x'|)(p'_{t'-r}(y-x')-p'_{t-r}(y-x))^2\\
&\phantom{\le C_{\ref{uglobalbound}}(\omega)\int  _0^{(s-\delta)^+}}\times e^{2 R_1|y|}e^{2|y-x|}(2^{-N}\vee(\sqrt{t-r}+|y-x|))^{\gamma3/2}dyds\\
&\le C_{\ref{uglobalbound}}(\omega)\int  _0^{(s-\delta)^+}\int 1(|y-x|>(t'-r)^{1/2-\eta_0}\vee 2|x-x'|)(p'_{t'-r}(y-x')-p'_{t-r}(y-x))^2\\
&\phantom{\le C_{\ref{uglobalbound}}(\omega)\int  _0^{(s-\delta)^+}}\times e^{2 R_1K}e^{2 (R_1+1)|y-x|} 
[ 2K^{\gamma 3/4}+2|y-x|^{\gamma 3/2}]dyds\\
&\le C_{\ref{uglobalbound}}(\omega) c_0(K,R)\int _0^{s-\delta}(t-r)^{-3/2}\exp\Bigl\{{-\eta_1(t'-r)^{-2\eta_0}\over 64}\Bigr\}\Bigl[1\wedge{d^2\over t-r}\Bigr]^{1-\eta_1/2}dr.
\end{align*}
In the last line we have used Lemma~\ref{pt'bnds}(b).  Use the trivial bound (recall $r\le s\le t\le t'$)
\[\exp\Bigl\{{-\eta_1(t'-r)^{-2\eta_0}\over 64}\Bigr\}\le \exp\Bigl\{{-\eta_1(t'-t)^{-2\eta_0}\over 128}\Bigr\}+\exp\Bigl\{{-\eta_1(t-r)^{-2\eta_0}\over 128}\Bigr\},\]
and then Lemma~\ref{Jbnd} in the above, to bound $Q_{S,1,\delta,\eta_0}(s,t,x,t',x')$ by
\begin{align*}
&C_{\ref{uglobalbound}}(\omega) c_0(K,R)\Bigl[\int_0^{s-\delta}(t-r)^{-3/2}\Bigl[1\wedge{d^2\over t-r}\Bigr]^{1-\eta_1/2}dr\exp\Bigl\{{-\eta_1(t'-t)^{-2\eta_0}\over 128}\Bigr\}\\
&\phantom{C_{\ref{uglobalbound}}(\omega) c_0(K,R)\Bigl[}+\int _0^{s-\delta}(t-r)^{-3/2}\Bigl(1\wedge {d^2\over t-r}\Bigr)^{1-\eta_1/2}\exp\Bigl\{{-\eta_1(t-r)^{-2\eta_0}\over 128}\Bigr\}dr\Bigr]\\
&\le C_{\ref{uglobalbound}}(\omega) c_1(K,R)\Bigl[(d^2\wedge\delta)^{1-\eta_1/2}\delta^{-3/2}\exp\{{-\eta_1(t'-t)^{-2\eta_0}\over 128}\Bigr\}+C_1(R)\int _0^{s-\delta}\Bigl(1\wedge {d^2\over t-r}\Bigr)^{1-\eta/2}dr\Bigr]\\
&\le  C_{\ref{uglobalbound}}(\omega) c_2(K,R)\Bigl[ (d\wedge \sqrt\delta)^{2-\eta_1}\delta^{-3/2}(d\wedge 1)^{4\gamma}+ d^{2-\eta_1}\Bigr].
\end{align*}
Now since we may set $\ve_0=0$ in the formula for $C_{\ref{uglobalbound}}$ by Remark~\ref{m0case}, the result follows.
\gdm

\begin{lemma}\label{QS2bound} Let $0\le m\le \om+1$ and assume $(P_m)$.  For any $K\in\N^{\ge K_1}, R>2, n\in\N, \ve_0\in(0,1)$, and $\beta\in [0,1/2]$ there is a $c_{\ref{QS2bound}}(K,R)$ and $N_{\ref{QS2bound}}=N_{\ref{QS2bound}}(m,n,R,\ve_0,K,\beta)(\omega)\in\N$ a.s. such that for any $\eta_1\in(R^{-1},1/2), \eta_0\in(0,\eta_1/32), \delta\in[{a_n},1]$, $N\in\N$, and $(t,x)\in\Rp\times\Re$, 
on
\begin{equation}\label {omegacond3}
\{\omega:(t,x)\in Z(N,n,K,\beta), N\ge N_{\ref{QS2bound}}\},
\end{equation}
\begin{align*}
&Q_{S,2,\delta,\eta_0}(s,t,x,t',x')\\
&\le c_{\ref{QS2bound}}(K,R)[a_n^{-2\ve_0}+2^{4N_{\ref{QS2bound}}}]\Bigl[d^{2-\eta_1}[{\bar{\delta}_N}^{(\gamma\gamma_m-3/2)\wedge 0}+a_n^{2\beta\gamma}\bar{\delta}_N^{(\gamma-{3\over 2})\wedge 0}]\\
&\phantom{\le C_{\ref{QS2bound}}(K,R)[a_n^{-2\ve_0}+2^{4N_3}]\Bigl[}+(d\wedge\sqrt\delta)^{2-\eta_1}\delta^{-3/2}[{\bar{d}_N}^{2\gamma\tilde\gamma_m}+a_n^{2\beta\gamma}\bar{d}_N^{2\gamma}]\Bigr]\\
&\phantom{\le C_{\ref{QS2bound}}(K,R)[a_n^{-2\ve_0}+2^{4N_3}]\Bigl[}\hbox{ for all }s\le t\le t'\le K, |x'|\le K+1.
\end{align*}
Here $d=d((t,x),(t',x'))$, $\bar d_N=d\vee 2^{-N}$ and $\bar\delta_N=\delta\vee {\bar d_N}^2$.  Moreover $N_{\ref{QS2bound}}$ is stochastically bounded uniformly in $(n,\beta)$.
\end{lemma}
\paragraph{Proof.} Let $\xi=1-(8R)^{-1}\in (15/16,1)$ and define $N_{\ref{QS2bound}}=N_1(m,n,\xi,\ve_0,K,\beta)$ so that the last statement is immediate from $(P_m)$.  We may assume $s\ge\delta$, or the left-hand side is $0$.  As $\delta\ge {a_n}$, when we use Lemma~\ref{uglobalbound} to bound $u(r,y)$ in the integral defining $Q_{S,2,\delta,\eta_0}$, we have $d((r,y),(t,x))\ge \sqrt{a_n}$ and so we may drop the $\max$ with $\sqrt{a_n}$.  So for $\omega$ as in \eqref{omegacond3}, $s\le t\le t'$ and $|x'|\le K+1$, Lemma~\ref{uglobalbound} implies that
\begin{align*}&Q_{S,2,\delta,\eta_0}(s,t,x,t',x')\\
&\le C_{\ref{uglobalbound}}\int_0^{s-\delta}\int(p'_{t'-r}(y-x')-p'_{t-r}(y-x))^2e^{2 R_1K}e^{2 (R_1+1)2(2K+1)}dy\\
&\phantom{\le C_{\ref{uglobalbound}}\int_0^{s-\delta}}\times[2^{-N}\vee((t-r)^{1/2}+(t'-r)^{1/2-\eta_0}\vee(2|x-x'|))]^{2\gamma\xi}\\
&\phantom{\le C_{\ref{uglobalbound}}\int_0^{s-\delta}}\times\Bigl\{[2^{-N}\vee((t-r)^{1/2}+(t'-r)^{1/2-\eta_0}\vee(2|x-x'|))]^{\tilde\gamma_m-1}+a_n^\beta\Bigr\}^{2\gamma}\,dr.
\end{align*}
Let $\gamma'=\gamma(1-2\eta_0)$.  Recall that $t\le t'\le K$, $|x|\le K$ and $|x'|\le K+1$, so that $\sqrt{t-r}\le K^{\eta_0}(t'-r)^{1/2-\eta_0}$ and $|x-x'|\le (2K+1)|x-x'|^{1-2\eta_0}$.  Use this and Lemma~\ref{pt'bnds}(a) to see that the above is at most
\begin{align}
\nonumber&c_1(K)C_{\ref{uglobalbound}}\int_0^{s-\delta}(t-r)^{-3/2}\Bigl(1\wedge {d^2\over t-r}\Bigr)(2^{-2N\gamma}\vee(t'-r)^{\gamma'}\vee|x-x'|^{2\gamma'})^\xi\\
\label{Qexp}&\phantom{c_1(K)C_{\ref{uglobalbound}}\int_0^{s-\delta}}\times\Bigl[2^{-2N\gamma(\tilde\gamma_m-1)}\vee(t'-r)^{\gamma'(\tilde\gamma_m-1)}\vee|x'-x|^{2\gamma'(\tilde \gamma_m-1)}+a_n^{2\beta\gamma}\Bigr]\,dr.
\end{align}
Note that
\begin{equation}
2^{-2N\gamma}\vee(t'-r)^{\gamma'}\vee|x'-x|^{2\gamma'}\le 2^{-2N\gamma'}\vee d^{2\gamma'}+(t-r)^{\gamma'}\le 2[\bar d_N^{2\gamma'}\vee(t-r)^{\gamma'}].
\end{equation}
Use this to bound the summands in \eqref{Qexp} and conclude that
\begin{align}
\nonumber&Q_{S,2,\delta,\eta_0}(s,t,x,t',x')\\
\nonumber&\le c_2(K)C_{\ref{uglobalbound}}\int_0^{s-\delta}(t-r)^{-3/2}\Bigl(1\wedge {d^2\over t-r}\Bigr)[\bar d_N^{2}\vee(t-r)]^{\gamma'\xi}\\
\nonumber&\phantom{\le c_2(K)C_{\ref{uglobalbound}}\int_0^{s-\delta}}\times\Bigl[[\bar d_N^2\vee (t-r)]^{\gamma'(\tilde \gamma_m-1)}+a_n^{2\beta\gamma}\Bigr]dr\\
\nonumber&\le c_2(K)C_{\ref{uglobalbound}}\Bigl\{\int_0^{t-\bar\delta_N}(t-r)^{\gamma'\xi-3/2}\Bigl(1\wedge {d^2\over t-r}\Bigr)\Bigl[(t-r)^{\gamma'(\tilde\gamma_m-1)}+a_n^{2\beta\gamma}\Bigr]dr\\
\nonumber&\phantom{\le c_2(K)C_{\ref{uglobalbound}}\Bigl\{}+\int_0^{t-\delta}(t-r)^{-3/2}\Bigl(1\wedge {d^2\over t-r}\Bigr)dr\,\bar d_N^{2\gamma'\xi}\Bigl[\bar d_N^{2\gamma'(\tilde\gamma_m-1)}+a_n^{2\beta\gamma}\Bigr]\Bigr\}\\
\label{QS2decomp}&\equiv c_2(K)C_{\ref{uglobalbound}}\Bigl\{I_1+I_2\Bigr\}.
\end{align}
Apply Lemma~\ref{Jbnd}(c) to see that 
\begin{align}\label{QS2decomp1}
I_2&\le (d\wedge\sqrt\delta)^2\delta^{-3/2}\bar d_N^{2\gamma'\xi}
\Bigl[\bar d_N^{2\gamma'(\tilde\gamma_m-1)}
+a_n^{2\beta\gamma}\Bigr].
\end{align}
In the integral defining $I_1$ we may drop the minimum with $1$ and, adding a $\log(1/\bar \delta_N)$ factor just in case the exponent on $u$ is $-1$, we arrive at 
\begin{align}
\nonumber I_1&\le d^2\Bigl[\int_{\bar \delta_N}^tu^{\gamma'(\tilde\gamma_m+\xi-1)-5/2}du+a_n^{2\beta\gamma}\int_{\bar \delta_N}^tu^{\gamma'\xi-5/2}du\Bigr]\\
\label{QS2decomp2}&\le c_3(K)d^2\log(1/\bar\delta_N)\Bigl[\bar\delta_N^{[\gamma'(\tilde\gamma_m+\xi-1)-3/2]\wedge0}+a_n^{2\beta\gamma}\bar\delta_N^{[\gamma'\xi-3/2]\wedge0}\Bigr]
\end{align}
The $\log 1/\bar\delta_N$ is bounded by $c(R)d^{-\eta_1/2}$.
A bit of arithmetic shows that our conditions $\eta_0\le \eta_1/32$, $1-\xi=1/8R$ and $\eta_1>1/R$ allow us to shift $\xi$ to $1$ and $\gamma'$ to $\gamma$ in the exponents on the right-hand sides of \eqref{QS2decomp1} and \eqref{QS2decomp2} at the cost of multiplying by $d^{-\eta_1/2}$.  So using this, \eqref{QS2decomp1} and \eqref{QS2decomp2} in \eqref{QS2decomp}, we get
\begin{align*}&Q_{S,2,\delta,\eta_0}(s,t,x,t',x')\\
&\le c_4(K,R)C_{\ref{uglobalbound}}\Bigl\{d^{2-\eta_1}\Bigl[\bar\delta_N^{[\gamma\tilde\gamma_m-3/2]\wedge0}+a_n^{2\beta\gamma}\bar\delta_N^{[\gamma-3/2]\wedge0}\Bigr]\\
&\phantom{\le c(K,R)C_{\ref{uglobalbound}}\Bigl\{}+(d\wedge\sqrt\delta)^{2-\eta_1}\delta^{-3/2}\Bigl[\bar d_N^{2\gamma\tilde\gamma_m}+a_n^{2\beta\gamma}\bar d_N^{2\gamma}\Bigr]\Bigr\}.
\end{align*}
The result follows from the definition of $C_{\ref{uglobalbound}}$ and the identity 
\begin{equation}\label{gatruncid}[\gamma(\gamma_m\wedge 2)-3/2]\wedge0=(\gamma\gamma_m-3/2)\wedge 0
\end{equation} 
(use $\gamma>3/4$ here).
\gdm

\begin{lemma}\label{QTbound}
Let $0\le m\le \om+1$ and assume $(P_m)$.  For any $K\in\N^{\ge K_1}, R>2, n\in\N, \ve_0\in(0,1)$, and $\beta\in [0,1/2]$ there is a $c_{\ref{QTbound}}(K)$ and $N_{\ref{QTbound}}=N_{\ref{QTbound}}(m,n,R,\ve_0,K,\beta)(\omega)\in\N$ a.s. such that for any $\eta_1\in(R^{-1},1/2), \delta\in[{a_n},1]$, $N\in\N$, and $(t,x)\in\Rp\times\Re$, 
on
\begin{equation}\label {omegacond4}
\{\omega:(t,x)\in Z(N,n,K,\beta), N\ge N_{\ref{QTbound}}\},
\end{equation}
\begin{align}
\nonumber&Q_{T,\delta}(s,s',t',x')\\
\nonumber&\le c_{\ref{QTbound}}(K)[a_n^{-2\ve_0}+2^{4N_{\ref{QTbound}}}]|s'-s|^{1-{\eta_1\over 2}}\Bigl[{\bar{\delta}_N}^{(\gamma\gamma_m-3/2)\wedge 0}+a_n^{2\beta\gamma}{\bar{\delta}_N}^{(\gamma-3/2)\wedge0}\\
\nonumber&\phantom{\le C_{\ref{QS2bound}}(K)[a_n^{-2\ve_0}+2^{4N_3}]\Bigl[}+1(\delta<{\bar d}_N^2)\delta^{-3/2}[{\bar{d}_N}^{2\gamma\tilde\gamma_m}+a_n^{2\beta\gamma}{{\bar d}_N}^{2\gamma}]\Bigr]\\
\label{varcond}&\phantom{\le C_{\ref{QS2bound}}(K)[a_n^{-2\ve_0}+2^{4N_3}]\Bigl[}\hbox{ for all }s\le t\le t', s'\le t'\le T_K,\hbox{ and } |x'|\le K+1.
\end{align}
Here $d=d((t,x),(t',x'))$, $\bar d_N=d\vee 2^{-N}$ and $\bar\delta_N=\delta\vee {\bar d_N}^2$.  Moreover $N_{\ref{QTbound}}$ is stochastically bounded uniformly in $(n,\beta)$.
\end{lemma}
\paragraph{Proof.} Let $\xi=1-(2R)^{-1}$ and define $N_{\ref{QTbound}}=N_1(m,n,\xi,\ve_0,K,\beta)$ so that the last statement is immediate from $(P_m)$.  We may assume $s\vee s'\equiv\bar s\ge\delta$, or the left-hand side is $0$.  Let $\underline s=s\wedge s'$.  We again use Lemma~\ref{uglobalbound} to bound $|u(r,y)|$ in the integrand defining $Q_{T,\delta}$ and the maximum with $\sqrt{a_n}$ can be ignored as it is less than $\sqrt{t'-r}$ in the calculation below.  So for $\omega$ as in \eqref{omegacond4} and $s,t,s',t',x'$ as in \eqref{varcond}, we have (note that $r\le \bar s\le t'\le T_K$ so that Lemma~\ref{uglobalbound} applies)
\begin{align*}
&Q_{T,\delta}(s,s',t',x')\\
&\le C_{\ref{uglobalbound}}\int_{(\underline s-\delta)^+}^{\bar s-\delta}\int p'_{t'-r}(y-x')^2e^{2R_1K}e^{2(R_1+1)|y-x|}[2^{-N}\vee(\sqrt{t'-r}+|y-x|)]^{2\gamma \xi}\\
&\phantom{\le c_{\ref{uglobalbound}} \int_{(\underline s-\delta)^+}^{\bar s-\delta}\int}\times \Bigl[(2^{-N}\vee(\sqrt{t'-r}+|y-x|))^{\tilde\gamma_m-1}+a_n^{\beta}\Bigr]^{2\gamma}dydr.
\end{align*}
Use Lemma~\ref{pt'triv}, the inequality
\begin{equation}
2^{-N}\vee(\sqrt{t'-r}+|y-x|)\le \Bigl(2^{-N}\vee|x-x'|\Bigr)+\sqrt{t'-r}+|y-x'|\le \bar d_N+\sqrt{t'-r}+|y-x'|,
\end{equation}
and $e^{2(R_1+1)|y-x|}\le c_0(K)e^{2(R_1+1)|y-x'|}$ to bound the above by
\begin{align}
\nonumber&c_1(K)C_{\ref{uglobalbound}}\int_{(\underline s-\delta)^+}^{\bar s-\delta}\int (t'-r)^{-1}p_{2(t'-r)}(z)^2e^{2(R_1+1)|z|}[\bar d_N^{2\gamma\xi}+(t'-r)^{\gamma\xi}+
|z|^{2\gamma\xi}]\\
\nonumber&\phantom{c_1(K)C_{\ref{uglobalbound}}\int_{(\underline s-\delta)^+}^{\bar s-\delta}\int }\times\Bigl[(\bar d_N^{2\gamma(\tilde\gamma_m-1)}+(t'-r)^{\gamma(\tilde\gamma_m-1)}+|z|^{2\gamma(\tilde\gamma_m-1)})+a_n^{2\beta\gamma}\Bigr]dzdr\\
\nonumber&\le c_2(K)C_{\ref{uglobalbound}}\int_{(\underline s-\delta)^+}^{\bar s-\delta}(t'-r)^{-3/2}[\bar d_N^{2\gamma\xi}+(t'-r)^{\gamma\xi}]\Bigl[(\bar d_N^{2\gamma(\tilde\gamma_m-1)}+(t'-r)^{\gamma(\tilde\gamma_m-1)})+a_n^{2\beta\gamma}\Bigr]\,dr\\
\nonumber&\le c_2(K)C_{\ref{uglobalbound}}\Bigl\{\int_{(\underline s-\delta)^+}^{\bar s-\delta}1(r\le t'-\bar d_N^2)\Bigl[(t'-r)^{\gamma(\tilde\gamma_m+\xi-1)-3/2}+a_n^{2\beta\gamma}(t'-r)^{\gamma\xi-3/2}\Bigr]dr\\
\nonumber&\phantom{\le c_2(K)C_{\ref{uglobalbound}}}+\int_{(\underline s-\delta)^+}^{\bar s-\delta}1(r> t'-\bar d_N^2)(t'-r)^{-3/2}dr\Bigl[\bar d_N^{2\gamma(\tilde\gamma_m+\xi-1)}+a_n^{2\beta\gamma}\bar d_N^{2\gamma\xi}\Bigr]\Bigr\}\\
\label{QTJDecomp}&\equiv c_2(K)C_{\ref{uglobalbound}}\Bigl\{J_1+J_2\Bigr\}.
\end{align}

Now
\begin{align*}
\int_{(\underline s-\delta)^+}^{\bar s-\delta}1(r> t'-\bar d_N^2)(t'-r)^{-3/2}dr&\le 1(\delta<\bar d_N^2)[(t'-\bar s+\delta)^{-3/2}|s'-s|\wedge 2(t'-\bar s+\delta)^{-1/2}]\\
&\le 1(\delta<\bar d_N^2)2\delta^{-3/2}(|s'-s|\wedge\delta),
\end{align*}
and so
\begin{align}
\label{J2expa} J_2&\le 1(\delta<\bar d_N^2)2\delta^{-3/2}(|s'-s|\wedge \delta)\bar d_N^{(-2\gamma(1-\xi))}\Bigl[\bar d_N^{2\gamma\tilde\gamma_m}+a_n^{2\beta\gamma}\bar d_N^{2\gamma}\Bigr]\\
\label{J2exp}&\le c_3(K)1(\delta<\bar d_N^2)\delta^{-3/2}(|s'-s|\wedge \delta)^{1-{\eta_1\over 2}}\Bigl[\bar d_N^{2\gamma\tilde\gamma_m}+a_n^{2\beta\gamma}\bar d_N^{2\gamma}\Bigr].
\end{align}
In the last line we used $\gamma(1-\xi)\le 1-\xi\le (2R)^{-1}<\eta_1/2$ and $|s'-s|\le 2K$.  

Turning to $J_1$, let $p=\gamma(\tilde \gamma_m+\xi-1)-3/2$ or $\gamma\xi-3/2$ for $0\le m-1\le\om$.  Our bounds on $\gamma$ and $\xi$ (both are bigger than $3/4$) imply $p\in[\gamma\xi-3/2,1/2]\subset[-15/16,1/2]$.  If $p'=p\wedge0$ and $0\le \ve\le -p'$, then
\begin{align}
\nonumber I(p)&\equiv\int_{(\underline s-\delta)^+}^{\bar s-\delta}1(r\le t'-\bar d_N^2)(t'-r)^pdr\\
\nonumber &\le \int_{(\underline s-\delta)^+}^{\bar s-\delta}1(r\le t'-\bar d_N^2)\sqrt K(t'-r)^{p'}dr\\
\nonumber &\le \sqrt K \min\Bigl(|s'-s|\bar \delta_N^{p'},\int_0^{|s'-s|}u^{p'}du\Bigr)\\
\nonumber &\le 16\sqrt K|s'-s|^{p'+1}\min\Bigl(\Bigl({|s'-s|\over \bar \delta_N}\Bigr)^{-p'},1\Bigr)\quad(\hbox{use }p'\ge -15/16)\\
\nonumber &\le 16\sqrt K|s'-s|^{p'+1}\Bigl({|s'-s|\over \bar \delta_N}\Bigr)^{-p'-\ve}\\
\label{Iform}&=16\sqrt K|s'-s|^{1-\ve}\bar \delta_N^{\ve+p'}.
\end{align}

Define $q=p+\gamma(1-\xi)$, so that $q=\gamma\tilde\gamma_m-3/2$ or $\gamma-3/2$.  \hfil\break
\noindent Case 1.  $q\le 0$.\hfil\break
Then $p'=p\le 0$.  If $\ve=\gamma(1-\xi)\le (2R)^{-1}<\eta_1/2$, then $\ve+p'=q\le 0$ and so \eqref{Iform} applies, and gives
\begin{equation}
I(p)\le 16\sqrt K|s'-s|^{1-\ve}\bar \delta_N^q\le 16 K|s'-s|^{1-{\eta_1\over 2}}\bar \delta_N^q.
\end{equation}

\noindent Case 2. $q>0$\hfil\break
Then $p'=(q-\gamma(1-\xi))\wedge 0\ge -\gamma(1-\xi)$.  Let $\ve=-p'\le \gamma(1-\xi)\le (2R)^{-1}<\eta_1/2$ in \eqref{Iform} and conclude
\begin{equation}
I(p)\le 16\sqrt K |s'-s|^{1-\ve}\le 16 K|s'-s|^{1-{\eta_1\over 2}}.
\end{equation}

In either case we have shown that $I(p)\le 16 K|s'-s|^{1-{\eta_1\over 2}}\bar \delta_N^{q\wedge 0}$.  This gives
\begin{align}
\label{J1exp}J_1\le c_4(K)|s'-s|^{1-{\eta_1\over 2}}\Bigl[\bar \delta_N^{(\gamma\tilde\gamma_m-3/2)\wedge 0}+a_n^{2\beta\gamma}\bar \delta_N^{(\gamma-3/2)\wedge 0}\Bigr].
\end{align}
Put \eqref{J2exp} and \eqref{J1exp} into \eqref{QTJDecomp} and use \eqref{gatruncid} to complete the proof.
\gdm

\noindent{\bf Notation.} $d((s,t,x), (s',t',x'))=\sqrt{|s'-s|}+\sqrt{|t'-t|}+|x'-x|$.

\begin{lemma} \label{dyadicexp}Let $c_0,c_1,c_2, k_0$ be positive (universal constants), $\eta\in(0,1/2)$, and 
$\Delta:\N\times(0,1]\to\Rp$ satisfy $\Delta(n,2^{-N+1})\le k_0\Delta(n,2^{-N})$ for all $n,N\in\N$.  For $n\in\N$ and $\tau$ in a set $S$ assume $\{Y_{\tau,n}(s,t,x):(s,t,x)\in\Rp^2\times\Re\}$ is a real-valued continuous process.  Assume for each $(n,\tau)$, $K\in\N$, and $\beta\in[0,1/2]$, there is an $N_0(\omega)=N_0(n,\eta,K,\tau, \beta)(\omega)\in\N$ a.s., stochastically bounded uniformly in $(n,\tau,\beta)$, such that for any $N\in\N$, $(t,x)\in \Rp\times\Re$, and $s\le K$, if $\tilde d=d((s,t,x),(s',t',x'))\le 2^{-N}$, then
\begin{align}\label{dyadicexphyp}
&P(|Y_{\tau,n}(s,t,x)-Y_{\tau,n}(s',t',x')|>{\tilde d} ^{1-\eta}\Delta(n,2^{-N}), (t,x)\in Z(N,n,K,\beta),N\ge N_0,t'\le T_K)\\
\nonumber&\le c_0\exp(-c_1{\tilde d}^{-\eta c_2}).
\end{align}
Then there is an $N'_0=N'_0(n,\eta,K,\tau,\beta)\in\N$ a.s., also stochastically bounded uniformly in $(n,\tau,\beta)$, such that for all $N\ge N'_0(\omega)$, $(t,x)\in Z(N,n,K,\beta)(\omega)$, $\tilde d=d((s,t,x),(s',t',x'))\le 2^{-N}$, $s\le K$ and $t'\le T_K$, 
\begin{equation*}
|Y_{\tau,n}(s,t,x)-Y_{\tau,n}(s',t',x')|\le 2^7k_0^3{\tilde d}^{1-\eta}\Delta(n,2^{-N}).
\end{equation*}
\end{lemma}
\paragraph{Proof.} Let
\begin{align*}
M_{\ell,N}=M^{\tau,n,K,\beta}_{\ell,N}=\max\{&|Y_{\tau,n}((i+e)2^{-2\ell},(j+f)2^{-2\ell},(k+g)2^{-\ell})-Y_{\tau,n}(i2^{-2\ell},j2^{-2\ell},k2^{-\ell})|: \\
&\phantom{}(j2^{-2\ell},k2^{-\ell})\in Z(N,n,K+1,\beta), e,f=-4,-3,\dots 4, i2^{-2\ell}\le K+1, \\
&g=-2,-1,\dots, 2, (j+f)2^{-2\ell}\le T_{K+1}, i, j, i+e,j+f\in\Z_+, k\in \Z\},
\end{align*}
and
\[A_N=\{\omega:\exists \ell\ge N+3\hbox{ s. t. }M_{\ell,N}\ge 2^{(3-\ell)(1-\eta)}\Delta(n,2^{-N}), N\ge N_0(n,\eta,K+1,\tau,\beta)\}.\]
For $i,j,k,e,f,g$ as in the definition of $M_{\ell,N}$ and $\ell\ge N+3$, 
\[ d(((i+e)2^{-2\ell},(j+f)2^{-2\ell},(k+g)2^{-\ell}),(i2^{-2\ell},j2^{-2\ell},k2^{-\ell}))\le 2^{3-\ell}\le 2^{-N}.\]
Therefore \eqref{dyadicexphyp} implies that for some $c'_1=c'_1(\eta)>0$,
\begin{align*}P(\cup_{N'=N}^\infty A_{N'})&\le \sum_{N'=N}^\infty\sum_{\ell=N'+3}^\infty5\cdot9^2 [2^{2\ell}(K+1)+1]^2[2^{\ell+1}(K+1)+1]c_0\exp(-c_1(2^{3-\ell})^{-\eta c_2})\\
&\le c_3(K)\exp(-c'_12^{N\eta c_2}).
\end{align*}
Let
\[N_2=N_2(n,\eta,K,\tau,\beta)=\min\{N:\omega\in\cap_{N'=N}^\infty A_{N'}^c\}.\]
The above implies that
\begin{equation}\label{N1tail}
P(N_2>N)=P(\cup_{N'=N}^\infty A_{N'})\le c_3(K)\exp(-c'_12^{N\eta c_2}).
\end{equation}

Define
\[N'_0(n,\eta,K,\tau,\beta)=(N_0(n,\eta,K+1,\tau,\beta)\vee N_2(n,\eta,K,\tau,\beta))+3.\]
$N'_0$ is stochastically bounded uniformly in $(n,\tau,\beta)$ by \eqref{N1tail} and the corresponding property of $N_0$.  Assume 
\begin{equation}\label{stxcond}
N\ge N'_0,\ (t,x)\in Z(N,n,K,\beta), \ d((s,t,x)(s',t',x'))\le 2^{-N}, s\le K\hbox{ and }t'\le T_K.
\end{equation}
Define dyadic approximations by $s_\ell=\lfloor 2^{2\ell}s\rfloor 2^{-2\ell}$, $t_\ell=\lfloor 2^{2\ell}t\rfloor 2^{-2\ell}$, $x_\ell=\sgn(x)\lfloor 2^{\ell}|x|\rfloor 2^{-\ell}$, and similarly define $s'_\ell$, $t'_\ell$ and $x'_\ell$ for $(s',t',x')$.  Choose $(\hat t_0,\hat x_0)$ as in the definition of $(t,x)\in Z(N,n,K,\beta)$.  
Then $|x_\ell|\le |x|\le K$, $|x'_\ell|\le |x'|\le K+1$, $s'_\ell\vee s_\ell\le s'\vee s\le K$, $t'_\ell\vee t_\ell\le t'\vee t\le T_K$, and if $\ell\ge N$, 
\begin{align*}
d((t'_\ell,x'_\ell),(\hat t_0,\hat x_0))&\le d((t'_\ell,x'_\ell),(t',x'))+d((t',x'),(t,x))+d((t,x),(\hat t_0,\hat x_0))\\
&\le \sqrt{|t'_\ell-t'|}+|x'_\ell-x'|+2^{1-N}\\
&\le 2^{2-N}.
\end{align*}
This proves that 
\begin{equation}\label{dyadicZN1}
(t'_\ell,x'_\ell)\in Z(N-2,n,K+1,\beta)\subset Z(N-3,n,K+1,\beta)\hbox{ for all }\ell\ge N,
\end{equation}
and even more simply one gets
\begin{equation}\label{dyadicZN2}
(t_\ell,x_\ell)\in Z(N-3,n,K,\beta)\hbox{ for all }\ell\ge N.
\end{equation}
In addition, the fact that $N\ge N'_0$ implies $\omega\in A^c_{N-3}$ and $N-3\ge N_0$, which in turn implies
\begin{equation}\label{MNbound1}
M_{\ell,N-3}\le 2^{(3-\ell)(1-\eta)}\Delta(n,2^{-(N-3)})\hbox{ for all }\ell\ge N.
\end{equation}
Choose $N'\ge N$ such that $2^{-N'-1}<d((s,t,x),(s',t',x'))\equiv \tilde d\le 2^{-N'}$.  Then $|x'-x|\le 2^{-N'}$ which implies  $x'_{N'}=x_{N'}+g 2^{-N'}$ for $g\in\{-1,0,1\}$.  Similarly $s'_{N'}=s_{N'}+e2^{-2N'}$ and $t'_{N'}=t_{N'}+f 2^{-2N'}$ for $e,f\in\{-1,0,1\}$.  In addition, $s_{\ell}=s_{\ell-1}+e 4^{-\ell}$, $t_{\ell}=t_{\ell-1}+f 4^{-\ell}$, and $x_{\ell}=x_{\ell-1}+g 2^{-\ell}$ for some $e,f\in\{0,\dots, 3\}$ and $g\in\{-1,0,1\}$, and similarly for $s'_\ell$, $t'_\ell$ and $x'_\ell$.  Let $w_\ell=(s_\ell,t_\ell,x_\ell)$ and $w'_\ell=(s'_\ell,t'_\ell,x'_\ell)$.
Now use \eqref{dyadicZN1}, \eqref{dyadicZN2}, \eqref{MNbound1}, the definition of $M_{\ell,N-3}$ and the continuity of $Y_{\tau,n}$ to see that for $(s,t,x), (s',t',x')$ as in \eqref{stxcond}, 
\begin{align*}
&|Y_{\tau,n}(s,t,x)-Y_{\tau,n}(s',t',x')|\cr
&\le |Y_{\tau,n}(w'_{N'})-Y_{\tau,n}(w_{N'})|+\sum_{\ell=N'+1}^\infty |Y_{\tau,n}(w'_{\ell})-Y_{\tau,n}(w'_{\ell-1})|+|Y_{\tau,n}(w_{\ell})-Y_{\tau,n}(w_{\ell-1})|\\
&\le M_{N',N-3}+\sum_{\ell=N'+1}^\infty 2M_{\ell,N-3}\\
&\le\Bigl[ 2^{(3-N')(1-\eta)}+2\sum_{\ell=N'+1}^\infty 2^{(3-\ell)(1-\eta)}\Bigr]\Delta(n,2^{-(N-3)})\quad\hbox{(by \eqref{MNbound1})}\\
&\le (36) 2^{-N'(1-\eta)}\Delta(n,2^{-(N-3)}) \\
&\le 2^7k_0^3{\tilde d}^{1-\eta}\Delta(n,2^{-N}). 
\end{align*}
\gdm

\noindent{\bf Notation.} Introduce
\begin{align}\label{barDelta}
\bar \Delta_{u'_1}(m,n,\alpha,\ve_0,2^{-N})&=a_n^{-\ve_0}\Bigl[a_n^{-3\alpha/4}2^{-N\gamma\tilde\gamma_m}+(a_n^{\alpha/2}\vee 2^{-N})^{(\gamma_{m+1}-2)\wedge 0}\\
\nonumber&\phantom{=a_n^{-\ve_0}\Bigl[}+a_n^{-3\alpha/4+\beta\gamma}(a_n^{\alpha/2}\vee 2^{-N})^{\gamma}\Bigr].
\end{align}
We often suppress the dependence on $\ve_0$ and $\alpha$. 
\medskip

\begin{proposition}\label{Fmodulus}  Let $0\le m\le \om+1$ and assume $(P_m)$.  For any  $n\in\N$, $\eta_1\in(0,1/2)$, $\ve_0\in(0,1)$, $K\in\N^{\ge K_1}$, $\alpha\in[0,1]$, and $\beta\in[0,1/2]$, there is an $N_{\ref{Fmodulus}}=N_{\ref{Fmodulus}}(m,n,\eta_1,\ve_0,K,\alpha,\beta)(\omega)\in\N^{\ge 2}$ a.s. such that for all $N\ge N_{\ref{Fmodulus}}$, $(t,x)\in Z(N,n,K,\beta)$, $t'\le T_K$, $s\le K$,
\begin{align*}&d((s,t,x),(s',t',x'))\le 2^{-N}\hbox{ implies that }\\
&|F_{a_n^\alpha}(s,t,x)-F_{a_n^\alpha}(s',t',x')|\le 2^{-86}d((s,t,x),(s',t',x'))^{1-\eta_1}\bar \Delta_{u'_1}(m,n,\alpha,\ve_0,2^{-N}).
\end{align*}
Moreover $N_{\ref{Fmodulus}}$ is stochastically bounded, uniformly in $(n,\alpha,\beta)$.
\end{proposition}

\paragraph{Proof.} Let $R=33\eta_1^{-1}$ and choose $\eta_0\in (R^{-1},\eta_1/32)$.  Let $d=d((t,x),(t',x'))$, $\tilde d=d+\sqrt{|s'-s|}$, $\bar d_N=d\vee 2^{-N}$, $\bar\delta_{n,N}=a_n^\alpha\vee{\bar d}^2_N$ and 
\[Q_{a_n^\alpha}(s,t,x,s',t',x')=Q_{T,a_n^\alpha}(s,s',t',x')+\sum_{i=1}^2 Q_{S,i,a_n^\alpha,\eta_0}(s,t,x,t',x').\]
By Lemmas \ref{QS1bound}, \ref{QS2bound} and \ref{QTbound} there are is a $c_1(K,\eta_1)$ and $N_2=N_2(m,n,\eta_1,\ve_0,K,\beta)(\omega)$ stochastically bounded uniformly in $(n,\beta)$, such that for all $N\in\N$ and $(t,x)$, on 
\begin{equation}\label{omegacond6}
\{\omega:(t,x)\in Z(N,n,K+1,\beta),N\ge N_2\},
\end{equation}
\begin{align}
\label{Qanbound}&R_0^\gamma Q_{a_n^\alpha}(s,t,x,s',t',x')^{1/2}\\
\nonumber&\le c_1(K,\eta_1)[a_n^{-\ve_0}+2^{2N_2}]\tilde d^{1-\eta_1/2}\Bigl\{a_n^{-3\alpha/4}[\bar d_N^{\gamma\tilde\gamma_m}+a_n^{\beta\gamma}\bar d_N^{\gamma}]\\
\nonumber&\phantom{\le  c_1(K,\eta_1)[a_n^{-\ve_0}+2^{2N_2}]}+\Bigl(\sqrt{\bar\delta_{n,N}}\Bigr)^{(\gamma\gamma_m-3/2)\wedge 0}+a_n^{\beta\gamma}\Bigl(\sqrt{\bar\delta_{n,N}}\Bigr)^{(\gamma-3/2)\wedge 0}\Bigr\}\\
\nonumber&\phantom{\le c_1(K,\eta_1)[a_n^{-\ve_0}+2^{2N_2}]\tilde d^{1-\eta_1/2}\}}\hbox{for all }s\le t\le t', s'\le t'\le T_K, |x'|\le K+2.
\end{align}
Let $N_3=(33/\eta_1)[N_2+N_4(K,\eta_1)]$, where $N_4(K,\eta_1)$ is chosen large enough so that
\begin{align}
\label{auxbound1} c_1(K,\eta_1)[a_n^{-\ve_0}+2^{2N_2}]2^{-\eta_1N_3/4}
&\le  c_1(K,\eta_1)[a_n^{-\ve_0}+2^{2N_2}]2^{-8N_2-8N_4}\\
\nonumber &\le a_n^{-\ve_0}2^{-104}.
\end{align}
Let 
\begin{align*}\Delta(m,n,\bar d_N)=&2^{-100}a_n^{-\ve_0}\Bigl\{a_n^{-3\alpha/4}[\bar d_N^{\gamma\tilde\gamma_m}+a_n^{\beta\gamma}\bar d_N^{\gamma}]\\
\nonumber&\phantom{2^{-100}{\om}^{-1}a_n^{-\ve_0}\Bigl\{}+\Bigl(\sqrt{\bar\delta_{n,N}}\Bigr)^{(\gamma\gamma_m-3/2)\wedge 0}+a_n^{\beta\gamma}\Bigl(\sqrt{\bar\delta_{n,N}}\Bigr)^{(\gamma-3/2)\wedge 0}\Bigr\}.
\end{align*}
Let $N'\in\N$ and assume $\tilde d\le 2^{-N'}$.  Use \eqref{Qanbound} and \eqref{auxbound1} to see that on
\[\{\omega:(t,x)\in Z(N,n,K+1,\beta), N\ge N_3, N'\ge N_3\}\]
(which implies $|x'|\le K+2$),
\begin{align}R_0^\gamma Q_{a_n^\alpha}(s,t,x,s',t',x')^{1/2}&\le c_1(K,\eta_1)[a_n^{-\ve_0}+2^{2N_2}]2^{-\eta_1N'/4}{\tilde d}^{1-(3\eta_1/4)} 2^{100}a_n^{\ve _0}\Delta(m,n,\bar d_N)\\
\nonumber&\le {\tilde d}^{1-(3\eta_1/4)}{\Delta(m,n,\bar d_N)\over 16}\quad\hbox{ for all }s\le t\le t', 
s'\le t'\le T_K, |x'|\le K+2.
\end{align}
Combine this with \eqref{Fdecomp}, \eqref{Dbound}, the definition of $Q_{a_n^\alpha}$, and the Dubins-Schwarz theorem, to conclude that for $s\le t\le t'$, $s'\le t'$, $d((s,t,x),(s',t',x'))\le 2^{-N'}$, 
\begin{align}
\nonumber&P(|F_{a_n^\alpha}(s,t,x)-F_{a_n^\alpha}(s',t',x')|\ge d((s,t,x),(s',t',s'))^{1-\eta_1}\Delta(m,n,\bar d_N)/8,\\
\nonumber&\phantom{P(|F_{a_n^\alpha}(s,t,x)-F_{a_n^\alpha}(s',t',x')|\ge d((s,t,}(t,x)\in Z(N,n,K+1,\beta), \ N'\wedge N\ge N_3,\ t'\le T_K)\\
\nonumber&\le 2P(\sup_{u\le {\tilde d}^{2-(3\eta_1/2)}(\Delta(m,n,\bar d_N)/16)^2}|B(u)|\ge {\tilde d}^{1-\eta_1}{\Delta(m,n,\bar d_N)\over 16})\\
\nonumber&\le 2P(\sup_{u\le 1}|B(u)|\ge {\tilde d}^{-\eta_1/4})\\
\label{gaussbound}&\le c_0\exp(-{\tilde d}^{-\eta_1/2}/2).
\end{align}
Here $B(u)$ is a one-dimensional Brownian motion.

To handle $s>t$, recall that 
$F_{a_n^\alpha}(s,t,x)=F_{a_n^\alpha}(s,s\vee t,x)$.  One easily checks that 
\[\sqrt{|s\vee t-s'\vee t'|}\le \sqrt{|s-s'|}+\sqrt{|t-t'|}\] 
and hence 
$d((s,s\vee t,x),(s',s'\vee t',x'))\le 2d((s,t,x),(s',t',x'))\equiv2\tilde d$. So \eqref{gaussbound} implies that for $t\le t'$ and all $s,s',x'$, if $d((s,t,x),(s',t',x'))\le 2^{-N'-1}$, then
\begin{align}
\nonumber &P(|F_{a_n^\alpha}(s,t,x)-F_{a_n^\alpha}(s',t',x')|\ge d((s,t,x),(s',t',x'))^{1-\eta_1}\Delta(m,n,\bar d_N)/4,\\
\nonumber&\phantom{P(|F_{a_n}^\alpha(s,t,x)-F_{a_n^\alpha}}(t,x)\in Z(N,n,K+1,\beta),\ t'\le T_K,\ N'+1\ge N_3+1,\ N\ge N_3)\\
\label{gaussbound2}&\le c_0\exp(-{\tilde d}^{-\eta_1/2}/2).
\end{align}

If $(t,x)\in Z(N,n,K,\beta)$, $t'\le t$ and $d=d((t,x),(t',x'))\le 2^{-N}$, then we claim that $(t',x')\in Z(N-1,n,K+1,\beta)$.  Indeed if $(\hat t_0,\hat x_0)$ is as in the definition of $(t,x)\in Z(N,n,K,\beta)$, then $d((\hat t_0,\hat x_0),(t',x'))\le 2^{-(N-1)}$.  Also $|x'|\le K+1$, $t'\le t\le T_K$, and the claim follows.  Note also that as $d\le 2^{-N}$, we have $\bar d_N=2^{-N}$.  An elementary argument using $\gamma\gamma_k\le 2$ for $k\le m-1\le \om$ and $\gamma\tilde \gamma_k\le 2$, shows that 
\begin{equation}\label{Deltaineq}
\Delta(m,n,2^{-N})\ge 4^{-1}\Delta(m,n,2^{-(N-1)}).
\end{equation}
So, by interchanging $(t',x')$ and $(t,x)$, and replacing $N$ with $N-1$, \eqref{gaussbound2} and \eqref{Deltaineq} imply that for $t'\le t$, $\tilde d\le 2^{-N'}$ and $d\le 2^{-N}$,
\begin{align}
\nonumber &P(|F_{a_n^\alpha}(s,t,x)-F_{a_n^\alpha}(s',t',x')|\ge {\tilde d}^{1-\eta_1}\Delta(m,n,2^{-N}), \\
\nonumber&\phantom{P(|F_{a_n^\alpha}(s,t,x)-F_{a_n^\alpha}(s',t',x')|}(t,x)\in Z(N,n,K,\beta),N'\wedge N\ge N_3+1)\\
\nonumber&\le P(|F_{a_n^\alpha}(s,t,x)-F_{a_n^\alpha}(s',t',x')|\ge {\tilde d}^{1-\eta_1}\Delta(m,n,2^{-(N-1)})/4, \\
\nonumber &\phantom{\le P(|F_{a_n^\alpha}(s,t,x)-F_{a_n^\alpha}}(t',x')\in Z(N-1,n,K+1,\beta), N'\ge N_3+1, N-1\ge N_3, t\le T_K)\\
\label{gaussbound3}&\le c_0\exp(-{\tilde d}^{-\eta_1/2}/2).
\end{align}
If $N_5(m,n,\eta_1\ve_0,K,\beta)(\omega)=N_3(\omega)+1$, then $N_5$ is stochastically bounded, uniformly in $(n,\beta)$.  We have shown, (taking $N'=N$ in the above) that for all $(s,t,x), (s',t',x')$, if $\tilde d=d((s,t,x),(s',t',x'))\le 2^{-N}$, then 
\begin{align}\label{gaussbound4}
&P(|F_{a_n^\alpha}(s,t,x)-F_{a_n^\alpha}(s',t',x')|\ge {\tilde d}^{1-\eta_1}\Delta(m,n,2^{-N}), (t,x)\in Z(N,n,K,\beta), N\ge N_5, t'\le T_K)\\
\nonumber &\le  c_0\exp(-{\tilde d}^{-\eta_1/2}/2).
\end{align}

Now apply Lemma~\ref{dyadicexp} with $\tau=\alpha\in[0,1]$, $Y_{\tau,n}=F_{a_n^\alpha}$ and $k_0=2^2$, the latter by \eqref{Deltaineq}.  \eqref{gaussbound4} shows that \eqref{dyadicexphyp} holds with $N_0=N_5$.  (The implicit restriction $K\ge K_1$ in \eqref{gaussbound4} from Lemmas~\ref{QS1bound}-\ref{QTbound} is illusory as increasing $K$ only strengthens \eqref{gaussbound4}.)  Therefore there is an $N_{\ref{Fmodulus}}(\omega)=N_{\ref{Fmodulus}}(m,n,\eta_1,\ve_0,K,\alpha,\beta)\ge 2$, stochastically bounded uniformly in $(n,\alpha,\beta)$, such that for $N\ge N_{\ref{Fmodulus}}$, $(t,x)\in Z(N,n,K,\beta)$, if $t'\le T_K$, $s\le K$ and $\tilde d=d((s,t,x),(s',t',x'))\le 2^{-N}$, then
\begin{equation}\label{Fmod1}
|F_{a_n^\alpha}(s,t,x)-F_{a_n^\alpha}(s',t',x')|\le 2^{13}\Delta(m,n,2^{-N}){\tilde d}^{1-\eta_1}.
\end{equation}
Note that 
\begin{align}\label{Deltaformula}
\Delta(m,n,2^{-N})=&{2^{-100}}a_n^{-\ve_0}\Bigl[a_n^{-3\alpha/4}2^{-N\gamma\tilde\gamma_m}+(a_n^{\alpha/2}\vee 2^{-N})^{(\gamma\gamma_m-3/2)\wedge 0}\\
\nonumber&\phantom{{2^{-100}\over \om}a_n^{-\ve_0}\Bigl[}+a_n^{\beta\gamma}(a_n^{-3\alpha/4}2^{-N\gamma}+(a_n^{\alpha/2}\vee 2^{-N})^{\gamma-3/2})\Bigr]\\
\nonumber\le&2^{-99}a_n^{\ve_0}\Bigl[a_n^{-3\alpha/4}2^{-N\gamma\tilde\gamma_m}+(a_n^{\alpha/2}\vee 2^{-N})^{(\gamma_{m+1}-2)\wedge 0}\\
\nonumber&\phantom{2^{-99}a_n^{\ve_0}\Bigl[}+a_n^{\beta\gamma-{3\alpha\over 4}}(2^{-N}\vee a_n^{\alpha/2})^\gamma\Bigr].
\end{align}
Use this in \eqref{Fmod1} to complete the proof.
\gdm

Since $F_\delta(t,t,x)=-u'_{1,\delta}(t,x)$ (see Remark~4.2), the following Corollary is immediate.
 
\begin{corollary}\label{u1'modulus} Let $0\le m\le \om+1$ and assume $(P_m)$. Let $n,\eta_1,\ve_0,K,\alpha$ and $\beta$ be as in Proposition~\ref{Fmodulus}. For all $N\ge N_{\ref{Fmodulus}}$, $(t,x)\in Z(N,n,K,\beta)$ and $t'\le T_K$,
\begin{align*}&d((t,x),(t',x'))\le 2^{-N}\hbox{ implies that }\\
&|u'_{1,a_n^\alpha}(t,x)-u'_{1,a_n^\alpha}(t',x')|\le 2^{-85}d((t,x),(t',x'))^{1-\eta_1}\bar\Delta_{u'_1}(m,n,\alpha,\ve_0,2^{-N}).
\end{align*}
\end{corollary}

We will need to modify the bound in Lemma~\ref{QTbound} to control $|u'_{1,\delta}-u'_{1, a_n}|$.  Note that if $\delta\ge a_n$ and $s=t-\delta+a_n$ then
\begin{align}\label{Fu1}u'_{1,\delta}(t,x)={d\over dx}P_\delta(u_{(t-\delta)^+})(x)&={d\over dx}P_{t-s+a_n}(u_{(s-a_n)^+})(x)\\
\nonumber&=-F_{a_n}(s,t,x)\\
\nonumber&=-F_{a_n}(t-\delta+a_n,t,x).
\end{align}
Therefore the key will be a bound on $|F_{a_n}(s,t,x)-F_{a_n}(t,t,x)|$ in which the hypothesis (for Proposition~\ref{Fmodulus}) $\sqrt{t-s}\le 2^{-N}$ is weakened substantially.

\begin{lemma}\label{QTbound2}
Let $0\le m\le \om+1$ and assume $(P_m)$.  For any $K\in\N^{\ge K_1}, R>2, n\in\N, \ve_0\in(0,1)$, and $\beta\in [0,1/2]$ there is a $c_{\ref{QTbound2}}(K)$ and $N_{\ref{QTbound2}}=N_{\ref{QTbound2}}(m,n,R,\ve_0,K,\beta)(\omega)\in\N$ a.s. such that for any $\eta_1\in(R^{-1},1/2)$, $N\in\N$, and $(t,x)\in\Rp\times\Re$, 
on
\begin{equation}\label {omegacond5}
\{\omega:(t,x)\in Z(N,n,K,\beta), N\ge N_{\ref{QTbound2}}\},
\end{equation}
\begin{align*}
\nonumber&Q_{T,a_n}(s,t,t,x)\\
\nonumber&\le c_{\ref{QTbound2}}(K)[a_n^{-2\ve_0}+2^{4N_{\ref{QTbound2}}}]\Bigl\{|t-s|^{1-\eta_1/4}[((t-s)\vee a_n)^{\gamma\tilde\gamma_m-3/2}+a_n^{2\beta\gamma}((t-s)\vee a_n)^{\gamma-3/2}]\\
&\phantom{\le c_{\ref{QTbound2}}(K)[a_n^{-2\ve_0}+2]}+1(a_n<2^{-2N})((t-s)\wedge a_n)a_n^{-3/2}2^{N\eta_1/2}[2^{-2N\gamma\tilde\gamma_m}+a_n^{2\beta\gamma}2^{-2N\gamma}]\Bigr\}\\
&\phantom{\le c_{\ref{QTbound2}}(K)[a_n^{-2\ve_0}+2^{4N_{\ref{QTbound2}}}]\Bigl\{|t-s|^{1-\eta_1/4}[((t-s)\vee a_n)^{\gamma(\gamma_m\wedge 2)-3/2}+\sum_{k=0}^{m-1}}\hbox{ for all }s\le t.
\end{align*}
Moreover $N_{\ref{QTbound2}}$ is stochastically bounded uniformly in $(n,\beta)$. 
\end{lemma}
\paragraph{Proof.}  Let $\xi=1-(4\gamma R)^{-1}$ and define $N_{\ref{QTbound2}}=N_1(m,n,\xi(R),\ve_0,K,\beta)$ so that the last statement is immediate from $(P_m)$.  We may assume $t\ge a_n$.  By Lemma~\ref{uglobalbound} (again the maximum with $\sqrt{a_n}$ may be ignored in the calculation below) and then Lemma~\ref{pt'triv}, we get for $\omega$ as in \eqref{omegacond5} and $s\le t$,
\begin{align}
\nonumber&Q_{T,a_n}(s,t,t,x)\\
\nonumber&\le C_{\ref{uglobalbound}}\int_{(s-a_n)^+}^{t-a_n} \int p'_{t-r}(y-x)^2 e^{2 R_1K}e^{2(R_1+1)|y-x|}[2^{-N}+\sqrt{t-r}+|y-x|]^{2\gamma\xi}\\
\nonumber&\phantom{\le C_{uglobalbound}}\times\Bigl[(2^{-N}+\sqrt{t-r}+|y-x|)^{\tilde\gamma_m-1}+a_n^{\beta}\Bigr]^{2\gamma} dydr\\
\nonumber&\le c_1(K)C_{\ref{uglobalbound}}\int_{(s-a_n)^+}^{t-a_n} (t-r)^{-1}\int p_{2(t-r)}(z)^2e^{2(R_1+1)|z|}[2^{-2N\gamma\xi}+(t-r)^{\gamma\xi}+|z|^{2\gamma\xi}]\\
\nonumber&\phantom{\le c_1(K)C_{\ref{uglobalbound}}\int }\times \Bigl[(2^{-2N\gamma(\tilde\gamma_m-1)}+(t-r)^{\gamma(\tilde\gamma_m-1)}+|z|^{2\gamma(\tilde\gamma_m-1)}+a_n^{2\beta\gamma}\Bigr]dzdr\\
\nonumber&\le c_2(K)C_{\ref{uglobalbound}}\int_{(s-a_n)^+}^{t-a_n}(t-r)^{-3/2}[2^{-2N\gamma\xi}+(t-r)^{\gamma\xi}]\\
\nonumber&\phantom{\le c_2(K)C_{\ref{uglobalbound}}\int_{(s-a_n)^+}^{t-a_n}}\times[(2^{-2N\gamma(\tilde\gamma_m-1)}+(t-r)^{\gamma(\tilde\gamma_m-1)})+a_n^{2\beta\gamma}]dr\\
\nonumber&\le c_3(K)C_{\ref{uglobalbound}}\Bigl[\int_{(s-a_n)^+}^{t-a_n}1(r<t-2^{-2N})[(t-r)^{\gamma(\tilde\gamma_m+\xi-1)-3/2}+a_n^{2\beta\gamma}(t-r)^{\gamma\xi-3/2}]dr\\
\nonumber&\phantom{le c_3(K)C_{\ref{uglobalbound}}}+\int_{(s-a_n)^+}^{t-a_n}1(r\ge t-2^{-2N})(t-r)^{-3/2}dr2^{-2N\gamma\xi}[2^{-2N\gamma(\tilde\gamma_m-1)}+a_n^{2\beta\gamma}]\Bigr]\\
\label{Jdecomp2}&\equiv c_3(K)C_{\ref{uglobalbound}}[J_1+J_2].
\end{align}
As in the derivation of \eqref{J2expa}, now with $\bar d_N=2^{-N}$, $\delta=a_n$ and $s'=t$, we get
\begin{equation}\label{J2bound2}J_2\le 1(a_n<2^{-2N})2(a_n\wedge(t-s))a_n^{-3/2}2^{2N\gamma(1-\xi)}[2^{-2N\gamma\tilde\gamma_m}+a_n^{2\beta\gamma}2^{-2N\gamma}].
\end{equation}

For $J_1$, let $p=\gamma(\tilde\gamma_m+\xi-1)-3/2$ or $p=\gamma\xi-3/2$.  Our choice of $\xi$ and $R$ implies $p\in[-15/16,1/2]$ and so, considering $p\ge 0$ and $p<0$ separately, we arrive at
\begin{align*}\int_{(s-a_n)^+}^{t-a_n}(t-r)^pdr\le 16[(t-s+a_n)^{p+1}-a_n^{p+1}]&\le 16(2^{p^+})(t-s)((t-s)\vee a_n)^p\\
&\le 24(t-s)((t-s)\vee a_n)^p.
\end{align*}
Therefore
\begin{align}
\nonumber J_1&\le 24(t-s)\Bigl[((t-s)\vee a_n)^{\gamma(\tilde\gamma_m+\xi-1)-3/2}+a_n^{2\beta\gamma}((t-s)\vee a_n)^{\gamma\xi-3/2}\Bigr]\\
\label{J1bound2}&\le 24(t-s)^{1-\gamma(1-\xi)}\Bigl[((t-s)\vee a_n)^{\gamma\tilde\gamma_m-3/2}+a_n^{2\beta\gamma}((t-s)\vee a_n)^{\gamma-3/2}\Bigr].
\end{align}
Put \eqref{J2bound2} and \eqref{J1bound2} into \eqref{Jdecomp2}, noting that $\gamma(1-\xi(R))=(4R)^{-1}< \eta_1/4$ and $t-s\le K$, to complete the proof.
\gdm

\begin{proposition}\label{Fmodulus2} Let $0\le m\le \om+1$ and assume $(P_m)$.  For any  $n\in\N$, $\eta_1\in(0,1/2)$, $\ve_0\in(0,1)$, $K\in\N^{\ge K_1}$, and $\beta\in[0,1/2]$, there is an $N_{\ref{Fmodulus2}}=N_{\ref{Fmodulus2}}(m,n,\eta_1,\ve_0,K,\beta)(\omega)\in\N$ a.s. such that for all $N\ge N_{\ref{Fmodulus2}}$, $(t,x)\in Z(N,n,K,\beta)$, $s\le t \hbox{ and }\sqrt{t-s}\le N^{-4/\eta_1}\hbox{ implies that }$\\
\begin{align*}
&|F_{a_n}(s,t,x)-F_{a_n}(t,t,x)|\\
&\le 2^{-81}a_n^{-\ve_0}\Bigl\{2^{-N(1-\eta_1)}(a_n^{1/2}\vee 2^{-N})^{(\gamma_{m+1}-2)\wedge 0}\\
&\phantom{\le 2^{-86}{(4+\om)}^{-1}}+2^{N\eta_1}a_n^{-1/4}\Bigl({2^{-N}\over \sqrt{a_n}}+1\Bigr)\Bigl(2^{-N\gamma\tilde\gamma_m}+a_n^{\beta\gamma}(\sqrt{a_n}\vee 2^{-N})^{\gamma}\Bigr)\\
&\phantom{\le 2^{-86}{(4+\om)}^{-1}}+(t-s)^{(1-\eta_1)/2}\Bigl((\sqrt{t-s}\vee\sqrt{a_n})^{\gamma\tilde\gamma_m-{3\over 2}}+a_n^{\beta\gamma}(\sqrt{t-s}\vee\sqrt{a_n})^{\gamma-{3\over 2}}\Bigr)\Bigr\}.
\end{align*}
Moreover $N_{\ref{Fmodulus2}}$ is stochastically bounded, uniformly in $(n,\beta)$.
\end{proposition}
\paragraph{Proof.} Apply Lemma~\ref{QTbound2} with $R=2/\eta_1$ so that on 
\[\{\omega:(t,x)\in Z(N,n,K,\beta),N\ge N_{\ref{QTbound2}}(m,n,2/\eta_1,\ve_0,K,\beta)\},\]
for $s\le t$,
\begin{align}
\nonumber&R_0^\gamma Q_{T,a_n}(s,t,t,x)^{1/2}\\
\label{sqfnbounda}&\le c_1(K)R_0^\gamma[a_n^{-\ve_0}+2^{2N_{\ref{QTbound2}}}]\Bigl\{(\sqrt{t-s})^{\eta_1/4}(\sqrt{t-s})^{1-{\eta_1\over 2}}[(\sqrt{t-s}\vee\sqrt{a_n})^{\gamma\tilde\gamma_m-{3\over 2}}\\
\nonumber&\phantom{\le c_1(K)R_0^\gamma[a_n^{-\ve_0}+2^{2N_{\ref{QTbound2}}}]\Bigl\{(t-s)^{\eta_1/8}(\sqrt{t-s})^{1-{\eta_1\over 2}}}+a_n^{\beta\gamma}(\sqrt{t-s}\vee\sqrt{a_n})^{\gamma-{3\over 2}}]\\
\nonumber&\phantom{\le c_1(K)R_0^\gamma[a_n^{-\ve_0}+2^{2N_{\ref{QTbound2}}}]\Bigl\{}+2^{-N\eta_1/4}a_n^{-1/4}2^{N\eta_1/2}[2^{-N\gamma\tilde\gamma_m}+a_n^{\beta\gamma}2^{-N\gamma}]\Bigr\}.
\end{align}
Let $N_2(m,n,\eta_1,\ve_0,K,\beta)(\omega)={8\over \eta_1}[N_{\ref{QTbound2}}+N_0(K)]$, where $N_0(K)\in\N$ is chosen large enough so that
\begin{align} \label{sqfnboundb}c_1(K)R_0^\gamma[a_n^{-\ve_0}+2^{2N_{\ref{QTbound2}}}]2^{-{\eta_1\over 4}N_2}&\le c_1(K)R_0^\gamma[a_n^{-\ve_0}+2^{2N_{\ref{QTbound2}}}]2^{-2N_{\ref{QTbound2}}-2N_0(K)}\\
\nonumber&\le 2^{-100}a_n^{-\ve_0}.
\end{align}
It follows from \eqref{sqfnbounda} and \eqref{sqfnboundb} that for $N\ge N_2$, $(t,x)\in Z(N,n,K,\beta)$, $s\le t$, and $\sqrt{t-s}\le 2^{-N_2}$, 
\begin{align*}
&R_0^\gamma Q_{T,a_n}(s,t,t,x)^{1/2}\\
&\le 2^{-100}a_n^{-\ve_0}\Bigl\{(\sqrt{t-s})^{1-{\eta_1\over 2}}[(\sqrt{t-s}\vee\sqrt{a_n})^{\gamma\tilde\gamma_m-{3\over 2}}+a_n^{\beta\gamma}(\sqrt{t-s}\vee\sqrt{a_n})^{\gamma-{3\over 2}}]\\
&\phantom{\le 2^{-100}a_n^{-\ve_0}\Bigl\{}+a_n^{-1/4}2^{N\eta_1/2}[2^{-N\gamma\tilde\gamma_m}+a_n^{\beta\gamma}2^{-N\gamma}]\Bigr\}\\
&\equiv(\sqrt{t-s})^{1-{\eta_1\over 2}}\Delta_1(m,n,\sqrt{t-s}\vee\sqrt{a_n})+2^{N\eta_1/2}\Delta_2(m,n,2^{-N}).
\end{align*}
Combine this with \eqref{Dbound}, \eqref{Fdecomp} (now with $t'=t=s'$, $x=x'$, so the second integral there is $0$) and the Dubins-Schwarz theorem to see that if $B(\cdot)$ is a standard $1$-dimensional Brownian motion, then
\begin{align}\label{Gaussy1}
&P(|F_{a_n}(s,t,x)-F_{a_n}(t,t,x)|\ge (\sqrt{t-s})^{1-\eta_1}\Delta_1(m,n,\sqrt{t-s}\vee\sqrt{a_n})+2^{N\eta_1}\Delta_2(m,n,2^{-N}),\\
\nonumber&\phantom{P(|F_{a_n}(s,t,x)-F_{a_n}(t,t,x)|\ge (\sqrt{t-s})}(t,x)\in Z(N,n,K,\beta),N\ge N_2,\sqrt{t-s}\le 2^{-N_2})\\
\nonumber&\le P(\sup_{u\le 1}|B(u)|\ge (\sqrt{t-s})^{-\eta_1/2}\wedge 2^{N\eta_1/2})1(t-s\le 1)\\
\nonumber&\le c_0\exp\Bigl\{-{1\over 2}[(t-s)^{-\eta_1/2}\wedge 2^{N\eta_1}]\Bigr\}.
\end{align}

Let $\ell_N=2^{2(N+3)}N^{-8/\eta_1}$ and set
\begin{align*}
M_N(\omega)=&\max\Bigl\{{|F_{a_n}(i2^{-2(N+2)},j2^{-2(N+2)},k2^{-(N+2)})-F_{a_n}(j2^{-2(N+2)},j2^{-2(N+2)},k2^{-(N+2)})|\over (\sqrt{j-i}2^{-(N+2)})^{1-\eta_1}\Delta_1(m,n,(\sqrt{j-i}2^{-(N+2)})\vee\sqrt{a_n})+2^{N\eta_1}\Delta_2(m,n,2^{-N})}:\\
&\phantom{\max\Bigl\{AAA}0\le j-i\le \ell_N,(j2^{-2(N+2)},k2^{-(N+2)})\in Z(N,n,K,\beta),i,j\in\Z_+,k\in\Z\Bigr\}.
\end{align*}
If $N_3=2^{N_2}$, then 
\begin{align}\label{N3fact}
N\ge N_3\Rightarrow N^{-4/\eta_1}\le 2^{-N_2-1}\Rightarrow \sqrt{\ell_N}2^{-N-2}
 =  2N^{-4/\eta_1}\le 2^{-N_2}.
\end{align}
The fact that $M_N=0$ if $\ell_N<1$, \eqref{Gaussy1}, and \eqref{N3fact} imply
\begin{align*}
&P(M_N\ge 1,N\ge N_3)\\
&\le (K+1)^22^{4(N+2)}(2K+1)2^{N+2}c_0\exp\Bigl\{-{1\over 2}((\ell_N2^{-2(N+2)})^{-\eta_1/2}\wedge 2^{N\eta_1})\Bigr\}1(\ell_N\ge 1)\\
&\le c_1K^32^{5N}\exp\Bigl\{-{1\over 2}((\sqrt{\ell_N}2^{-N})\vee 2^{-N})^{-\eta_1}\Bigr\}1(\ell_N\ge 1)\\
&\le c_1K^32^{5N}\exp\Bigl\{-2^{-5/2}N^4\Bigr\}\quad\hbox{(recall $\eta_1<1/2$).}
\end{align*}
If $A_N=\{M_N\ge 1, N\ge N_3\}$ and
\[N_4=N_4(m,n,\eta_1,\ve_0,K,\beta)(\omega)=\min\{N:\omega\in\cap_{N'=N}^\infty A_{N'}^c\},\]
then
\begin{align}
\nonumber P(N_4>N)=P(\cup_{N'=N}^\infty A_{N'})&\le c_1K^3\sum_{N'=N}^\infty 2^{5N'}\exp\Bigl\{-2^{-5/2}(N')^4\Bigr\}\\
\label{N4bounda}&\le c_2(K)\exp(-N^4/6).
\end{align}
Let $N_5(\eta_1)$ be large enough so that 
\begin{equation}\label{N5bounda}
N\ge N_5\Rightarrow 2^{1-N}\le N^{-4/\eta_1}.
\end{equation}

Define 
\[N_{\ref{Fmodulus2}}(m,n,\eta_1,\ve_0,K,\beta)=(N_{\ref{Fmodulus}}\vee N_3\vee N_4 \vee N_5) +2.\]
It follows from Proposition~\ref{Fmodulus}, \eqref{N4bounda}, and the definition of $N_3$ that 
$N_{\ref{Fmodulus2}}$ is stochastically bounded uniformly in $(n,\beta)$.  Assume 
\begin{equation}\label{casehyps}
N\ge N_{\ref{Fmodulus2}}, (t,x)\in Z(N,n,K,\beta), s\le t\hbox{ and }\sqrt{t-s}\le N^{-4/\eta_1}.
\end{equation}
{\bf Case 1.}  $\sqrt{t-s}\ge 2^{1-N}$.

\noindent The condition $N\ge N_{\ref{Fmodulus2}}$ implies $\omega\in A_{N-2}^c$ and $N-2\ge N_3$, which in turn implies 
\begin{equation}\label{MNbounda} M_{N-2}<1.
\end{equation}
Let $s_\ell =\lfloor2^{2\ell} s\rfloor2^{-2\ell}$, $t_\ell=\lfloor 2^{2\ell} t\rfloor2^{-2\ell}$ and $x_\ell=\sgn(x)\lfloor 2^\ell |x|\rfloor 2^{-\ell}$, be the usual dyadic approximations to $s,t$ and $x$, respectively, and let $(\hat t_0,\hat x_0)$ be as in the definition of $(t,x)\in Z(N,n,K,\beta)$.  Then
\[ d((t_N,x_N),(\hat t_0,\hat x_0))\le 2^{-N}+\sqrt{t-t_N}+|x-x_N|\le 2^{2-N},\ t_N\le t\le T_K,\ |x_N|\le |x|\le K,
\]
and so 
\begin{equation}\label {dyadicZN}
(t_N,x_N)\in Z(N-2,n,K,\beta).
\end{equation}
Write
\begin{align}
\nonumber|F_{a_n}(s,t,x)-F_{a_n}(t,t,x)|&\le \Bigl[|F_{a_n}(s,t,x)-F_{a_n}(s_N,t_N,x_N)|+|F_{a_n}(t,t,x)-F_{a_n}(t_N,t_N,x_N)|\Bigr]\\
\nonumber&\quad+\Bigl[|F_{a_n}(s_N,t_N,x_N)-F_{a_n}(t_N,t_N,x_N)|\Bigr]\\
\label{Fdecompb}&\equiv T_1+T_2.
\end{align}
The fact that $(t,x)\in Z(N,n,K,\beta)$, $t_N\le t\le T_K$, $s\le t\le K$,
\[d((t,t,x),(t_N,t_N,x_N))\vee d((s,t,x),(s_N,t_N,x_N))\le 3 (2^{-N})\le 2^{-(N-2)},\] 
and $N-2\ge N_{\ref{Fmodulus}}$, allows us to use Proposition~\ref{Fmodulus} and infer that 
\begin{equation}\label{T1bounda}
T_1\le 2^{-85}2^{-(N-2)(1-\eta_1)}\bar \Delta_{u'_1}(m,n,1,\ve_0,2^{-(N-2)}).
\end{equation}
For $T_2$ we have from $N\ge N_5$, \eqref{N5bounda}, and the last part of \eqref{casehyps},
\[\sqrt{t_N-s_N}\le \sqrt{t-s}+2^{1-N}\le 2 N^{-4/\eta_1}\le \sqrt{\ell_{N-2}}2^{-N}.\]
In view of \eqref{dyadicZN} and \eqref{MNbounda}, this implies
\begin{align}\nonumber T_2&\le M_{N-2}\Bigl[\sqrt{t_N-s_N}^{1-\eta_1}\Delta_1(m,n,\sqrt{t_N-s_N}\vee\sqrt{a_n})+2^{(N-2)\eta_1}\Delta_2(m,n,2^{-(N-2)})\Bigl]\\
\label{T2bounda}&\le \Bigl[\sqrt{t_N-s_N}^{1-\eta_1}\Delta_1(m,n,\sqrt{t_N-s_N}\vee\sqrt{a_n})+2^{(N-2)\eta_1}\Delta_2(m,n,2^{-(N-2)})\Bigl].
\end{align}
As $t-s\ge 2^{2-2N}$ (recall this defines Case 1), we have
\[t-s\le(t-t_N)+(t_N-s_N)\le 2^{-2N}+(t_N-s_N)\le {1\over 4}(t-s)+(t_N-s_N),\]
and so 
\begin{equation}\label{tNone}
t_N-s_N\ge {1\over 2}(t-s).
\end{equation}
More simply,
\begin{equation}\label{tNtwo}
t_N-s_N\le t-s+2^{1-2N}\le 2(t-s).
\end{equation}
Use \eqref{tNone} and \eqref{tNtwo} in \eqref{T2bounda} and then combine the result with \eqref{T1bounda} and \eqref{Fdecompb} to conclude that
\begin{align}
\nonumber &|F_{a_n}(s,t,x)-F_{a_n}(t,t,,x)|\\
\nonumber&\le 2^{-85}2^{-(N-2)(1-\eta_1)}a_n^{-\ve_0}\Bigl[(a_n^{1/2}\vee 2^{-(N-2)})^{(\gamma_{m+1}-2)\wedge 0}\\
\nonumber&\phantom{\le 2^{-85}(4+\om)^{-1}2^{-(N-2)(1-\eta_1)}a_n^{-\ve_0}}+a_n^{-3/4}\Bigl[2^{-(N-2)\gamma\tilde\gamma_m}+a_n^{\beta\gamma}(a_n^{1/2}\vee 2^{-(N-2)})^{\gamma}\Bigr]\Bigr]\\
&\label {Bounda}\ +2^{-99}a_n^{-\ve_0}(\sqrt{t-s})^{1-\eta_1}\Bigl[(\sqrt{t-s}\vee\sqrt{a_n})^{\gamma\tilde\gamma_m-{3\over 2}}+a_n^{\beta\gamma}(\sqrt{t-s}\vee\sqrt{a_n})^{\gamma-{3\over 2}}\Bigr]\\
\nonumber&\ +2^{-100}a_n^{-\ve_0-{1\over 4}}2^{(N-2)\eta_1}\Bigl[2^{-(N-2)\gamma\tilde\gamma_m}+a_n^{\beta\gamma}2^{-(N-2)\gamma}\Bigr].
\end{align}
Next use 
\[2^{-N(1-\eta_1)}a_n^{-3/4}+2^{N\eta_1}a_n^{-1/4}=a_n^{-1/4}2^{N\eta_1}\Bigl[{2^{-N}\over \sqrt{a_n}}+1\Bigr]\]
to combine the first and third terms in \eqref{Bounda} and conclude, after a bit of arithmetic, that
\begin{align}
\nonumber &|F_{a_n}(s,t,x)-F_{a_n}(t,t,x)|\\
\nonumber&\le 2^{-81}a_n^{-\ve_0}\Bigl\{\Bigl[a_n^{-1/4}2^{N\eta_1}\Bigl({2^{-N}\over \sqrt{a_n}}+1\Bigr)\Bigl(2^{-N\gamma\tilde\gamma_m}+a_n^{\beta\gamma}(a_n^{1/2}\vee 2^{-N})^{\gamma}\Bigl)\\
\nonumber&\phantom{\le 2^{-81}(4+\om)^{-1}a_n^{-\ve_0}\Bigl\{\Bigl[a_n}+2^{-N(1-\eta_1)}(a_n^{1/2}\vee 2^{-N})^{(\gamma_{m+1}-2)\wedge 0}\Bigl]\\
\label {Fboundb}&\phantom{\le 2^{-81}(4+\om)^{-1}a_n^{-\ve_0}\Bigl\{}+(\sqrt{t-s})^{1-\eta_1}\Bigl((\sqrt{t-s}\vee\sqrt{a_n})^{\gamma\tilde\gamma_m-{3\over 2}}+a_n^{\beta\gamma}(\sqrt{t-s}\vee\sqrt{a_n})^{\gamma-{3\over 2}}\Bigr)\Bigr\}.
\end{align}

\noindent{\bf Case 2.} $\sqrt{t-s}\le 2^{1-N}$.\hfil\break
\noindent As $(t,x)\in Z(N-1,n,K,\beta)$ (by \eqref{casehyps}), $s\le t\le K$, $N-1\ge N_{\ref{Fmodulus}}$, and $d((s,t,x),(t,t,x))\le 2^{-(N-1)}$, we may use Proposition~\ref{Fmodulus} with $\alpha=1$ to conclude
\begin{align*}
&|F_{a_n}(s,t,x)-F_{a_n}(t,t,x)|\\
&\le 2^{-86}(\sqrt{t-s})^{1-\eta_1}\bar\Delta_{u'_1}(m,n,1,\ve_0,2^{-(N-1)})\\
&\le 2^{-83}2^{-N(1-\eta_1)}a_n^{-\ve_0}\Bigl[a_n^{-3/4}2^{-N\gamma\tilde\gamma_m}+(a_n^{1/2}\vee 2^{-N})^{(\gamma_{m+1}-2)\wedge 0}+a_n^{-3/4+\beta\gamma}(a_n^{1/2}\vee 2^{-N})^{\gamma}\Bigr]\\
&\le 2^{-83}a_n^{-\ve_0}\Bigl\{a_n^{-1/4}2^{N\eta_1}\Bigl({2^{-N}\over \sqrt{a_n}}\Bigr)\Bigl(2^{-N\gamma\tilde\gamma_m}+a_n^{\beta\gamma}(a_n^{1/2}\vee 2^{-N})^{\gamma}\Bigl)+2^{-N(1-\eta_1)}(a_n^{1/2}\vee 2^{-N})^{(\gamma_{m+1}-2)\wedge 0}\Bigr\}
\end{align*}
which is bounded by the first term on the right-hand side of \eqref{Fboundb}.
\gdm

We also need an analogue of Proposition~\ref{Fmodulus} for $G_{a_n^\alpha}$.  A subset of the arguments in Lemma~\ref{G'} shows that
\begin{equation}\label{Gform}
G_\delta(s,t,x)=\int_0^{(s-\delta)^+}\int p_{(t\vee s)-r}(y-x)D(r,y)W(dr,dy)\hbox{ for all }s\hbox{ a.s. for all } (t,x),
\end{equation}
which is just the analogue of the expression for $F_\delta$, \eqref{G'form}, with $p_{t-r}$ in place of  $p'_{t-r}$.  
Although we only will need bounds on $G_{a_n^\alpha}(s,t,x)-G_{a_n^\alpha}(t,t,x)$ (and for $\sqrt{t-s}$  small as in Proposition~\ref{Fmodulus}), this seems to require bounds on the analogues of the three types of square functions handled in Lemmas~\ref{QS1bound}, \ref{QS2bound} and \ref{QTbound}, but now with no derivatives on the Gaussian densities. This results in some simplification and a smaller singularity in $a_n^\alpha$.  We omit the proof of the following result as the details are quite similar to those used to establish Proposition~\ref{Fmodulus}.  

\begin{proposition}\label{Gmodulus1} Let $0\le m\le \om+1$ and assume $(P_m)$.  For any $n\in\N$, $\eta_1\in(0,1/2)$, $\ve_0\in(0,1)$, $K\in\N^{\ge K_1}$, $\alpha\in[0,1]$, and $\beta\in[0,1/2]$, there is an $N_{\ref{Gmodulus1}}=N_{\ref{Gmodulus1}}(m,n,\eta_1,\ve_0,K,\alpha,\beta)\in\N$ a.s. 
such that for all $N\ge N_{\ref{Gmodulus1}}$, $(t,x)\in Z(N,n,K,\beta)$, $s\le t$ and $\sqrt{t-s}\le 2^{-N}$, 
\begin{align*}
|G_{a_n^\alpha}(s,t,x)-G_{a_n^\alpha}(t,t,x)|&\le 2^{-92}(t-s)^{{1\over 2}(1-\eta_1)}a_n^{-\ve_0}a_n^{-\alpha/4}\\
&\phantom{\le}\times\Bigl[(a_n^{\alpha/2}\vee 2^{-N})^{\gamma\tilde\gamma_m}+a_n^{\beta\gamma}(a_n^{\alpha/2}\vee 2^{-N})^{\gamma}\Bigr].
\end{align*}
\end{proposition}

We need to use our global modulus of continuity for $u'_{1,a_n^\alpha}$ (Corollary~\ref{u1'modulus}) to get a modulus for $u_{1,a_n^\alpha}$ itself.  This is of course easy for spatial increments, but a key observation is that it is possible to also use control of the spatial derivatives to get a better modulus on the temporal increments. 

\medskip

\noindent{\bf Notation.} Define
\begin{align}\label{Deltau1def}
\bar\Delta_{u_1}(m,n,\alpha,\ve_0,2^{-N})&=a_n^{-\ve_0-3\alpha/4}\Bigl[a_n^\beta a_n^{3\alpha/4}+ a_n^{\beta\gamma}(a_n^{\alpha/2}\vee 2^{-N})^{\gamma+1}\\
\nonumber&\quad+(a_n^{\alpha/2}\vee 2^{-N})^{\gamma\tilde\gamma_m+1}+1(m\ge\om)a_n^{3\alpha/4}(a_n^{\alpha/2}\vee 2^{-N})\Bigr].
\end{align}
Dependence on $\alpha$ or $\ve_0$ is often suppressed.

\noindent If $\eta>0$ let $N'_{\ref{u1modulus}}(\eta)$ be the smallest natural number such that $2^{1-N}\le N^{-{4\over \eta}}$ whenever $N\ge N'_{\ref{u1modulus}}(\eta)$.

\medskip

\begin{proposition}\label{u1modulus} Let $0\le m\le \om+1$ and assume $(P_m)$.  For any $n\in\N$, $\eta_1\in(0,1/2)$, $\ve_0,\ve_1\in(0,1)$, $K\in\N^{\ge K_1}$, $\alpha\in[0,1]$ and $\beta\in[0,1/2]$, there is an $N_{\ref{u1modulus}}=N_{\ref{u1modulus}}(m,n,\eta_1,\ve_0,K,\alpha,\beta)\in\N$ a.s.  so that for all $N\ge N_{\ref{u1modulus}}$,  $n, \alpha$ satisfying
\begin{equation}\label{cond*}
a_n\le 2^{-2(N_{\ref{Fmodulus2}}(m,n,\eta_1/2,\ve_0,K,\beta)+1)}\wedge 2^{-2(N'_{\ref{u1modulus}}(\eta_1\ve_1)+1)}, \hbox{ and }\alpha\ge \ve_1,
\end{equation}
$(t,x)\in Z(N,n,K,\beta)$, $t'\le T_K$, if $d((t,x),(t',x'))\le 2^{-N}$, then
\begin{align*}
|u_{1,a_n^\alpha}(t,x)-u_{1,a_n^\alpha}(t',x')|\le 2^{-90}d((t,x),(t',x'))^{1-\eta_1}\bar \Delta_{u_1}(m,n,\alpha,\ve_0,2^{-N}).
\end{align*}
Moreover $N_{\ref{u1modulus}}$ is stochastically bounded uniformly in $n\in \N,\alpha\in [0,1]$ and $\beta\in [0,1/2]$.
\end{proposition}

\noindent{\bf Remark.}  Although $n$ appears on both sides of \eqref{cond*}, the stochastic boundedness of $N_{\ref{Fmodulus2}}$ ensures it will hold for infinitely many $n$.  This condition becomes stronger as $\alpha$ goes to $0$ and $a_n^\alpha$ moves away from the value $a_n$ where the definition of $Z(N,n,K,\beta)$ ensures some control on $u'_{1,a_n}$. This effectively rules out $\alpha=0$ from the above conclusion. 

\paragraph{Proof.}  Let 
\begin{equation*}N''_{\ref{u1modulus}}(m,n,\eta_1,\ve_0,K,\alpha,\beta)=((2N_{\ref{Fmodulus}})(m,n,\eta_1/2,\ve_0,K+1,\alpha,\beta)\vee N_{\ref{Gmodulus1}}(m,n,\eta_1,\ve_0,K+1,\alpha,\beta))+1.
\end{equation*}
Clearly $N''_{\ref{u1modulus}}$ is stochastically bounded uniformly in $(n,\alpha,\beta)$.  Assume \eqref{cond*} and 
\begin{equation}\label{u1hyps}
N\ge N''_{\ref{u1modulus}},\ (t,x)\in Z(N,n,K,\beta),\ t'\le T_K\hbox{ and }d((t,x),(t',x'))\le 2^{-N}.
\end{equation}
As in the proof of Proposition~\ref{Fmodulus}, $(t',x')\in Z(N-1,n,K+1,\beta)$, and by interchanging $(t,x)$ with $(t',x')$, $N$ with $N-1$ and $K$ with $K+1$, in the argument below (again as in the proof of Proposition~\ref{Fmodulus}) we may assume without loss of generality that $t'\le t$. Indeed, this is the reason for having $K+1$ and adding $1$ in our definition of $N''_{\ref{u1modulus}}$.  

Recall that 
\begin{equation}\label{Gandef}
G_{a_n^\alpha}(t',t,x)=P_{t-t'+a_n^\alpha}(u((t'-a_n^\alpha,\cdot))(x)=P_{t-t'}(u_{1,a_n^\alpha}(t',\cdot))(x),
\end{equation}
and so 
\begin{align}
\nonumber |u_{1,a_n^\alpha}(t',x')-u_{1,a_n^\alpha}(t,x)|&\le  
|u_{1,a_n^\alpha}(t',x')-u_{1,a_n^\alpha}(t',x)|+|u_{1,a_n^\alpha}(t',x)-P_{t-t'}(u_{1,a_n^\alpha}(t',\cdot))(x)|\\
\nonumber&\quad+|G_{a_n^\alpha}(t',t,x)-G_{a_n^\alpha}(t,t,x)|\\
\label{u1andecomp}&\equiv T_1+T_2+T_3.
\end{align}
For $T_1$, let $(\hat t_0,\hat x_0)$ be as in the definition of $(t,x)\in Z(N,n,K,\beta)$.  For $y$ between $x$ and $x'$, $d((t',y),(t,x))\le 2^{-N}$, and also $d((\hat t_0,\hat x_0),(t,x))\le 2^{-N}$.  Therefore by Corollary~\ref{u1'modulus} (twice) with $\eta_1/2$ in place of $\eta_1$,
\begin{align}
\nonumber |u'_{1,a_n^\alpha}(t',y)|&\le |u'_{1,a_n^\alpha}(t',y)-u'_{1,a_n^\alpha}(t,x)|+|u'_{1,a_n^\alpha}(t,x)-u'_{1,a_n^\alpha}(\hat t_0,\hat x_0)|\\
\nonumber&\phantom{|u'_{1,a_n^\alpha}(t',y)|\le }+|u'_{1,a_n^\alpha}(\hat t_0,\hat x_0)-u'_{1,a_n}(\hat t_0,\hat x_0)|+a_n^\beta\\
\label{u'1bounda}&\le 2^{-84}2^{-N(1-{\eta_1\over 2})}\bar \Delta_{u'_1}(m,n,\alpha,\ve_0,2^{-N})\\
\nonumber&\phantom{|u'_{1,a_n^\alpha}(t',y)|\le }+|F_{a_n}(\hat t_0-a_n^\alpha+a_n,\hat t_0, \hat x_0)-F_{a_n}(\hat t_0,\hat t_0, \hat x_0)|+a_n^\beta.
\end{align}
We have used \eqref{Fu1} in the last line. 

We now use Proposition~\ref{Fmodulus2} to control the $F$ increment in \eqref{u'1bounda}.  Choose $N'$ so that 
\begin{equation}\label{N'an}
2^{-N'-1}\le \sqrt{a_n}\le 2^{-N'}.
\end{equation}
\eqref{cond*} implies 
$\sqrt{a_n}\le 2^{-N_{\ref{Fmodulus2}}(m,n,{\eta_1\over 2},\ve_0,K,\beta)-1}$ and so
\begin{equation}\label{N'lower}
N'\ge N_{\ref{Fmodulus2}}(m,n,\eta_1/2,\ve_0,K,\beta).
\end{equation}
In addition, \eqref{cond*} implies $2^{-N'-1}\le \sqrt{a_n}\le 2^{-N'_{\ref{u1modulus}}(\eta_1\ve_1)-1}$ and so
$N'\ge N'_{\ref{u1modulus}}$ which in turn implies
\begin{equation}\label {413cond}
a_n^{\alpha/2}\le 2^{-N'\alpha}\le 2^{-N'\ve_1}\le {N'}^{-{4\ve_1\over \eta_1\ve_1}}={N'}^{-{4\over \eta_1}}.
\end{equation}
Since 
\[|u(\hat t_0,\hat x_0)|\le a_n=a_n\wedge(\sqrt{a_n}2^{-N'})\quad\hbox{ (by \eqref{N'an}}),\]
we see that $(\hat t_0,\hat x_0)\in Z(N',n,K,\beta)$.  \eqref{N'lower} and \eqref{413cond} allow us to apply Proposition~\ref{Fmodulus2} with $N'$ in place of $N$, $(\hat t_0,\hat x_0)$ in place of $(t,x)$, $\eta_1/2$ in place of $\eta_1$, and $s=\hat t_0-a_n^\alpha+a_n$, and deduce
\begin{align*}&|F_{a_n}(\hat t_0-a_n^\alpha+a_n,\hat t_0, \hat x_0)-F_{a_n}(\hat t_0,\hat t_0, \hat x_0)|\\
&\le 2^{-78}a_n^{-\ve_0}\Bigl[\sqrt{a_n}^{(1-{\eta_1\over 2})}\sqrt{a_n}^{(\gamma_{m+1}-2)\wedge 0}+\sqrt{a_n}^{-{\eta_1+1\over 2}}\Bigl(\sqrt{a_n}^{\gamma\tilde\gamma_m}+a_n^{\beta\gamma}\sqrt{a_n}^{\gamma}\Bigr)\\
&\phantom{\le 2^{-79}(4+\om)^{-1}a_n^{-\ve_0}\Bigl[}+a_n^{{\alpha\over 2}(1-{\eta_1\over 2})}\Bigl(a_n^{{\alpha\over 2}(\gamma\tilde\gamma_m-{3\over 2})}+a_n^{\beta\gamma}a_n^{{\alpha\over 2}(\gamma-{3\over 2})}\Bigr)\Bigr].
\end{align*}
The middle term in the square brackets is bounded by the last term because $\sqrt{a_n}\le a_n^{{\alpha\over 2}}$.
Therefore
\begin{align}
\nonumber &|F_{a_n}(\hat t_0-a_n^\alpha+a_n,\hat t_0, \hat x_0)-F_{a_n}(\hat t_0,\hat t_0, \hat x_0)|\\
\label{Fincbound}&\le  2^{-77}a_n^{-\ve_0}\Bigl[\sqrt{a_n}^{(1-{\eta_1\over 2})}\sqrt{a_n}^{(\gamma_{m+1}-2)\wedge 0}+a_n^{{\alpha\over 2}(1-{\eta_1\over 2})}a_n^{{-3\alpha\over 4}}\Bigl(a_n^{{\alpha\over 2}\gamma\tilde\gamma_m}+a_n^{\beta\gamma}a_n^{{\alpha\over 2}\gamma}\Bigr)\Bigr].
\end{align}
We also have 
\[\sqrt{a_n}^{1-{\eta_1\over 2}+(\gamma_{m+1}-2)\wedge 0}\le (a_n^{{\alpha\over 2}}\vee 2^{-N})^{1-{\eta_1\over 2}+(\gamma_{m+1}-2)\wedge 0},\]
because the above exponent is positive since $\eta_1<1/2$ and $\gamma>3/4$.  Use this bound in \eqref{Fincbound} and then insert the result into \eqref{u'1bounda} to conclude that for any $y$ between $x$ and $x'$,
\begin{align}
\nn&|u'_{1,a_n^\alpha}(t',y)|\\
\nn&\le 2^{-77}a_n^{-\ve_0}\Bigl\{2^{-N(1-{\eta_1\over 2})}(a_n^{{\alpha\over 2}}\vee 2^{-N})^{(\gamma_{m+1}-2)\wedge 0}\\
\nn&\phantom{\le 2^{-78}(4+\om)^{-1}}+2^{-N(1-{\eta_1\over 2})}a_n^{-{3\alpha\over 4}}\Bigl[2^{-N\gamma\tilde\gamma_m}+a_n^{\beta\gamma}(a_n^{{\alpha\over 2}}\vee 2^{-N})^{\gamma}\Bigr]\\
\nn&\phantom{\le 2^{-78}(4+\om)^{-1}}+(a_n^{{\alpha\over 2}}\vee 2^{-N})^{1-{\eta_1\over 2}+(\gamma_{m+1}-2)\wedge 0}\\
\nn&\phantom{\le 2^{-78}(4+\om)^{-1}} +a_n^{{\alpha\over 2}(1-{\eta_1\over 2})}a_n^{{-3\alpha\over 4}}\Bigl(a_n^{{\alpha\over 2}\gamma\tilde\gamma_m}+a_n^{\beta\gamma}a_n^{{\alpha\over 2}\gamma}\Bigr)\Bigr\}+a_n^\beta\\
\nn&\le 2^{-76}a_n^{-\ve_0}(a_n^{{\alpha\over 2}}\vee 2^{-N})^{1-{\eta_1\over 2}}\Bigr[(a_n^{{\alpha\over 2}}\vee 2^{-N})^{(\gamma_{m+1}-2)\wedge 0}+a_n^{-{3\alpha\over 4}}(a_n^{{\alpha\over 2}}\vee 2^{-N})^{\gamma\tilde\gamma_m}\\
\nn&\phantom{\le 2^{-77}(4+\om)^{-1}a_n^{-\ve_0}(a_n^{{\alpha\over 2}}\vee 2^{-N})}+a_n^{-{3\alpha\over 4}}a_n^{\beta\gamma}(a_n^{\alpha\over 2}\vee 2^{-N})^{\gamma}\Bigr]+a_n^\beta\\
\label{firstcalc}&\equiv  2^{-76}\tilde\Delta_{u_1}(m,n,\alpha,\ve_0,\eta_1, a_n^{{\alpha\over 2}}\vee2^{-N})+a_n^\beta.
\end{align}
Note that $\tilde\Delta_{u_1}$ is monotone increasing 
in the $2^{-N}\vee a_n^{{\alpha\over 2}}$ variable due to the positivity of the exponents (since $\eta_1<1/2$).  
The Mean Value Theorem now shows that 
\begin{equation}\label{T1bnda}
T_1\le \Bigl[a_n^\beta+ 2^{-76}\tilde\Delta_{u_1}(m,n,\alpha,\ve_0,\eta_1, a_n^{{\alpha\over 2}}\vee2^{-N})\Bigr]|x-x'|.
\end{equation}
Recalling that $t'\le t$ and that (from \eqref{u1hyps}) $N\ge N_{\ref{Gmodulus1}}$ and $\sqrt{t-t'}\le 2^{-N}$, we may apply Proposition~\ref{Gmodulus1} and infer
\begin{equation}\label{T3bnda}
T_3\le2^{-92}(t-t')^{{1\over 2}(1-\eta_1)}a_n^{-\ve_0-\alpha/4}\Bigl[(a_n^{\alpha/2}\vee 2^{-N})^{\gamma\tilde\gamma_m}+a_n^{\beta\gamma}(a_n^{\alpha/2}\vee 2^{-N})^{\gamma}\Bigr].
\end{equation}

For $T_2$, let $\{B(s):s\ge 0\}$ be a one-dimensional Brownian motion, starting at $x$ under $P_x$. Assume first that
\begin{equation}\label{IC} |B(t-t')-x|\le 2^{-{3\over 2}N_{\ref{Fmodulus}}}.
\end{equation}
Recalling \eqref{u1hyps} and that $N_{\ref{Fmodulus}}\ge 2$, we have 
\begin{align*}
d((t',B(t-t')),(t,x))\le \sqrt{t-t'}+2^{-{3\over 2}N_{\ref{Fmodulus}}}&\le 2^{-N}+2^{-{3\over 2}N_{\ref{Fmodulus}}}\\
&\le 2^{-2N_{\ref{Fmodulus}}}+2^{-{3\over 2}N_{\ref{Fmodulus}}}\\
&\le 2^{-N_{\ref{Fmodulus}}}.
\end{align*}
Define a random $N'\in\{N_{\ref{Fmodulus}},\dots,N\}$ by
\begin{align*}
&(i)\hbox{ if }d((t',B(t-t')),(t,x))\le 2^{-N}\hbox{ then }N'=N;\\
&(ii)\hbox{ if }d((t',B(t-t')),(t,x))> 2^{-N}\hbox{ then }2^{-N'-1}< d((t',B(t-t')),(t,x))\le 2^{-N'}.
\end{align*}
In case (ii) we have $2^{-N'-1}\le 2^{-N}+|B(t-t')-x|$, and so
\begin{equation}\label{N'bound}
2^{-N'}\le 2^{1-N}+2|B(t-t')-x|,
\end{equation}
a result which is trivial in case (i).   If $y$ is between $x$ and $B(t-t')$ we may argue as in \eqref{firstcalc}, but now using $(t,x)\in Z(N',n,K,\beta)$, to see that 
\begin{align}
\label{u'1boundb}  |u'_{1,a_n^\alpha}(t',y)|&\le 2^{-76}\tilde\Delta_{u_1}(m,n,\alpha,\ve_0,\eta_1,a_n^{{\alpha\over 2}}\vee2^{-N'})+a_n^\beta.
\end{align}  
Use \eqref{N'bound} and the monotonicity of $\tilde \Delta_{u_1}$ observed above to see that
\begin{align}
\label{barDeltabnd}
&a_n^{\ve_0}\tilde\Delta_{u_1}(m,n,\alpha,\ve_0,\eta_1, a_n^{{\alpha\over 2}}\vee2^{-N'})\\
\nonumber&\le 8\Bigl[[a_n^{{\alpha\over 2}}+2^{-N}+|B(t-t')-x|]^{1-{\eta_1\over 2}}\\
\nonumber&\phantom{\le\Bigl[}\times\Bigl\{(a_n^{{\alpha\over 2}}+2^{-N}+|B(t-t')-x|)^{(\gamma_{m+1}-2)\wedge 0}+a_n^{-{3\alpha\over 4}}\Bigl[(a_n^{{\alpha\over 2}}+2^{-N}+|B(t'-t)-x|)^{\gamma\tilde\gamma_m}\\
\nonumber&\phantom{\le 8\Bigl[\times\Bigl\{}+a_n^{\beta\gamma}(a_n^{{\alpha\over 2}}+2^{-N}+|B(t'-t)-x|)^{\gamma}\Bigr]\Bigr\}\Bigr].
\end{align}
Use \eqref{barDeltabnd} in \eqref{u'1boundb} and then the Mean Value Theorem to obtain (the expectation is over $B$ alone--$N_{\ref{Fmodulus}}$ remains fixed--and we are dropping a number of small constants)
\begin{align}
\nonumber &E_x(1(|B(t-t')-x|\le 2^{-{3\over 2}N_{\ref{Fmodulus}}})|u_{1,a_n^\alpha}(t',B(t-t'))-u_{1,a_n^\alpha}(t',x)|)\\
\nonumber &\le E_0(|B(t-t')|\Bigl\{a_n^\beta+a_n^{-\ve_0}[a_n^{{\alpha\over 2}(1-{\eta_1\over 2})}+2^{-N(1-{\eta_1\over 2})}+|B(t-t')|^{1-{\eta_1\over 2}}]\Bigl[a_n^{-3\alpha/4}[a_n^{{\alpha\over 2}}+2^{-N}+|B(t-t')|]^{\gamma\tilde\gamma_m}\\
\nonumber&\phantom{le E_0(|B(t-t')|\Bigl\{a_n^\beta}+(a_n^{\alpha/2}\vee 2^{-N})^{(\gamma_{m+1}-2)\wedge 0} +a_n^{-3\alpha/4}a_n^{\beta\gamma}(a_n^{\alpha/2}+2^{-N}+|B(t-t')|)^{\gamma}\Bigr]\Bigr\}\\
\nonumber&\le c_1\sqrt{t-t'}\Bigl\{a_n^\beta+a_n^{-\ve_0}[(a_n^{{\alpha\over 2}}\vee2^{-N})^{(1-{\eta_1\over 2})}+(t-t')^{{1\over 2}(1-{\eta_1\over 2})}]\Bigl[a_n^{-3\alpha/4}(a_n^{{\alpha\over 2}}\vee2^{-N}+\sqrt{t-t'})^{\gamma\tilde\gamma_m}\\
\nonumber &\phantom{\le c_1\sqrt{t-t'}\Bigl\{a_n^\beta+a_n^{-\ve_0}}+(a_n^{\alpha/2}\vee 2^{-N})^{(\gamma_{m+1}-2)\wedge 0}+a_n^{-3\alpha/4}a_n^{\beta\gamma}(a_n^{\alpha/2}\vee2^{-N}+\sqrt{t-t'})^{\gamma}\Bigr\}\\
\nonumber&\le c_2\sqrt{t-t'}\Bigl\{a_n^\beta+a_n^{-\ve_0}(a_n^{{\alpha\over 2}}\vee2^{-N})^{(1-{\eta_1\over 2})}\Bigl[(a_n^{\alpha/2}\vee 2^{-N})^{(\gamma_{m+1}-2)\wedge 0}+a_n^{-3\alpha/4}(a_n^{\alpha/2}\vee2^{-N})^{\gamma\tilde\gamma_m}\\
\nonumber&\phantom{\le c_2\sqrt{t-t'}\Bigl\{a_n^\beta+a_n^{-\ve_0}2^{-N(1-{\eta_1\over 2})}\Bigl[}+a_n^{-3\alpha/4}a_n^{\beta\gamma}(a_n^{\alpha/2}\vee 2^{-N})^{\gamma}\Bigr]\Bigr\}\quad\hbox{(since }\sqrt{t-t'}\le 2^{-N})\\
\label{BBound1}&=c_2\sqrt{t-t'}[a_n^\beta+\tilde\Delta_{u_1}(m,n,\alpha,\ve_0,\eta_1,a_n^{{\alpha\over 2}}\vee2^{-N})].
\end{align}
To handle the complementary set to that on the left-hand side of \eqref{BBound1}, note that for $K\ge K_1$ and $t'\le T_K$,
\begin{align*}
|u_{1,a_n^\alpha}(t',y)|\le E_y(|u((t'-a_n^\alpha)^+,B(a_n^\alpha))|)
\le 2KE_y(e^{|B(a_n^\alpha)|})\le 2Ke^{1+|y|}.
\end{align*}
This and the fact that $\sqrt{t-t'}\le 2^{-2N_{\ref{Fmodulus}}}$ imply that
\begin{align}
\nonumber &E_x(1(|B(t-t')-x|>2^{-{3\over 2}N_{\ref{Fmodulus}}})|u_{1,a_n^\alpha}(t',B(t-t'))-u_{1,a_n^\alpha}(t',x)|)\\
\nonumber&\le P_0(|B(t-t')|>2^{-{3\over 2}N_{\ref{Fmodulus}}})^{1/2}8KeE_x(e^{2|B(t-t')|}+e^{2|x|})^{1/2}\\
\nonumber&\le c_3(K)P_0(|B(1)|>(t-t')^{-1/8})^{1/2}\qquad\hbox{  (since }|x|\le K\hbox{ by \eqref{u1hyps})}\\
\nonumber&\le c_4(K)(t-t')\\
\label{BBound2}&\le c_5(K)\sqrt{t-t'}\tilde\Delta_{u_1}(m,n,\alpha,\ve_0,\eta_1,a_n^{{\alpha\over 2}}\vee2^{-N}),
\end{align}
where in the last line we use 
\[\tilde\Delta_{u_1}(m,n,\alpha,\ve_0,\eta_1,a_n^{{\alpha\over 2}}\vee2^{-N})\ge (a_n^{{\alpha\over 2}}\vee 2^{-N})^{1-{\eta_1\over 2}}\ge 2^{-N}\ge \sqrt{t-t'}.\]
\eqref{BBound1} and \eqref{BBound2} imply
\begin{equation}\label{T2bnda}
T_2\le c_6(K)\sqrt{t-t'}[a_n^\beta+\tilde\Delta_{u_1}(m,n,\alpha,\ve_0,\eta_1,a_n^{{\alpha\over 2}}\vee2^{-N})].
\end{equation}

Use \eqref{T1bnda}, \eqref{T3bnda} and \eqref{T2bnda} in \eqref{u1andecomp} to conclude
\begin{align*}
\nonumber&|u_{1,a_n^\alpha}(t',x')-u_{1,a_n^\alpha}(t,x)|\\
\nonumber&\le c_7(K)d((t,x),(t',x'))[a_n^\beta+\tilde\Delta_{u_1}(m,n,\alpha,\ve_0,\eta_1,a_n^{{\alpha\over 2}}\vee2^{-N})]\\
\nonumber&\quad+2^{-92}(t'-t)^{{1\over 2}(1-\eta_1)}a_n^{-\ve_0}a_n^{\alpha/2}a_n^{-3\alpha/4}\Bigl[(a_n^{\alpha/2}\vee 2^{-N})^{\gamma\tilde\gamma_m}+a_n^{\beta\gamma}(a_n^{\alpha/2}\vee 2^{-N})^{\gamma}\Bigr].
\end{align*}
since $d((t,x),t',x'))\le 2^{-N}$, a bit of arithmetic shows the above is at most
\begin{align}
\nonumber & (c_7(K)2^{-N{\eta_1\over 2}}+2^{-92})d((t,x),(t',x'))^{1-\eta_1}a_n^{-\ve_0-3\alpha/4}\\
\label{u1prebound}&\quad\times\Bigl\{(a_n^{\alpha/2}\vee2^{-N})\Bigl[(a_n^{\alpha/2}\vee 2^{-N})^{\gamma\tilde\gamma_m}+a_n^{\beta\gamma}(a_n^{\alpha/2}\vee 2^{-N})^{\gamma}\\
\nonumber &\phantom{\quad\times\Bigl\{(a_n^{\alpha/2}\vee2^{-N})\Bigl[}+a_n^{3\alpha/4}(a_n^{\alpha/2}\vee 2^{-N})^{(\gamma_{m+1}-2)\wedge 0}\Bigr]+a_n^{3\alpha/4}a_n^\beta\Bigr\}.
\end{align}
Choose $N_1(K,\eta_1)$ so that 
\begin{equation*} 2^{-N_1\eta_1/2}c_7(K)\le 2^{-92},
\end{equation*}
and define $N_{\ref{u1modulus}}=N''_{\ref{u1modulus}}\vee N_1$, which is clearly stochastically bounded uniformly in $(n,\alpha,\beta)\in\N\times[0,1]\times[0,{1\over 2}]$.  Assume $N\ge N_{\ref{u1modulus}}$.  Note that if $m<\om$, then
\begin{align} 
\nonumber a_n^{3\alpha/4}(a_n^{\alpha/2}\vee 2^{-N})^{(\gamma_{m+1}-2)\wedge 0}&=a_n^{3\alpha/4}(a_n^{\alpha/2}\vee 2^{-N})^{\gamma\gamma_m-{3\over 2}}\\
\label{extrabnd}&\le(a_n^{\alpha/2}\vee 2^{-N})^{\gamma\gamma_m}\\
\nonumber&\le (a_n^{\alpha/2}\vee 2^{-N})^{\gamma\tilde\gamma_m},
\end{align}
and so the last term in square brackets in \eqref{u1prebound} is bounded by the first term in the same brackets.  If $m\ge \om$, the left-hand side of \eqref{extrabnd} is $a_n^{3\alpha/4}$. Therefore if $N\ge N_{\ref{u1modulus}}$ we conclude that
\begin{align*}
&|u_{1,a_n^\alpha}(t',x')-u_{1,a_n^\alpha}(t,x)|\\
&\le 2^{-90}d((t,x),(t',x'))^{1-\eta_1}a_n^{-\ve_0-3\alpha/4}\\
&\quad \times\Bigl[a_n^{\beta+3\alpha/4}+a_n^{\beta\gamma}(a_n^{\alpha/2}\vee 2^{-N})^{\gamma+1}+(a_n^{\alpha/2}\vee 2^{-N})^{\gamma\tilde\gamma_m+1}+1(m\ge \om)(a_n^{\alpha/2}\vee 2^{-N})a_n^{3\alpha/4}\Bigr]\\
&=2^{-90}d((t,x),(t',x'))^{1-\eta_1}\bar\Delta_{u_1}(m,n,\alpha,\ve_0,2^{-N}).
\end{align*}
\gdm

We also require  an analogue of the bound on increments on $u_{1,a_n^\alpha}$ (Proposition~\ref{u1modulus}) for $u_{2,a_n^\alpha}$.  
\medskip

\noindent{\bf Notation.}\begin{align*}
 \bar\Delta_{1,u_2}(m,n,\ve_0,2^{-N})=&a_n^{-\ve_0}2^{-N\gamma}\Bigl[(a_n^{1/2}\vee 2^{-N})^{\gamma(\tilde\gamma_m-1)}+a_n^{\beta\gamma}\Bigr]\\
 \bar\Delta_{2,u_2}(m,n,\alpha,\ve_0)=&a_n^{-\ve_0}\Bigl[a_n^{{\alpha\over 2}(\gamma\tilde\gamma_m-{1\over 2})}+a_n^{\beta\gamma}a_n^{{\alpha\over 2}(\gamma-{1\over 2})}\Bigr].
 \end{align*}
We often will suppress the dependence on $\ve_0$ and $\alpha$ below.

\medskip

\begin{proposition}\label{u2modulus}
Let $0\le m\le \om+1$ and assume $(P_m)$.  For any  $n\in\N$, $\eta_1\in(0,1/2)$, $\ve_0\in(0,1)$, $K\in\N^{\ge K_1}$, $\alpha\in[0,1]$, and $\beta\in[0,1/2]$, there is an $N_{\ref{u2modulus}}=N_{\ref{u2modulus}}(m,n,\eta_1,\ve_0,K,\alpha,\beta)(\omega)\in\N$ a.s. such that for all $N\ge N_{\ref{u2modulus}}$, $(t,x)\in Z(N,n,K,\beta)$, and $t'\le T_K$, 
\begin{align*}&d\equiv d((t,x),(t',x'))\le 2^{-N}\hbox{ implies that }\\
&|u_{2,a_n^\alpha}(t,x)-u_{2,a_n^\alpha}(t',x')|\le 2^{-89}\Bigl[d^{{1-\eta_1\over 2}}\bar\Delta_{1,u_2}(m,n,\ve_0,2^{-N})+d^{1-\eta_1}\bar \Delta_{2,u_2}(m,n,\alpha,\ve_0)\Bigr].
\end{align*}
Moreover $N_{\ref{u2modulus}}$ is stochastically bounded, uniformly in $(n,\alpha,\beta)$.
\end{proposition}

The proof is more straightforward than that of Proposition~\ref{u1modulus} and is given in Section~\ref{sec6} below.  

\begin{lemma}\label{missed}
For all $n\in \N$, $0\le \beta\le 1/2$ and $0<d\le 1$, 
$$a_n^{\beta\gamma}(\sqrt{a_n}\vee d)^{\gamma_1-1}\le d\vee a_n^\beta.$$
\end{lemma}
\paragraph{Proof.}
Recall that $\gamma_1-1=\gamma-{1\over 2}$ and $2\gamma-{1\over 2}>1$.

\noindent Case 1. $d\ge a_n^\beta$.\hfil\break
$a_n^{\beta\gamma}(\sqrt{a_n}\vee d)^{\gamma_1-1}\le d^{2\gamma-{1\over 2}}\le d$.

\noindent Case 2.  $d<a_n^\beta$.\hfil\break
$a_n^{\beta\gamma}(\sqrt{a_n}\vee d)^{\gamma_1-1}\le a_n^{\beta\gamma+\beta(\gamma-{1\over 2})}=a_n^{\beta(2\gamma-{1\over 2})}\le a_n^\beta$.
\gdm

We are finally ready to complete the

\paragraph{Proof of Proposition~\ref{Pminduction}.} Let $0\le m \le \om$ and assume $(P_m)$.  We must  derive $(P_{m+1})$.  Let $\ve_0\in(0,1)$, $M=\lceil{2\over \ve_0}\rceil$, $\ve_1={1\over M}\le \ve_0/2$ and set $\alpha_i=i\ve_1$ for $i=0,1,\dots,M$, so that $\alpha_i\in[\ve_1,1]$ for $i\ge 1$.  Let $n$, 
$\xi$, $K$, and $\beta$ be as in $(P_m)$ where we may assume $\xi>1/2$ without loss of generality.   Define $\eta_1=1-\xi\in (0,1/2)$, $\xi'=\xi+(1-\xi)/2\in(\xi,1)$,
\[N_2(m,n,\xi,\ve_0,K,\beta)(\omega)=\vee_{i=1}^MN_{\ref{u1modulus}}(m,n,\eta_1,\ve_0/2,K+1,\alpha_i,\beta)(\omega),\]
\[N_3(m,n,\xi,\ve_0,K,\beta)(\omega)=\vee_{i=1}^MN_{\ref{u2modulus}}(m,n,\eta_1,\ve_0/2,K+1,\alpha_i,\beta)(\omega),\]
\begin{align*}
N_4(m,n,\xi,\ve_0,K,\beta)&=\lceil{2\over 1-\xi}((N_{\ref{Fmodulus2}}(m,n,\eta_1/2,\ve_0/2,K+1,\beta)\vee 
N'_{\ref{u1modulus}}(\eta_1\ve_1))+1)\rceil\\
&\equiv \lceil{1\over 1-\xi}N_5(m,n,\eta_1,\ve_0,K,\beta)\rceil,
\end{align*}
(recall $\ve_1$ is a function of $\ve_0$) and
\begin{equation}\label{N1def}
N_1(m,n,\xi,\ve_0,K,\beta)(\omega)=(N_2\vee N_3\vee N_4(m,n,\xi,\ve_0,K,\beta))\vee N_1(0,\xi',K)+1\in\N\hbox{ a.s.}
\end{equation}
Recall that in the verification of $(P_0)$, we may take $\ve_0=0$ and $N_1=N_1(0,\xi',K)$ was independent of $n$ and $\beta$.  
Then $N_1=N_1(m,n,\xi,\ve_0,K,\beta)$
 is stochastically bounded uniformly in $(n,\beta)$ because $N_{\ref{Fmodulus2}}$, $N_{\ref{u1modulus}}$ and 
 $N_{\ref{u2modulus}}$ all are. 

Assume 
\begin{equation*}N\ge N_1,\ (t,x)\in Z(N,n,K,\beta),\ t'\le T_K,\ \hbox{ and }d\equiv d((t,x),(t',x'))\le 2^{-N}.
\end{equation*}

Suppose first that 
\begin{equation} \label{anbig}a_n>2^{-N_5(m,n,\eta_1,\ve_0,K,\beta)}.
\end{equation}
Since $N\ge N_1(0,\xi',K)$, we have by $(P_0)$, with $\ve_0=0$ and $\xi'$ in place of $\xi$, and the fact that $\tilde\gamma_{m+1}-1\le 1$, 
\begin{align*}
|u(t',x')|&\le 2^{-N\xi'}\\
&\le 2^{-N\xi'}\Bigl[(\sqrt{a_n}\vee 2^{-N})^{\tilde\gamma_{m+1}-1}\Bigr] 2^{N_5/2}\\
&\le 2^{-N(1-\xi)/2}2^{N_5/2}2^{-N\xi}\Bigl[(\sqrt{a_n}\vee 2^{-N})^{\tilde\gamma_{m+1}-1}+a_n^{\beta}\Bigr] \\
&\le 2^{-N\xi}\Bigl[(\sqrt{a_n}\vee 2^{-N})^{\tilde\gamma_{m+1}-1}+ a_n^{\beta}\Bigr],
\end{align*}
where in the last line we used $N\ge N_4\ge (1-\xi)^{-1}N_5$.  This proves $(P_{m+1})$.

So assume now that 
\begin{equation}\label{ansmall}
a_n\le2^{-N_5(m,n,\eta_1,\ve_0,K,\beta)}.
\end{equation}
Let $N'=N-1\ge N_2\vee N_3$.  Note  that $(\hat t_0,\hat x_0)$ (the point near $(t,x)$ in the definition of $Z(N,n,K,\beta)$) is in $Z(N,n,K+1,\beta)\subset Z(N',n,K+1,\beta)$ and by the triangle inequality $d((\hat t_0,\hat x_0),(t',x'))\le 2^{-N'}$. \eqref{ansmall} shows that \eqref{cond*} holds with $(\ve_0/2,K+1)$ in place  of $(\ve_0,K)$.  Therefore the inequality $N'\ge N_2$ allows us to apply the conclusion of Proposition~\ref{u1modulus}  for $\alpha=\alpha_i\ge\ve_1$, $i=1,\dots,M$ with $(\hat t_0,\hat x_0)$ in place of $(t,x)$, $\ve_0/2$ in place of $\ve_0$, and $N'$ in place of $N$.  Simpler reasoning, using $N'\ge N_3$, allows us to apply the conclusion of Proposition~\ref{u2modulus} with the same parameter values.  

Choose $i\in\{1,\dots,M\}$ so that 
\begin{align}
\label{idef} (i)& \hbox{ if }2^{-N'}>\sqrt{a_n},\hbox{ then }a_n^{{\alpha_i\over 2}}< 2^{-N'}\le a_n^{{\alpha_{i-1}\over 2}}=a_n^{{\alpha_i\over 2}}a_n^{-{\ve_1\over 2}},\\
\nonumber (ii)& \hbox{ if }2^{-N'}\le \sqrt{a_n},\hbox{ then } i=M\hbox{ and so }a_n^{{\alpha_i\over 2}}=\sqrt{a_n}\ge 2^{-N'}.
\end{align}
In either case we have
\begin{equation}\label{veebnd1}
a_n^{{\alpha_i\over 2}}\vee 2^{-N'}\le \sqrt{a_n}\vee 2^{-N'},
\end{equation}
and
\begin{equation}
\label{veebnd2}
a_n^{-{3\alpha_i\over 4}}(\sqrt{a_n}\vee 2^{-N'})^{3/2}\le a_n^{-{3\ve_1\over 4}}.
\end{equation}

Now apply Propositions~\ref{u1modulus} and \ref{u2modulus}, as described above, as well as \eqref{veebnd1}, and the facts that $\tilde\gamma_m=\gamma_m$ for $m\le\om$, $\gamma_{m+1}=\gamma\gamma_m+{1\over 2}$ and $d((\hat t_0,\hat x_0),(t',x'))\le 2^{-N'}$, to conclude
\begin{align}
\nonumber &|u(\hat t_0,\hat x_0)-u(t',x')|\\
\nonumber &\le |u_{1,a_n^{\alpha_i}}(\hat t_0,\hat x_0)-u_{1,a_n^{\alpha_i}}(t',x')|+|u_{2,a_n^{\alpha_i}}(\hat t_0,\hat x_0)-u_{2,a_n^{\alpha_i}}(t',x')|\\
\nonumber&\le 2^{-89}a_n^{-{\ve_0\over 2}}\Bigl\{2^{-N'\xi}\Bigl[a_n^\beta+a_n^{-{3\alpha_i\over 4}}a_n^{\beta\gamma} (\sqrt{a_n}\vee 2^{-N'})^{(\gamma+1)}+a_n^{-{3\alpha_i\over 4}}(\sqrt{a_n}\vee 2^{-N'})^{(\gamma\gamma_m+1)}\\
\nonumber&\phantom{\le 2^{-89}}+1(m=\om)(\sqrt{a_n}\vee 2^{-N'})\Bigr]+2^{-N'({\xi\over 2}+\gamma)}\Bigl[(\sqrt{a_n}\vee 2^{-N'})^{\gamma(\gamma_m-1)}+a_n^{\beta\gamma}\Bigr]\\
\nonumber&\phantom{\le 2^{-89}a_n^{-{\ve_0\over 2}}\Bigl\{}+ 2^{-N'\xi}\Bigl[(\sqrt{a_n}\vee 2^{-N'})^{(\gamma_{m+1}-1)}+a_n^{\beta\gamma}(\sqrt{a_n}\vee 2^{-N'})^{(\gamma-{1\over 2})}\Bigr]\Bigr\}\\
\nonumber&\le 2^{-89}a_n^{-{\ve_0\over 2}} 2^{-N'\xi}\Bigl\{a_n^\beta+a_n^{-{3\alpha_i\over 4}}(\sqrt{a_n}\vee 2^{-N'})^{{3\over 2}}a_n^{\beta\gamma}(\sqrt{a_n}\vee 2^{-N'})^{(\gamma-{1\over 2})}\\
\nn&\phantom{\le 2^{-89}a_n^{-{\ve_0\over 2}} 2^{-N'\xi}\Bigl\{}+a_n^{-{3\alpha_i\over 4}}(\sqrt{a_n}\vee 2^{-N'})^{{3\over 2}}(\sqrt{a_n}\vee 2^{-N'})^{(\gamma_{m+1}-1)}+1(m=\om)(\sqrt{a_n}\vee 2^{-N'})\\
\label{firstubound}&\phantom{\le 2^{-89}a_n^{-{\ve_0\over 2}} 2^{-N'\xi}\Bigl\{}+(\sqrt{a_n}\vee 2^{-N'})^{(\gamma-{1\over 2})}\Bigl[(\sqrt{a_n}\vee 2^{-N'})^{\gamma(\gamma_m-1)}+a_n^{\beta\gamma}\Bigr]\\
\nonumber&\phantom{\le 2^{-89}a_n^{-{\ve_0\over 2}} 2^{-N'\xi}\Bigl\{}+(\sqrt{a_n}\vee 2^{-N'})^{(\gamma_{m+1}-1)}+a_n^{\beta\gamma}(\sqrt{a_n}\vee 2^{-N'})^{(\gamma-{1\over 2})}\Bigr\}.
\end{align}
Now apply \eqref{veebnd2} and combine some duplicate terms to bound $|u(\hat t_0,\hat x_0)-u(t',x')|$ by
\begin{align*}
&2^{-87}a_n^{-{\ve_0\over 2}-{3\ve_1\over 4}}2^{-N'\xi}\Bigl[a_n^{\beta}+a_n^{\beta\gamma}(\sqrt{a_n}\vee 2^{-N'})^{\gamma-{1\over 2}}+(\sqrt{a_n}\vee 2^{-N'})^{(\gamma_{m+1}-1)}+1(m=\om)(\sqrt{a_n}\vee 2^{-N'})\Bigr].
\end{align*}
Use the fact that 
\[(\sqrt{a_n}\vee 2^{-N'})^{(\gamma_{m+1}-1)}+1(m=\om)(\sqrt{a_n}\vee 2^{-N'})\le 2(\sqrt{a_n}\vee 2^{-N'})^{((\gamma_{m+1}\wedge 2)-1)}\]
(consider $m<\om$ and $m=\om$ separately) and $\ve_1\le \ve_0/2$  in the above to derive
\begin{align}
\nonumber |u(\hat t_0,\hat x_0)-u(t',x')|&\le 2^{-86}a_n^{-\ve_0}2^{-N'\xi}\Bigl[a_n^\beta+ a_n^{\beta\gamma}(\sqrt{a_n}\vee 2^{-N'})^{(\gamma-{1\over 2})}+(\sqrt{a_n}\vee 2^{-N'})^{((\gamma_{m+1}\wedge 2)-1)}\Bigr]\\
\label{secondubound}&\le 2^{-84}a_n^{-\ve_0}2^{-N\xi}\Bigl[a_n^\beta+ a_n^{\beta\gamma}(\sqrt{a_n}\vee 2^{-N})^{(\gamma-{1\over 2})}+(\sqrt{a_n}\vee 2^{-N})^{((\gamma_{m+1}\wedge 2)-1)}\Bigr].
\end{align}

Finally combine \eqref{secondubound} and $|u(\hat t_0,\hat x_0)|\le \sqrt{a_n}2^{-N}$ to conclude
\begin{align}
\nn |u(t',x')|&\le \sqrt{a_n}2^{-N}+2^{-84}a_n^{-\ve_0}2^{-N\xi}\Bigl[\sum_{k=0}^1 a_n^{\beta\gamma^k}(2^{-N}\vee\sqrt{a_n})^{\gamma_k-1}+(2^{-N}\vee \sqrt{a_n})^{\tilde\gamma_{m+1}-1}\Bigr]\\
\label{fourthubound}&\le a_n^{-\ve_0}2^{-N\xi}[\sqrt{a_n}2^{-N(1-\xi)}+2^{-84}\Bigl[\sum_{k=0}^1a_n^{\beta\gamma^k}(2^{-N}\vee\sqrt{a_n})^{\gamma_k-1}+(2^{-N}\vee \sqrt{a_n})^{\tilde\gamma_{m+1}-1}\Bigr]\Bigr].
\end{align}
Our definition of $N_1$ (and especially $N_4$) ensures that $N(1-\xi)\ge 1$ and hence
\[\sqrt{a_n}2^{-N(1-\xi)}\le{\sqrt{a_n}\over 2}\le {a_n^\beta\over 2}
.\]
In addition by Lemma~\ref{missed}
\begin{align*}
a_n^{\beta\gamma}(\sqrt{a_n}\vee 2^{-N})^{\gamma_1-1}\le a_n^\beta\vee 2^{-N}\le a_n^\beta+(2^{-N}\vee\sqrt{a_n})^{\tilde\gamma_{m+1}-1}.
\end{align*}
Substitute the last bounds into \eqref{fourthubound} to obtain $(P_{m+1})$ and hence complete the induction.\gdm

\section{Proof of Proposition~\ref{tildeJ}}\label{sec5}
\setcounter{equation}{0}
\setcounter{theorem}{0} 
We continue to assume $b\equiv0$ in this Section.  Having established the bound $(P_{\om+1})$ in Proposition~\ref{Pminduction}, we are now free to use the conclusions of  Corollary~\ref{u1'modulus}, 
Proposition~\ref{Fmodulus2} and Proposition~\ref{u2modulus}, with $m=\om+1$, to derive local moduli of continuity for $u'_{1,a_n^\alpha}$ and $u_{2,a_n^\alpha}$.  In view of our main goal,  
Proposition~\ref{tildeJ},
it is the space modulus we will need--the time modulus was only needed to carry out the induction leading to $(P_{\om+1})$.

We fix a $K_0\in\N^{\ge K_1}$ and positive constants $\ve_0,\ve_1$ as in~\eqref{veconditions}.
For $M,n\in\N$ and $0<\beta\le {1\over2}-\ve_1$, define
\begin{equation}\label{alphadef}
\alpha=\alpha(\beta)=2(\beta+\ve_1)\in[0,1].
\end{equation}
and
\begin{align*}
U^{(1)}_{M,n,\beta}&=\inf\Bigl\{t:\hbox{there are }\ve\in[0,2^{-M}],\ |x|\le K_0+1,\ \hat x_0, \ x'\in\Re,\ s.t.\ |x-x'|\le 2^{-M},\\
&\phantom{=\inf\{t:}\ |x-\hat x_0|\le \ve, \ |u(t,\hat x_0)|\le a_n\wedge(\sqrt{a_n}\ve),\ |u'_{1,{a_n}}(t,\hat x_0)|\le a_n^\beta, \hbox{ and }\\
&\phantom{=\inf\{t:}|u'_{1,a_n^\alpha}(t,x)-u'_{1,a_n^\alpha}(t,x')|>2^{-82}a_n^{-\ve_0-{3\ve_1\over 2}}|x-x'|^{1-\ve_0}\Bigl[a_n^{-3\beta/2}(\ve\vee|x'-x|)^{2\gamma}+1\\
&\phantom{|u'_{1,a_n^\alpha}(t,x)-u'_{1,a_n^\alpha}(t,x')||x-x'|}+a_n^{\beta(\gamma-{3\over 2})}(\ve\vee|x'-x|)^{\gamma}\Bigr]\Bigr\}\wedge T_{K_0}.
\end{align*}
Define $U^{(1)}_{M,n,0}$ by the same expression (with $\beta=0$) but with the condition on $|u'_{1,a_n^\alpha}(t,\hat x_0)|$ omitted.  

$\{U_{M,n,\beta}^{(1)}<t\}$ is the projection onto $\Omega$ of a Borel$\times \cF_t$-measurable set in $K\times[0,t]\times\Omega$ where $K$ is a compact subset of $\Re^3$, and so $U_{M,n,\beta}^{(1)}$ is an $(\cF_t)$-stopping time as in IV.T52 of \cite{bib:mey66}. 

\begin{lemma}\label{U1conv}For each $n\in \N$ and $\beta$ as in \eqref{betarange}, $U^{(1)}_{M,n,\beta}\uparrow T_{K_0}$ as $M\uparrow\infty$, and in fact
\[\lim_{M\to\infty} \sup_{n,0\le \beta\le {1\over 2}-\ve_1}P(U^{(1)}_{M,n,\beta}<T_{K_0})=0.\]
\end{lemma}
\paragraph{Proof.}  By monotonicity in $M$ the first assertion is immediate from the second.  Proposition~\ref{Pminduction} allows us to apply Corollary~\ref{u1'modulus} with $m=\om+1$, $\eta_1=\ve_0$, $K=K_0+1$, and $\alpha,\beta$ as in \eqref{alphadef} and \eqref{betarange}, respectively.  Hence there is an $N_0=N_0(n,\ve_0,\ve_1,K_0+1,\beta)\in\N$ a.s., stochastically bounded uniformly in $(n,\beta)$ (as in \eqref{betarange}), and such that if
\begin{equation}\label{u'1cond}
N\ge N_0(\omega),\ (t,x)\in Z(N,n,K_0+1,\beta),\ |x-x'|\le 2^{-N},
\end{equation}
then
\begin{align*}
&|u'_{1,a_n^\alpha}(t,x)-u'_{1,a_n^\alpha}(t,x')|\\
&\le 2^{-85}|x-x'|^{1-\ve_0}a_n^{-\ve_0-{3\ve_1\over 2}}
\Bigl[a_n^{-{3\beta\over 2}}2^{-2N\gamma}+1+
a_n^{\beta(\gamma-{3\over 2})}(a_n^{\beta+\ve_1}\vee 2^{-N})^{\gamma}
\Bigr].
\end{align*}
Note that $a_n^{-{3\beta\over 2}+\beta\gamma+(\beta+\ve_1)\gamma}\le 1$ since $\gamma>3/4$, and so by the above we have
\begin{align}\label{U1bnd1}
&|u'_{1,a_n^\alpha}(t,x)-u'_{1,a_n^\alpha}(t,x')|\\
\nn&\le 2^{-84}
a_n^{-\ve_0-{3\ve_1\over 2}}|x-x'|^{1-\ve_0}\Bigl[a_n^{-{3\beta\over 2}}2^{-2N\gamma}+1+ a_n^{\beta(\gamma-{3\over 2})}2^{-N\gamma}\Bigr].
\end{align}

Let us assume $\beta>0$ for if $\beta=0$ we can just omit the bound on $|u'_{1,a_n^\alpha}(t,\hat x_0)|$ in what follows.  Assume $M\ge N_0(n,\ve_0,\ve_1,K_0+1,\beta)$.  Suppose for some $t<T_{K_0} (\le T_{K_0+1})$ there are $\ve\in[0,2^{-M}], |x|\le K_0+1, \hat x_0,x'\in\Re$ satisfying $|x-x'|\le 2^{-M}$, $|\hat x_0-x|\le \ve$, $|u(t,\hat x_0)|\le a_n\wedge(\sqrt{a_n}\ve)$, and $|u'_{1,a_n}(t,\hat x_0)|\le a_n^\beta$. If $x\neq x'$, then $0<|x-x'|\vee\ve\le 2^{-M}\le 2^{-N_0}$ and we may choose $N\ge N_0$ so that $2^{-N-1}< \ve\vee|x-x'|\le 2^{-N}$. Then \eqref{u'1cond} holds and so by \eqref{U1bnd1},
\begin{align*}
|u'_{1,a_n^\alpha}(t,x)-u'_{1,a_n^\alpha}(t,x')|\le 2^{-82}a_n^{-\ve_0-{3\ve_1\over 2}}|x-x'|^{1-\ve_0}\Bigl[&a_n^{-{3\beta\over 2}}(\ve\vee |x-x'|)^{2\gamma}+1\\
&+a_n^{\beta(\gamma-{3\over 2})}(\ve\vee |x-x'|)^{\gamma}\Bigr].
\end{align*}
If $x=x'$ the above is trivial.  This implies $U^{(1)}_{M,n,\beta}=T_{K_0}$ by its definition.  We have therefore shown
\[P(U^{(1)}_{M,n,\beta}<T_{K_0})\le P(M<N_0).\]
This completes the proof because $N_0(n,\ve_0,\ve_1,K_0+1,\beta)$ is stochastically bounded uniformly in $(n,\beta)$ (as in \eqref{betarange}).\gdm

Turning next to $u_{2,a_n^\alpha}$, for $0<\beta\le {1\over 2}-\ve_1$, define
\begin{align*}
U^{(2)}_{M,n,\beta}&=\inf\Bigl\{t:\hbox{there are }\ve\in[0,2^{-M}],\ |x|\le K_0+1,\ \hat x_0, x'\in\Re,\  s.t.\  |x-x'|\le 2^{-M},\\
&\phantom{=\inf\{t:}\ |x-\hat x_0|\le \ve, \ |u(t,\hat x_0)|\le a_n\wedge(\sqrt{a_n}\ve),\ |u'_{1,{a_n}}(t,\hat x_0)|\le a_n^\beta, \hbox{ and }\\
&\phantom{=\inf\{t:}|u_{2,a_n^\alpha}(t,x)-u_{2,a_n^\alpha}(t,x')|>2^{-87}a_n^{-\ve_0}\Bigl(|x-x'|^{{1-\ve_0\over 2}}\Bigl[(\sqrt{a_n}\vee\ve\vee|x'-x|)^{2\gamma}\\
&\phantom{|u'_{1,a_n^\alpha}(t,x)-u'_{1,a_n^\alpha}(t,x')>2^{-87}a_n^{-\ve_0}}+a_n^{\beta\gamma}(\sqrt{a_n}\vee\ve\vee|x'-x|)^{\gamma}\Bigr]\\
&\phantom{=\inf\Bigl\{t:\hbox{there are }\ve\in[0,2^{-M}],\ |x}\ +|x-x'|^{1-\ve_0}a_n^{\beta+{\ve_1\over 4}}\Bigr)\Bigr\}\wedge T_{K_0}.
\end{align*}
Define $U^{(2)}_{M,n,0}$ by the same expression (with $\beta=0$) but with the condition on $|u'_{1,a_n^\alpha}(t,\hat x_0)|$ omitted.  Just as for $U^{(1)}$, $U^{(2)}_{M,n,\beta}$ is an $\cF_t$-stopping time. 

\begin{lemma}\label{U2conv}For each $n\in \N$ and $\beta$ as in \eqref{betarange}, $U^{(2)}_{M,n,\beta}\uparrow T_{K_0}$ as $M\uparrow\infty$, and in fact
\[\lim_{M\to\infty} \sup_{n,0\le \beta\le {1\over 2}-\ve_1}P(U^{(2)}_{M,n,\beta}<T_{K_0})=0.\]
\end{lemma}
\paragraph{Proof.}  As before we only need to show the second assertion.  Proposition~\ref{Pminduction} allows us to apply Proposition~\ref{u2modulus} with $m=\om+1$, $\eta_1=\ve_0$, $K=K_0+1$, and $\alpha,\beta$ as in \eqref{betarange}, \eqref{alphadef}.  Hence there is an $N_0=N_0(n,\ve_0,\ve_1,K_0+1,\beta)\in\N$ a.s., stochastically bounded uniformly in $(n,\beta)$ (as in \eqref{betarange}), and such that if
\begin{equation}\label{U2cond}
N\ge N_0(\omega),\ (t,x)\in Z(N,n,K_0+1,\beta),\ |x-x'|\le 2^{-N},
\end{equation}
then
\begin{align}
\nonumber&|u_{2,a_n^\alpha}(t,x)-u_{2,a_n^\alpha}(t,x')|\\
\nn&\le 2^{-89}a_n^{-\ve_0}\Bigl\{|x-x'|^{{1-\ve_0\over 2}}\Bigl[(\sqrt{a_n}\vee 2^{-N})^{2\gamma}+a_n^{\beta\gamma}(\sqrt{a_n}\vee 2^{-N})^{\gamma}\Bigr]\\
\label{U2inc1}&\phantom{\le 2^{-89}|a_n^{-\ve_0}\Bigl\{}+|x-x'|^{1-\ve_0}\Bigl[a_n^{(\beta+\ve_1)(2\gamma-{1\over 2})}+a_n^{\beta\gamma}a_n^{(\beta+\ve_1)(\gamma-{1\over 2})}\Bigr]\Bigr\}.
\end{align}
Since $\gamma>{3\over 4}$, 
\begin{align*}
a_n&^{(\beta+\ve_1)(2\gamma-{1\over 2})}+a_n^{\beta\gamma}a_n^{(\beta+\ve_1)(\gamma-{1\over 2})}\le a_n^{\beta+\ve_1}+a_n^{\beta+{\ve_1\over 4}}\le 2a_n^{\beta+{\ve_1\over 4}}.
\end{align*}
Therefore \eqref{U2inc1} shows that \eqref{U2cond} implies
\begin{align}\nn&|u_{2,a_n^\alpha}(t,x)-u_{2,a_n^\alpha}(t,x')|\\
\nn&\le 2^{-89}a_n^{-\ve_0}\Bigl\{|x-x'|^{{1-\ve_0\over 2}}\Bigr[(\sqrt{a_n}\vee 2^{-N})^{2\gamma}+a_n^{\beta\gamma}(\sqrt{a_n}\vee 2^{-N})^{\gamma}\Bigr]\\
&\label{U2inc2}\phantom{\le 2^{-89}a_n^{-\ve_0}\Bigl\{}+2|x-x'|^{1-\ve_0}a_n^{\beta+{\ve_1\over 4}}\Bigr\}.
\end{align}
The proof is now completed just as for Lemma~\ref{U1conv} where \eqref{U2inc2} is used in place of \eqref{U1bnd1}.
\gdm

\noindent{\bf Notation.} 
\begin{equation}\label{tildel}
\tilde\Delta_{u'_1}(n,\ve,\ve_0,\beta)=a_n^{-\ve_0}\ve^{-\ve_0}\Bigl\{\ve+(\ve a_n^{-3/4}+a_n^{-1/4})\Bigl(\ve^{2\gamma}+ a_n^{\beta\gamma}(\ve\vee\sqrt{a_n})^{\gamma}\Bigr)\Bigr\}.
\end{equation}
\medskip

For $0<\beta\le {1\over 2}-\ve_1$, define
\begin{align*}
U^{(3)}_{M,n,\beta}&=\inf\Bigl\{t:\hbox{there are }\ve\in[2^{-a_n^{-(\beta+\ve_1)\ve_0/4}},2^{-M}],\ |x|\le K_0+1,\ \hat x_0\in\Re,\  s.t.\\\
&\phantom{=\inf\{t:}\ |x-\hat x_0|\le \ve, \ |u(t,\hat x_0)|\le a_n\wedge(\sqrt{a_n}\ve),\ |u'_{1,{a_n}}(t,\hat x_0)|\le a_n^\beta, \hbox{ and }\\
&\phantom{=\inf\{t:}|u'_{1,a_n}(t,x)-u'_{1,a_n^\alpha}(t,x)|>2^{-78}(\tilde\Delta_{u_1'}(n,\ve,\ve_0,\beta)+a_n^{\beta+{\ve_1\over 8}})\Bigr\}\wedge T_{K_0}.
\end{align*}
Define $U^{(3)}_{M,n,0}$ by the same expression (with $\beta=0$) but with the condition on $|u'_{1,a_n^\alpha}(t,\hat x_0)|$ omitted.  Just as for $U^{(1)}$, $U^{(3)}_{M,n,\beta}$ is an $\cF_t$-stopping time. 

\begin{lemma}\label{U3conv}For each $n\in \N$ and $\beta$ as in \eqref{betarange}, $U^{(3)}_{M,n,\beta}\uparrow T_{K_0}$ as $M\uparrow\infty$, and in fact
\[\lim_{M\to\infty} \sup_{n,0\le \beta\le {1\over 2}-\ve_1}P(U^{(3)}_{M,n,\beta}<T_{K_0})=0.\]
\end{lemma}
\paragraph{Proof.} It suffices to prove the second assertion.  By Proposition~\ref{Pminduction} we may apply Proposition~\ref{Fmodulus2} with $m=\om+1$, $\eta_1=\ve_0$, $K=K_0+1$ and $\beta$ as in \eqref{betarange}.  Note also that if $s=t-a_n^\alpha+a_n$, then 
\begin{equation}\label{for*}
\sqrt{t-s}\le a_n^{\alpha/2}=a_n^{\beta+\ve_1},
\end{equation}
and
\begin{equation}\label{for*2}
a_n^{\beta+\ve_1}\le N^{-4/\ve_0}\iff 2^{-N}\ge 2^{-a_n^{-(\beta+\ve_1)\ve_0/4}}.
\end{equation}
So Proposition~\ref{Fmodulus2} shows there is an $N_0=N_0(n,\ve_0,K_0+1,\beta)\in\N$ a.s., stochastically bounded uniformly in $(n,\beta)$, and so that if 
\begin{equation}\label{U3cond}
N\ge N_0,\ (t,x)\in Z(N,n,K_0+1,\beta)\hbox{ and }2^{-N}\ge 2^{-a_n^{-(\beta+\ve_1)\ve_0/4}},
\end{equation}
then
\begin{align}
\nn&|u'_{1,a_n}(t,x)-u'_{1,a_n^\alpha}(t,x)|\\
\nn&=|F_{a_n}(t-a_n^\alpha+a_n,t,x)-F_{a_n}(t,t,x)|\quad\hbox{ (by \eqref{Fu1})}\\
\label{U3inc1}&\le 2^{-81}a_n^{-\ve_0}\Bigl\{2^{N\ve_0}\Bigl[2^{-N}+a_n^{-1/4}(2^{-N}a_n^{-1/2}+1)(2^{-2\gamma N}+a_n^{\beta\gamma}(2^{-N}\vee\sqrt{a_n})^{\gamma})\Bigr]\\
\nn&\phantom{\le 2^{-81}(4+\om)^{-1}a_n^{-\ve_0}\Bigl\{2^{N\ve_0}\Bigl[}+a_n^{(\beta+\ve_1)(1-\ve_0)}\Bigl[a_n^{(\beta+\ve_1)(2\gamma-{3\over 2})}+a_n^{\beta\gamma}a_n^{(\beta+\ve_1)(\gamma-{3\over 2})}\Bigr]\Bigr\}.
\end{align}
We have used $\beta+\ve_1\le {1\over 2}$ (from \eqref{betarange}) in the last line.
The fact that 
$$(\beta+\ve_1)(2\gamma-{3\over 2})>\beta\gamma+(\beta+\ve_1)(\gamma-{3\over 2})$$ 
implies that
\begin{align*}
&a_n^{-\ve_0}a_n^{(\beta+\ve_1)(1-\ve_0)}\Bigl[a_n^{(\beta+\ve_1)(2\gamma-{3\over 2})}+ a_n^{\beta\gamma+(\beta+\ve_1)(\gamma-{3\over 2})}\Bigr]\\
&\le 2 a_n^{-\ve_0}a_n^{(\beta+\ve_1)(1-\ve_0)+\beta(2\gamma-{3\over 2})+\ve_1(\gamma-{3\over 2})}\\
&\le 2a_n^{\beta(2\gamma-{1\over 2}-\ve_0)+\ve_1(\gamma-{1\over 2}-\ve_0)-{\ve_1\over 100}}\\
&\le 2a_n^{\beta+{\ve_1\over 8}},
\end{align*}
where \eqref{veconditions} is used in the last two inequalities. This allows us to simplify \eqref{U3inc1} and show that \eqref{U3cond} implies
\begin{equation}\label{U3inc2}
|u'_{1,a_n}(t,x)-u'_{1,a_n^\alpha}(t,x)|\le 2^{-81}[\tilde\Delta_{u'_1,\ve_0}(n,2^{-N},\ve_0,\beta)+2a_n^{\beta+{\ve_1\over 8}}].
\end{equation}

The proof now is similar to that of Lemma~\ref{U1conv}.  As before, we may assume $\beta>0$.  Assume $M\ge N_0(n,\ve_0,K_0+1,\beta)$.  Suppose for some $t<T_{K_0}$ there are $\ve\in[2^{-a_n^{-(\beta+\ve_1)\ve_0/4}},2^{-M}]$, $|x|\le K_0+1$, and $\hat x_0\in\Re$, such that $|\hat x_0-x|\le \ve$, $|u(t,\hat x_0)|\le a_n\wedge (\sqrt{a_n}\ve)$, and $|u'_{1,a_n}(t,\hat x_0)|\le a_n^\beta$.  We may choose $N\ge M\ge N_0(\omega)$ so that $2^{-N-1}<\ve \le 2^{-N}$. Then $(t,x)\in Z(N,n,K_0+1,\beta)$ and $2^{-N}\ge \ve\ge 2^{-a_n^{-(\beta+\ve_1)\ve_0/4}}$, and therefore \eqref{U3cond} holds.  Therefore we may use \eqref{U3inc2} and the fact that $\tilde\Delta(n,2\ve,\ve_0,\beta)\le 8\tilde\Delta(n,\ve,\ve_0,\beta)$ to see that
\begin{align*}
|u'_{1,a_n}(t,x)-u'_{1,a_n^\alpha}(t,x)|&\le 2^{-78}[\tilde\Delta_{u'_1}(n,\ve,\ve_0,\beta)+a_n^{\beta+{\ve_1\over 8}}].
\end{align*}
This shows that $M\ge N_0(n,\ve_0,K_0+1,\beta)$ implies $U^{(3)}_{M,n,\beta}=T_{K_0}$ and so
\[\sup_{n,0\le \beta\le {1\over 2}-\ve_1}P(U^{(3)}_{M,n,\beta}<T_{K_0})\le \sup_{n,0\le \beta\le {1\over 2}-\ve_1}P(N_0>M)\to 0\hbox{ as }M\to\infty\]
by the stochastic boundedness of $N_0$ uniformly in $(n,\beta)$.\gdm

Finally for $M\in\N$, define
\begin{align*}
U^{(4)}_{M}&=\inf\Bigl\{t:\hbox{there are }\ve\in[0,2^{-M}],\ |x|\le K_0+1,\ \hat x_0, x'\in\Re,\  s.t.\  |x-x'|\le 2^{-M},\\
&\phantom{=\inf\{t:}\ |x-\hat x_0|\le \ve, \ |u(t,\hat x_0)|\le \ve, \hbox{ and }|u(t,x)-u(t,x')|>(\ve\vee|x'-x|)^{1-\ve_0}\Bigr\}\wedge T_{K_0}.
\end{align*}

\begin{lemma}\label{U4conv} $U^{(4)}_{M}\uparrow T_{K_0}$ as $M\uparrow\infty$, and 
\[\lim_{M\to\infty} P(U^{(4)}_{M}<T_{K_0})=0.\]
\end{lemma}
\paragraph{Proof.} It suffices to prove the second result.  This follows easily from Theorem~\ref{theorem:lipmod} as in the proof of Lemma~\ref{U1conv}.  
The constant multiplicative factors arising in the proof can easily be handled by applying Theorem~\ref{theorem:lipmod} with $\xi=1-\ve_0/2$ in place of $\xi=1-\ve_0$.\gdm

Let
\[U_{M,n,\beta}=\wedge_{j=1}^3 U^{(j)}_{M,n,\beta},\]
and
\begin{equation}\label{UMdef}
U_{M,n}=\Bigl(\wedge_{i=0}^{L(\ve_0,\ve_1)}U_{M,n,\beta_i}\Bigr)\wedge U_M^{(4)},
\end{equation}
where we recall that $\{\beta_i:i\le L\}$ were introduced in \eqref{betai}.
We have suppressed the dependence of $U_{M,n}$ on our fixed values of $K_0,\ve_0$ and $\ve_1$. 
Note that $\beta_i\in [0,{1\over 2}-\ve_1]$ for $i=0,\dots,L$ by \eqref{betarange} and $\alpha_i=\alpha(\beta_i)$. 
Lemmas~\ref{U1conv}, \ref{U2conv}, \ref{U3conv} and \ref{U4conv} therefore show that $\{U_{M,n}\}$ satisfy hypothesis 
$(H_1)$ of Proposition~\ref{prop:Inbound}.  Hence to complete the proof of Proposition~\ref{tildeJ} it suffices to establish compactness of $\tilde J_{n,i}(s)$, and the inclusion $\tilde J_{n,i}(s)\supset J_{n,i}(s)$ for all $s< U_{M,n}$, $(n,M)$ as in \eqref{ncond}, and $i=0,\ldots, L$. The next lemmas will show the inclusion part of the proof.  We assume $(n,M)$ satisfies \eqref{ncond} throughout the rest of this Section.

\medskip

\begin{lemma}\label{u1'lem}
If $i\in\{0,\dots,L\}$, $0\le s<U_{M,n}$, and $x\in J_{n,i}(s)$, then
\noindent(a) $|u'_{1,a_n}(s,\hat x_n(s,x))-u'_{1,a_n^{\alpha_i}}(s,\hat x_n(s,x))|\le 2^{-74}a_n^{\beta_i+{\ve_1\over 8}}$,

\noindent(b) for $i>0$, $|u'_{1,a_n^{\alpha_i}}(s,\hat x_n(s,x))|\le a_n^{\beta_i}/2$,

\noindent(c) for $i<L$, $u'_{1,a_n^{\alpha_i}}(s,\hat x_n(s,x))\ge a_n^{\beta_{i+1}}/8$.
\end{lemma}
\paragraph{Proof.} (a) Assume $(n,i,s,x)$ are as above and set $\ve=\sqrt{a_n}$.  We have $|\langle u_s,\Phi^{m_{n+1}}_x\rangle|\le a_n$ and \[Supp(\Phi_x^{m_{n+1}})\subset [x-\sqrt{a_n},x+\sqrt{a_n}].\]  
Using the continuity of $u(s,\cdot)$, we conclude that
\begin{equation}\label{hatxcond}
|u(s,\hat x_n(s,x))|\le a_n=a_n\wedge(\sqrt{a_n}\ve),\quad |\hat x_n(s,x)-x|\le \ve.
\end{equation}
The definition of $J_{n,i}$ also implies 
\begin{equation}\label{u1'bndipos}
|u'_{1,a_n}(s,\hat x_n(s,x))|\le a_n^{\beta_i}/4\quad\hbox{for }i>0.
\end{equation}
In addition, \eqref{ncond}  and $\ve_1<1/2$ (by\eqref{veconditions}) imply 
\begin{equation}\label{cond*ch}
2^{-M}\ge \sqrt{a_n}=\ve\ge 2^{-a_n^{-\ve_0\ve_1/4}}.
\end{equation}
Combine \eqref{hatxcond}, \eqref{u1'bndipos} and \eqref{cond*ch} with $|\hat x_n(s,x)|\le K_0+1$ and $s<U_{M,n}\le U^{(3)}_{M,n,\beta_i}$, and take $x=\hat x_0=\hat x_n(s,x)$ in the definition of $U^{(3)}$, to conclude that 
\begin{equation}\label{u1'incipos}
|u'_{1,a_n}(s,\hat x_n(s,x))-u'_{1,a_n^{\alpha_i}}(s,\hat x_n(s,x))|\le 2^{-78}(\tilde\Delta_{u'_1}(n,\sqrt{a_n},\ve_0,\beta_i)+a_n^{\beta_i+{\ve_1\over 8}}).
\end{equation}
Now 
\begin{align*}
\tilde\Delta_{u'_1}(n,\sqrt{a_n},\ve_0,\beta_i)
&\le a_n^{-3\ve_0/2}\Bigl[\sqrt{a_n}+2a_n^{-1/4}(a_n^\gamma+a_n^{\gamma(\beta_i+{1\over 2})})\Bigr]\\
&\le 4a_n^{-3\ve_0/2}\Bigl[ \sqrt{a_n} +a_n^{\gamma(\beta_i+{1\over 2})-{1\over 4}}\Bigr]\quad\hbox{(since $\beta_i<{1\over 2}$)}\\
&\le 4\Bigl[ a_n^{{1-3\ve_0\over 2}}+a_n^{\beta_i-{3\ve_0\over 2}+{3\ve_1\over 2}}\Bigr]\quad\hbox{(using $\gamma>{3\over 4}$ and $\beta_i\le {1\over 2}-6\ve_1$)}\\
&\le 8 a_n^{\beta_i+{\ve_1\over 8}}.
\end{align*}
The last line follows from \eqref{veconditions} and a bit of arithmetic. Use the above bound in \eqref{u1'incipos} and conclude that
\[|u'_{1,a_n}(s,\hat x_n(s,x))-u'_{1,a_n^{\alpha_i}}(s,\hat x_n(s,x))|\le 2^{-74}a_n^{\beta_i+{\ve_1\over 8}}.\]

\noindent(b) This is immediate from (a) and the fact that $|u'_{1,a_n}(s,\hat x_n(s,x))|\le a_n^{\beta_i}/4$ (by the definition of $J_{n,i}$ for $i>0$).

\noindent(c)  Since $\ve_0\le \ve_1/8$ by \eqref{veconditions}, $a_n^{\beta_i+{\ve_1\over 8}}\le a_n^{\beta_{i+1}}$.  For $i<L$ we have $u'_{1,a_n}(s,\hat x_n(s,x))\ge a_n^{\beta_{i+1}}/4$. The result is now clear from (a) and the triangle inequality.
\gdm

\begin{lemma}\label{u'1lem2}
If $i\in\{0,\dots,L\}$, $0\le s<U_{M,n}$, $x\in J_{n,i}(s)$ and $|x-x'|\le 5\oln(\beta_i)$, then

\noindent (a) for $i>0$, $|u'_{1, a_n^{\alpha_i}}(s,x')|\le a_n^{\beta_i},$

\noindent (b) for $i<L$, $u'_{1, a_n^{\alpha_i}}(s,x')\ge a_n^{\beta_{i+1}}/16$.
\end{lemma}
\paragraph{Proof.} Let $(n,i,s,x,x')$ be as above and set $\ve=|x-x'|+\sqrt{a_n}$.  Then we have \eqref{hatxcond}, \eqref{u1'bndipos}, and also (by \eqref{ncond})
\begin{equation}\label{veBnd} \ve\le 5a_n^{5\ve_1}+\sqrt{a_n}\le2^{-M},
\end{equation}
and
\begin{equation}\label{hatxcond2}
|x'-\hat x_n(s,x)|\le \ve\le 2^{-M},\quad |\hat x_n(s,x)|\le |x|+1\le K_0+1.
\end{equation}
\eqref{hatxcond}, \eqref{u1'bndipos}, \eqref{veBnd} and \eqref{hatxcond2} allow us to use $s<U_{M,n}\le U^{(1)}_{M,n,\beta_i}$ and the definition of $U^{(1)}_{M,n,\beta_i}$ (for $i>0$ or $i=0$), with $\hat x_n(s,x)$ playing the role of $x$, and deduce that
\begin{align*}
&|u'_{1,a_n^{\alpha_i}}(s,x')-u'_{1,a_n^{\alpha_i}}(s,\hat x_n(s,x))|\\
&\le 2^{-82}a_n^{-\ve_0-{3\ve_1\over 2}}(|x-x'|+\sqrt{a_n})^{1-\ve_0}\Bigl[a_n^{-3\beta_i/2}(|x'-x|+\sqrt{a_n})^{2\gamma}\\
&\phantom{\le 2^{-82}a_n^{-\ve_0-{3\ve_1\over 2}}(|x-x'|+\sqrt{a_n})^{1-\ve_0}\Bigl[}+1+a_n^{\beta_i(\gamma-{3\over 2})}(|x-x'|+\sqrt{a_n})^{\gamma}\Bigr].
\end{align*}
Use the fact that $\beta_i+5\ve_1\le 1/2$ (recall \eqref{betai}) to infer $|x-x'|+\sqrt{a_n}\le 6 a_n^{\beta_i+5\ve_1}\le a_n^{\beta_i}$ (by \eqref{ncond}) and so bound the above by
\begin{align*}
2^{-79}&a_n^{-\ve_0-{3\ve_1\over 2}}a_n^{(\beta_i+5\ve_1)(1-\ve_0)}\Bigl[a_n^{\beta_i(2\gamma-{3\over 2})}+1+a_n^{\beta_i(2\gamma-{3\over 2})}\Bigr]\\
&\le 2^{-79}a_n^{-\ve_0-{3\ve_1\over 2}}a_n^{(\beta_i+5\ve_1)(1-\ve_0)}3\\
&\le 2^{-77}a_n^{\beta_{i+1}},
\end{align*}
provided that $\beta_{i+1}\le (\beta_i+5\ve_1)(1-\ve_0)-\ve_0-{3\ve_1\over 2}$, or equivalently
\[\ve_0+(\beta_i+1)\ve_0\le ((7/2)-5\ve_0)\ve_1.\]
This follows easily from \eqref{veconditions}.  We have therefore shown that 
\begin{equation*} 
|u'_{1,a_n^{\alpha_i}}(s,x')-u'_{1,a_n^{\alpha_i}}(s,\hat x_n(s,x))|\le 2^{-77}a_n^{\beta_{i+1}},
\end{equation*}
 and so both (a) and (b) are now immediate from Lemma~\ref{u1'lem} (b), (c).
\gdm

\begin{lemma}\label{u2lem}
If $i\in\{0,\dots,L\}$, $0\le s<U_{M,n}$, $x\in J_{n,i}(s)$, and $|x-x'|\le 4\sqrt{a_n}$, then
\begin{align*}
|u_{2,a_n^{\alpha_i}}(s,x')-u_{2,a_n^{\alpha_i}}(s,x'')|\le 2^{-75}a_n^{\beta_{i+1}}(|x'-x''|\vee a_n^{\gamma-2\beta_i(1-\gamma)-\ve_1})\hbox{ whenever }|x'-x''|\le \oln(\beta_i).
\end{align*}
\end{lemma}
\paragraph{Proof.} Assume $(i,n,s,x,x')$ are as above and set $\ve=5\sqrt{a_n}\le 2^{-M}$, by \eqref{ncond}.  Then 
\[|x'-\hat x_n(s,x)|\le |x'-x|+\sqrt{a_n}\le \ve,\quad|x'|\le |x|+1\le K_0+1,\quad |u(s,\hat x_n(s,x))|\le a_n=a_n\wedge (\sqrt{a_n}\ve),\]
and the definition of $(s,x)\in J_{n,i}$ implies that for $i>0$,
\begin{equation*}
|u'_{1,a_n}(s,\hat x_n(s,x))|\le a_n^{\beta_i}/4\le a_n^{\beta_i}.
\end{equation*}
Let 
\[Q(n,\ve_0,\beta_i,r)=a_n^{-\ve_0}r^{{1-\ve_0\over 2}}\Bigl[(\sqrt{a_n}\vee r)^{2\gamma}+a_n^{\beta_i\gamma}(\sqrt{a_n}\vee r)^{\gamma}\Bigr].\]
Assume $|x'-x''|\le \oln(\beta_i)\le 2^{-M}$, the last by \eqref{ncond}.  The condition $s<U_{M,n}\le U^{(2)}_{M,n,\beta_i}$ and the definition of $U^{(2)}$, with $(x',x'')$ playing the role of $(x,x')$, ensures that 
\begin{align}
\nn |u_{2,a_n^{\alpha_i}}(s,x'')-u_{2,a_n^{\alpha_i}}(s,x')|&\le 2^{-87}a_n^{-\ve_0}\Bigl[|x''-x'|^{{1-\ve_0\over 2}}\Bigl[((5\sqrt{a_n})\vee|x''-x'|)^{2\gamma}\\
\nn&\phantom{\le 2^{-87}a_n^{-\ve_0}\Bigl[}+a_n^{\beta_i\gamma}((5\sqrt{a_n})\vee|x''-x'|)^{\gamma}\Bigr]+|x''-x'|^{1-\ve_0}a_n^{\beta_i+{\ve_1\over 4}}\Bigr]\\
\label{u2lem2}&\le 2^{-82}[Q(n,\ve_0,\beta_i,|x''-x'|)+|x''-x'|^{1-\ve_0}a_n^{\beta_i+{\ve_1\over 4}-\ve_0}].
\end{align}

We first show that 
\begin{equation}\label{QBnd}
Q(n,\ve_0,\beta_i,r)\le 2a_n^{\beta_{i+1}}(r\vee a_n^{\gamma-2\beta_i(1-\gamma)-\ve_1})\hbox{ for }0\le r\le \oln(\beta_i).
\end{equation}
{\bf Case 1.} $\sqrt{a_n}\le r\le \oln(\beta_i)$.
\[Q(n,\ve_0,\beta_i,r)= a_n^{-\ve_0}\Bigl[r^{2\gamma+{1\over 2}-{\ve_0\over 2}}+a_n^{\beta_i\gamma}r^{\gamma+{1\over 2}-{\ve_0\over 2}}\Bigr],\]
and so \eqref{QBnd}, will hold if 
\begin{equation}\label{u2lem3}r^{2\gamma-{1\over 2}-{\ve_0\over 2}}\le a_n^{\beta_{i+1}+\ve_0},
\end{equation}
and
\begin{equation}\label{u2lem4}
a_n^{\beta_i\gamma}r^{\gamma-{1\over 2}-{\ve_0\over 2}}\le a_n^{\beta_{i+1}+\ve_0}.
\end{equation}
\eqref{veconditions} implies $2\gamma-{1\over 2}-{\ve_0\over 2}>1$ and so \eqref{u2lem3} would follow from
\begin{equation*} r\le a_n^{\beta_{i+1}+\ve_0}.
\end{equation*}
Hence, by the upper bound on $r$ in this case, it suffices to show that $a_n^{{\beta_i}+5\ve_1}\le a_n^{{\beta_i}+2\ve_0}$ and this is immediate from \eqref{veconditions}.  

Turning to \eqref{u2lem4}, note that 
\begin{align*}
a_n^{\beta_i\gamma}r^{\gamma-{1\over 2}-{\ve_0\over 2}}a_n^{-\beta_{i+1}-\ve_0}
&\le a_n^{\beta_i\gamma+(\beta_i+5\ve_1)(\gamma-{1\over 2}-{\ve_0\over 2})-\beta_{i+1}-\ve_0}\\
&\le a_n^{\beta_i(2\gamma-{3\over 2}-{\ve_0\over 2})+5\ve_1(\gamma-{1\over 2}-{\ve_0\over 2})-2\ve_0}\\
&\le a_n^{5\ve_1(\gamma-{1\over 2}-{\ve_0\over 2})-2\ve_0}\le 1,
\end{align*}
where in the last two inequalities we are using \eqref{veconditions}.   This proves \eqref{u2lem4} and hence completes the derivation of \eqref{QBnd} in this case.

\noindent {\bf Case 2.} $a_n^{\gamma-2\beta_i(1-\gamma)-\ve_1}\le r<\sqrt{a_n}$.\hfil\break
Then
\[Q(n,\ve_0,\beta_i,r)=a_n^{-\ve_0}r^{{1-\ve_0\over 2}}\Bigl[a_n^\gamma+a_n^{\gamma(\beta_i+{1\over 2})}\Bigr]\le 2a_n^{-\ve_0}r^{{1-\ve_0\over 2}}a_n^{\gamma(\beta_i+{1\over 2})},\]
and so \eqref{QBnd} will hold if 
\begin{equation}\label{u2lem5}
r^{{1+\ve_0\over 2}}\ge a_n^{-\beta_{i+1}-\ve_0+\gamma(\beta_i+{1\over 2})}=a_n^{-\beta_i(1-\gamma)+{\gamma\over 2}-2\ve_0}.
\end{equation}
Our lower bound on $r$ implies that
\[r^{{1+\ve_0\over 2}}\ge a_n^{(\gamma-2(\beta_i(1-\gamma)-\ve_1)(1+\ve_0)/2}\ge a_n^{-\beta_i(1-\gamma)+{\gamma\over 2}+{\gamma\ve_0\over 2}-{\ve_1\over 2}}\]
which implies \eqref{u2lem5} by \eqref{veconditions}.

\noindent {\bf Case 3.} $r<a_n^{\gamma-2\beta_i(1-\gamma)-\ve_1}$.\hfil\break
This case follows from Case 2 and the monotonicity of $Q(n,\ve_0,\beta_i,r)$ in $r$.  Strictly speaking we also need the fact that Case 2 is non-empty as was done in Lemma~\ref{intervalbounds}.

This completes the proof of \eqref{QBnd}.  Consider next the second term in \eqref{u2lem2}.  If $r\ge a_n$, then
\begin{align*}
r^{1-\ve_0}a_n^{\beta_i+{\ve_1\over 4}-\ve_0}(a_n^{\beta_{i+1}}r)^{-1}&=r^{-\ve_0}a_n^{-2\ve_0+{\ve_1\over 4}}\\
&\le a_n^{-3\ve_0+{\ve_1\over 4}}<1
\end{align*}
by \eqref{veconditions}.  It follows that 
\begin{equation}\label{u2lem7}
r^{1-\ve_0}a_n^{\beta_i+{\ve_1\over 4}-\ve_0}\le a_n^{\beta_{i+1}}(r\vee a_n)\le a_n^{\beta_{i+1}}(r\vee a_n^{\gamma-2\beta_i(1-\gamma)-\ve_1}),
\end{equation}
the last inequality being trivial. 

Insert \eqref{QBnd} and \eqref{u2lem7} into \eqref{u2lem2} to derive the desired bound.
\gdm

\begin{lemma}\label{i0lem} If $0\le s<U_{M,n}$ and $x\in J_{n,0}(s)$, then
\begin{equation}\label{i0lem1} |u(s,x)-u(s,x')|\le (\sqrt{a_n}\vee |x'-x|)^{1-\ve_0}\hbox{ whenever }|x-x'|\le 2^{-M},
\end{equation}
and
\begin{equation}\label{i0lem2}
|u(s,x')|\le 3(\sqrt{a_n})^{1-\ve_0}\hbox{ whenever }|x'-x|\le \sqrt{a_n}.
\end{equation}
\end{lemma}
\paragraph{Proof.}  As in \eqref{hatxcond},
 if $\ve=\sqrt{a_n}$, $(s,x)\in J_{n,0}$ implies
\begin{equation}\label{i0lem3}|u(s,\hat x_n(s,x))|\le a_n\le \ve,\ \ |\hat x_n(s,x)-x|\le \ve,\hbox{ and }|x|\le K_0.\end{equation}
In addition, \eqref{ncond} ensures that $\ve\le 2^{-M}$, and so $s<U_{M,n}\le U^{(4)}_{M}$ means that 
\[|u(s,x')-u(s,x)|\le (\sqrt{a_n}\vee |x'-x|)^{1-\ve_0}\quad\hbox{for all }|x'-x|\le 2^{-M}.\]
This proves \eqref{i0lem1}.  Next take $x'=\hat x_n(s,x)$ in the above inequality and use \eqref{i0lem3} to obtain for $|x'-x|\le \sqrt{a_n}$,
\begin{align*}
|u(s,x')|&\le |u(s,x')-u(s,x)|+|u(s,x)-u(s,\hat x_n(s,x))|+|u(s,\hat x_n(s,x))|\\
&\le 2\sqrt{a_n}^{1-\ve_0}+a_n\le 3\sqrt{a_n}^{1-\ve_0}.
\end{align*}
This proves \eqref{i0lem2}.\gdm

\paragraph{Proof of Proposition~\ref{tildeJ}} The compactness is elementary and left for the reader--note here that continuity allows us to replace the closed intervals on which the inequalities defining $\tilde J_{n,i}(s)$ hold with open intervals.

The inclusions $J_{n,i}(s)\subset \tilde J_{n,i}(s)$ for $0\le s<U_{M,n}$ are immediate from Lemmas~\ref{u'1lem2}, \ref{u2lem} and \ref{i0lem}.\gdm

This finishes the proof of Proposition~\ref{tildeJ} except for the proof of  Proposition \ref{u2modulus} which is the objective of the next section.
\section{Proof of Proposition \ref{u2modulus}}\label{sec6}
\setcounter{equation}{0}
\setcounter{theorem}{0} 

We now continue to assume $b\equiv 0$ and give the proof of Proposition~\ref{u2modulus}. Assume first $t'\ge t$ and use \eqref{u2eq} to write
\begin{align}\label{u2decompa}
&|u_{2,\delta}(t',x')-u_{2,\delta}(t,x)|\\
\nonumber&\le \Bigl|\int_{(t-\delta)^+}^{t\wedge(t'-\delta)^+}\int p_{t-s}(y-x)D(s,y)W(ds,dy)\Bigr|+\Bigl|\int_{{(t'-\delta)^+}\vee t}^{t'}\int p_{t'-s}(y-x')D(s,y)W(ds,dy)\Bigr|\\
\nonumber &\quad +1(t'-t<\delta)\Bigl|\int_{(t'-\delta)^+}^t\int(p_{t'-s}(y-x')-p_{t-s}(y-x))D(s,y)W(ds,dy)\Bigr|.
\end{align}
This decomposition and \eqref{Dbound} suggests we define the following square functions for $\delta\in(0,1]$ and $\eta_0\in(0,1/2)$ (noting also that $(t'-\delta)^+\vee t\ge t'-(\delta\wedge (t'-t))$):

\begin{align*}
&\hat Q_{T,1,\delta}(t,t',x)=\int_{(t-\delta)^+}^{t\wedge(t'-\delta)^+}\int p_{t-s}(y-x)^2e^{2R_1|y|}|u(s,y)|^{2\gamma}dyds,\\
&\hat Q_{T,2,\delta}(t,t',x')=\int_{t'-(\delta\wedge(t'-t))}^{t'}\int p_{t'-s}(y-x')^2e^{2R_1|y|}|u(s,y)|^{2\gamma}dyds,\\
&\hat Q_{S,1,\delta,\eta_0}(t,x,t',x')=1(t'-t<\delta)\int _{(t'-\delta)^+}^t\int 1(|y-x|>(t'-s)^{1/2-\eta_0}\vee (2|x-x'|))\\
&\phantom{Q_{S,1,\eta_0}(s,t,x,t',x')=\int _0^{(s-\delta)^+}\int}\times(p_{t'-s}(y-x')-p_{t-s}(y-x))^2e^{2 R_1|y|}|u(s,y)|^{2\gamma}dyds,\\
&\hat Q_{S,2,\delta,\eta_0}(t,x,t',x')=1(t'-t<\delta)\int _{(t'-\delta)^+}^t\int 1(|y-x|\le(t'-s)^{1/2-\eta_0}\vee (2|x-x'|))\\
&\phantom{Q_{S,1,\eta_0}(s,t,x,t',x')=\int _0^{(s-\delta)^+}\int}\times(p_{t'-s}(y-x')-p_{t-s}(y-x))^2e^{2 R_1|y|}|u(s,y)|^{2\gamma}dyds.
\end{align*}

\begin{lemma}\label{hQS1bound}
For any $K\in\N^{\ge K_1}$ and $R>2$ there is a $c_{\ref{hQS1bound}}(K,R)>0$ and an $N_{\ref{hQS1bound}}=N_{\ref{hQS1bound}}(K,\omega)\in\N$ a.s. such that for all $\eta_0,\eta_1\in(1/R,1/2)$, $\delta\in(0,1]$, $N,n\in\N$, $\beta\in[0,1/2]$ and $(t,x)\in\Rp\times\Re$, on
\begin{equation}\label{omegacond9}
\{\omega:(t,x)\in Z(N,n,K,\beta),N\ge N_{\ref{hQS1bound}}\},
\end{equation}
\begin{equation*}
\hat Q_{S,1,\delta,\eta_0}(t,x,t',x')\le c_{\ref{hQS1bound}}(K,R)2^{4N_{\ref{hQS1bound}}}[d((t,x),(t',x'))\wedge \sqrt \delta]^{2-\eta_1}\delta^{3/2}\quad\hbox{for all }t\le t' \hbox{ and }
x'\in\Re.
\end{equation*}
\end{lemma}
\paragraph{Proof.} The proof is quite similar to that of Lemma~\ref{QS1bound}.  We let $d=d((t,x,),(t',x'))$ and $N_{\ref{hQS1bound}}=N_1(0,3/4,K)$, where $N_1$ is as in $(P_0)$.  Recall here from Remark~\ref{m0case} that for $m=0$, $N_1$ depends only on $(\xi,K)$ and we take $\xi=3/4$.  We may assume $t'-t<\delta$.  Then for $\omega$ as in \eqref{omegacond9} and $t\le t'$, Lemma~\ref{uglobalbound}, with $m=0$, implies
\begin{align*}
&\hat Q_{S,1,\delta,\eta_0}(t,x,t',x')\\
&\le C_{\ref{uglobalbound}}(\omega)\int _{(t'-\delta)^+}^t\int 1(|y-x|>(t'-s)^{1/2-\eta_0}\vee (2|x-x'|))(p_{t'-s}(y-x')-p_{t-s}(y-x))^2\\
&\phantom{\le C_{\ref{uglobalbound}}(\omega)\int _{(t'-\delta)^+}^t\int 1(}\times e^{2(|y-x|+R_1|y|)}(2^{-N}\vee(\sqrt{t-s}+|y-x|))^{3\gamma/2}\,dyds\\
&\le  C_{\ref{uglobalbound}}(\omega)\int _{(t'-\delta)^+}^t\int 1(|y-x|>(t'-s)^{1/2-\eta_0}\vee (2|x-x'|))(p_{t'-s}(y-x')-p_{t-s}(y-x))^2\\
&\phantom{\le C_{\ref{uglobalbound}}(\omega)\int _{(t'-\delta)^+}^t\int 1(}e^{2R_1K}e^{2(R_1+1)|y-x|}[1+|y-x|]^{3\gamma/2}\,dyds\\
&\le c_1(K,R)C_{\ref{uglobalbound}}(\omega)\int_{t'-\delta}^t(t-s)^{-1/2}\exp\Bigl(-{\eta_1(t'-s)^{-2\eta_0}\over 33}\Bigr)\Bigl[1\wedge {d^2\over t-s}\Bigr]^{1-{\eta_1\over 2}}ds,
\end{align*}
where we have used Lemma~\ref{ptbnds}(b) in the last line. 
Now use 
\[\exp\Bigl(-{\eta_1(t'-s)^{-2\eta_0}\over 33}\Bigr)\le \exp\Bigl(-{\eta_1(t'-t)^{-2\eta_0}\over 66}\Bigr)+\exp\Bigl(-{\eta_1(t-s)^{-2\eta_0}\over 66}\Bigr)\]
to bound the above by
\begin{align*}
&c_1(K,R)C_{\ref{uglobalbound}}(\omega)\Bigl[\exp\Bigl(-{\eta_1(t'-t)^{-2\eta_0}\over 66}\Bigr)\int_{t'-\delta}^t(t-s)^{-1/2}\Bigl[1\wedge {d^2\over t-s}\Bigr]^{1-{\eta_1\over 2}}ds\\
&\phantom{c_1(K,R)C_{\ref{uglobalbound}}(\omega)\Bigl[}+\int_{t'-\delta}^t(t-s)^{-1/2}\exp\Bigl(-{\eta_1(t-s)^{-2\eta_0}\over 66}\Bigr)\Bigl[1\wedge{d^2\over t-s}\Bigr]^{1-{\eta_1\over 2}}ds\Bigr]\\
&\le c_2(K,R)C_{\ref{uglobalbound}}(\omega)\Bigl[\exp\Bigl(-{\eta_1(t'-t)^{-2\eta_0}\over 66}\Bigr)(d^2\wedge\delta)^{1/2}\\
&\phantom{c_1(K,R)C_{\ref{uglobalbound}}(\omega)\Bigl[}+\int_{t'-\delta}^t(t-s)^{3/2}\Bigl[1\wedge{d^2\over t-s}\Bigr]^{1-{\eta_1\over 2}}ds\Bigr]\quad\hbox{(use \eqref{J2})}\\
&\le c_3(K,R)C_{\ref{uglobalbound}}(\omega)(d^2\wedge \delta)^{1-{\eta_1\over 2}}\delta^{3/2}\qquad(\hbox{use \eqref{J1} and }t'-t\le d^2\wedge\delta).
\end{align*}
As we may take $\ve_0=0$ in the formula for $C_{\ref{uglobalbound}}(\omega)$ (by Remark~\ref{m0case}), the result follows.
\gdm

\begin{lemma}\label{hQS2bound} Let $0\le m\le \om+1$ and assume $(P_m)$.  For any $K\in\N^{\ge K_1}, R>2, n\in\N, \ve_0\in(0,1)$, and $\beta\in [0,1/2]$ there is a $c_{\ref{hQS2bound}}(K)$ and $N_{\ref{hQS2bound}}=N_{\ref{hQS2bound}}(m,n,R,\ve_0,K,\beta)(\omega)\in\N$ a.s. such that for any $\eta_1\in(R^{-1},1/2), \eta_0\in(0,\eta_1/24), \delta\in[a_n,1]$, $N\in\N$, and $(t,x)\in\Rp\times\Re$, 
on
\begin{equation}\label {omegacond10}
\{\omega:(t,x)\in Z(N,n,K,\beta), N\ge N_{\ref{hQS2bound}}\},
\end{equation}
\begin{align*}
&\hat Q_{S,2,\delta,\eta_0}(t,x,t',x')\\
&\le c_{\ref{hQS2bound}}(K)[a_n^{-2\ve_0}+2^{4N_{\ref{hQS2bound}}}]\Bigl[(d\wedge\sqrt\delta)^{1-{\eta_1\over 2}}\bar d_N^{2\gamma}[{\bar d_{n,N}}^{2\gamma(\tilde\gamma_m-1)}+a_n^{2\beta\gamma}]\\
&\phantom{\le C_{\ref{QS2bound}}(K,R)[a_n^{-2\ve_0}+2^{4N_3}]\Bigl[}+(d\wedge\sqrt\delta)^{2-\eta_1}[\delta^{\gamma\tilde\gamma_m-{1\over 2}}+a_n^{2\beta\gamma}\delta^{\gamma-{1\over 2}}]\Bigr]\\
&\phantom{\le C_{\ref{QS2bound}}(K,R)[a_n^{-2\ve_0}+2^{4N_3}]\Bigl[}\hbox{ for all } t\le t'\le K, |x'|\le K+1.
\end{align*}
Here $d=d((t,x),(t',x'))$, $\bar d_N=d\vee 2^{-N}$ and $\bar d_{n,N}=\sqrt{a_n}\vee {\bar d_N}$.  Moreover $N_{\ref{hQS2bound}}$ is stochastically bounded uniformly in $(n,\beta)$.
\end{lemma}
\paragraph{Proof.} Set $\xi=1-(24R)^{-1}$ and $N_{\ref{hQS2bound}}(m,n,R,\ve_0,K,\beta)=N_1(m,n,\xi,\ve_0,K,\beta)$, which is clearly stochastically bounded uniformly in $(n,\beta)$ by $(P_m)$.  We may assume $t'-t\le \delta$.  For $\omega$ as in \eqref{omegacond10}, $t\le t'\le K$ and $|x'|\le K+1$, Lemma~\ref{uglobalbound} implies
\begin{align}
\nonumber&\hat Q_{S,2,\delta,\eta_0}(t,x,t',x')\\
\nonumber&\le C_{\ref{uglobalbound}}(\omega)c_1(K)\int_{(t'-\delta)^+}^t\Bigr[\int (p_{t'-s}(y-x')-p_{t-s}(y-x))^2dy\Bigl]e^{2R_1K}e^{2(R_1+1)4(K+1)}\\
\nonumber&\phantom{\le C_{\ref{uglobalbound}}(\omega)c_1(K)\int}\times((2^{-N}\vee|x'-x|)+(t'-s)^{{1\over 2}-\eta_0})^{2\gamma\xi}\Bigl[((\sqrt {a_n}\vee 2^{-N}\vee|x'-x|)\\
\nonumber&\phantom{\le C_{\ref{uglobalbound}}(\omega)c_1(K)\int\times((2^{-N}\vee|x'-x|)+(t'-s)^{{1\over 2}-\eta_0})}+(t'-s)^{{1\over 2}-\eta_0})^{2\gamma(\tilde\gamma_m-1)}+a_n^{2\beta\gamma}\Bigr]ds\\
\nonumber&\le C_{\ref{uglobalbound}}(\omega)c_2(K)\int_{(t'-\delta)^+}^t(t-s)^{-1/2}\Bigl[1\wedge{d^2\over t-s}\Bigr]((2^{-N}\vee d^{1-2\eta_0})+(t-s)^{{1\over 2}-\eta_0})^{2\gamma\xi}\\
\label{hqS2bound1}&\phantom{\le C_{\ref{uglobalbound}}(\omega)c_1(K)\int}\times\Bigl[((\sqrt{a_n}\vee 2^{-N}\vee d^{1-2\eta_0})+(t-s)^{{1\over 2}-\eta_0})^{2\gamma(\tilde\gamma_m-1)}+a_n^{2\beta\gamma}\Bigl]ds.
\end{align}
We have used Lemma~\ref{ptbnds}(a) in the last line.  
Below we will implicitly use the conditions on $\eta_0,\eta_1$, $R$ and $\gamma$ to see that 
\[(1-2\eta_0)\gamma\xi>\Bigl({23\over 24}\Bigr)\Bigl({3\over 4}\Bigr)\Bigl({47\over 48}\Bigr)>{1\over 2},\]
and also use the conditions on $t,x,t',x'$ which imply $d\le c(K)$ (the latter was also used in \eqref{hqS2bound1}).  By considering separately the cases 
\begin{align*}
&(t-s)^{{1\over 2}-\eta_0}< 2^{-N}\vee d^{1-2\eta_0}, \quad (t-s)^{{1\over 2}-\eta_0}\ge \sqrt{a_n}\vee2^{-N}\vee d^{1-2\eta_0}\\
&\hbox{and }2^{-N}\vee d^{1-2\eta_0}\le (t-s)^{{1\over 2}-\eta_0}<\sqrt{a_n}\vee2^{-N}\vee d^{1-2\eta_0}
\end{align*}(the latter case implies $\sqrt{a_n}>2^{-N}\vee d^{1-2\eta_0}$), and then using Lemma~\ref{Jbnd} we may bound \eqref{hqS2bound1} by 
\begin{align}
\nonumber&C_{\ref{uglobalbound}}(\omega)c_3(K)\Bigl\{\int_{(t'-\delta)^+}^t(t-s)^{-1/2}\Bigl[1\wedge{d^2\over t-s}\Bigr]ds(2^{-N}\vee d^{{1}-2\eta_0})^{2\gamma\xi}\\
\nonumber&\phantom{C_{\ref{uglobalbound}}(\omega)c_3(K)\{\int_{(t'-\delta)}^t}\times\Bigl[(\sqrt{a_n}\vee 2^{-N}\vee d^{1-2\eta_0})^{2\gamma(\tilde\gamma_m-1)}+a_n^{2\beta\gamma}\Bigr]\\
\nonumber&\phantom{C_{\ref{uglobalbound}}(\omega)c_3(K)\Bigl\{}+\int_{(t'-\delta)^+}^t\Bigl[(t-s)^{-{1\over 2}+(1-2\eta_0)\gamma(\tilde\gamma_m-1+\xi)}+a_n^{2\beta\gamma}(t-s)^{-{1\over 2}+(1-2\eta_0)\gamma\xi}\Bigr]\Bigl[1\wedge{d^2\over t-s}\Bigr]ds
\\
\nonumber&\phantom{C_{\ref{uglobalbound}}(\omega)c_3(K)\Bigl\{}+\int_{(t'-\delta)^+}^t(t-s)^{(1-2\eta_0)\gamma\xi-{1\over 2}}\Bigl[1\wedge{d^2\over t-s}\Bigr]ds[a_n^{\gamma(\tilde\gamma_m-1)}+a_n^{2\beta\gamma}]
\Bigr\}\\
\nonumber&\le C_{\ref{uglobalbound}}(\omega)c_4(K)\Bigl\{(d^2\wedge\delta)^{1/2}(2^{-N}\vee d^{{1}-2\eta_0})^{2\gamma\xi}\\
\nonumber&\phantom{C_{\ref{uglobalbound}}(\omega)c_3(K)\{\int AA}\times\Bigl[(\sqrt{a_n}\vee 2^{-N}\vee d^{1-2\eta_0})^{2\gamma(\tilde\gamma_m-1)}+a_n^{2\beta\gamma}\Bigr]\\
\nonumber&\phantom{C_{\ref{uglobalbound}}(\omega)c_3(K)\{AA}+(d^2\wedge\delta)\Bigr[\delta^{(1-2\eta_0)\gamma(\tilde\gamma_m-1+\xi)-{1\over 2}}+a_n^{2\beta\gamma}\delta^{(1-2\eta_0)\gamma\xi-{1\over 2}}\Bigr]
\\
\nonumber&\phantom{C_{\ref{uglobalbound}}(\omega)c_3(K)\{AA}+(d^2\wedge\delta)\delta^{(1-2\eta_0)\gamma\xi-{1\over 2}}[a_n^{\gamma(\tilde\gamma_m-1)}+a_n^{2\beta\gamma}] \Bigr\}
\\
\nonumber&\phantom{C_{\ref{uglobalbound}}(\omega)c_3(K)\Bigl\{(d^2\wedge\delta)^{1/2}\Bigl[1\wedge{d^2\over t-s}\Bigr]ds(2^{-N}\vee d^{{1\over 2}-\eta_0})^{2\gamma\xi}}\hbox{(by \eqref{J2} and \eqref{J1}, respectively).}
\end{align}
The last term is less than the middle term because $\delta\in[a_n,1]$.  Therefore $\hat Q_{S,2,\delta,\eta_0}(t,x,t',x')$ is at most
\begin{align}
\nonumber& C_{\ref{uglobalbound}}(\omega)c_5(K)\Bigl\{(d\wedge\sqrt\delta)^{1-{\eta_1\over 2}}(2^{-N}\vee d)^{(1-2\eta_0)2\gamma\xi+{\eta_1\over 4}}\\
\label{hqS2bound2}&\phantom{C_{\ref{uglobalbound}}(\omega)c_3(K)\{\int}\times\Bigl[(\sqrt{a_n}\vee 2^{-N}\vee d)^{(1-2\eta_0)2\gamma(\tilde\gamma_m-1)+{\eta_1\over 4}}+a_n^{2\beta\gamma}(\sqrt{a_n}\vee 2^{-N}\vee d)^{{\eta_1\over 4}}\Bigr]\\
\nonumber&\phantom{C_{\ref{uglobalbound}}(\omega)c_3(K)\{AA}+(d^2\wedge\delta)^{1-{\eta_1\over 2}}\Bigr[\delta^{(1-2\eta_0)\gamma(\tilde\gamma_m-1+\xi)-{1\over 2}+{\eta_1\over 2}}+a_n^{2\beta\gamma}\delta^{(1-2\eta_0)\gamma\xi-{1\over 2}+{\eta_1\over 2}}\Bigr]\Bigr\}.
\end{align}
Our conditions on $\eta_0$, $\eta_1$, and $R$ imply\begin{align*}
&(1-2\eta_0)2\gamma\xi+{\eta_1\over 4}\ge 2\gamma,\ \ (1-2\eta_0)2\gamma(\tilde\gamma_m-1)+{\eta_1\over 4}\ge 2\gamma(\tilde\gamma_m-1),\\
 &(1-2\eta_0)\gamma(\tilde\gamma_m-1+\xi)+{\eta_1\over 2}\ge \gamma\tilde\gamma_m,\hbox{ and }(1-2\eta_0)\gamma\xi+{\eta_1\over 2}\ge \gamma.
 \end{align*}
Finally, insert the above bounds into \eqref{hqS2bound2} to derive the required bound on $\hat Q_{S,2,\delta,\eta_0}(t,x,t',x')$.  
\gdm

\begin{lemma}\label{hQTbound} Let $0\le m\le \om+1$ and assume $(P_m)$.  For any $K\in\N^{\ge K_1}, R>2, n\in\N, \ve_0\in(0,1)$, and $\beta\in [0,1/2]$ there is a $c_{\ref{hQTbound}}(K)$ and $N_{\ref{hQTbound}}=N_{\ref{hQTbound}}(m,n,R,\ve_0,K,\beta)(\omega)\in\N$ a.s. such that for any $\eta_1\in(R^{-1},1/2), \delta\in[a_n,1]$, $N\in\N$, and $(t,x)\in\Rp\times\Re$, 

\begin{equation}\label {omegacond11}
\hbox{on }\{\omega:(t,x)\in Z(N,n,K,\beta), N\ge N_{\ref{hQTbound}}\}, \hbox{ and for all }t\le t'\le T_K,\ |x'|\le K+1,
\end{equation}
\begin{align*}
&\hat Q_{T,1,\delta}(t,t',x)\\
&\le c_{\ref{hQTbound}}(K)[a_n^{-2\ve_0}+2^{4N_{\ref{hQTbound}}}]\Bigl[(\sqrt{t'-t}\wedge\sqrt\delta)^{1-{\eta_1\over 2}}\bar d_N^{2\gamma}[{\bar d_{n,N}}^{2\gamma(\tilde\gamma_m-1)}+a_n^{2\beta\gamma}]\\
&\phantom{\le C_{\ref{QS2bound}}(K,R)[a_n^{-2\ve_0}+2^{4N_3}]\Bigl[}+(\sqrt{t'-t}\wedge\sqrt\delta)^{2-\eta_1}[\delta^{\gamma\tilde\gamma_m-{1\over 2}}+a_n^{2\beta\gamma}\delta^{\gamma-{1\over 2}}]\Bigr],\\
\end{align*}
and
\begin{align*}
&\hat Q_{T,2,\delta}(t,t',x')\\
&\le c_{\ref{hQTbound}}(K)[a_n^{-2\ve_0}+2^{4N_{\ref{hQTbound}}}](\sqrt{t'-t}\wedge\sqrt\delta)^{1-{\eta_1\over 2}}\bar d_N^{2\gamma}\Bigl[{\bar d_{n,N}}^{2\gamma(\tilde\gamma_m-1)}+a_n^{2\beta\gamma}\Bigr].\\
\end{align*}
Here $d$, $\bar d_N$ and $\bar d_{n,N}$ are as in 
Lemma~\ref{hQS2bound}.  Moreover, $N_{\ref{hQTbound}}$ is stochastically bounded uniformly in $(n,\beta)$.
\end{lemma}
\paragraph{Proof.} Set $\xi=1-(4R)^{-1}$ and $N_{\ref{hQTbound}}(m,n,R,\ve_0,K,\beta)=N_1(m,n,\xi,\ve_0,K,\beta)$, which is clearly stochastically bounded uniformly in $(n,\beta)$ by $(P_m)$.  For $\omega$, $t$, $x$, $t'$ and $x'$ as in \eqref{omegacond11}, Lemma~\ref{uglobalbound} gives
\begin{align}
\nonumber&\hat Q_{T,2,\delta}(t,t',x')\\
\nonumber&\le C_{\ref{uglobalbound}}(\omega)c_1\int_{t'-(\delta\wedge(t'-t))}^{t'}\int p_{t'-s}(y-x')^2e^{2R_1K+(2(R_1+1))|y-x'|+2(R_1+1)(2K+1)}\\
\nonumber&\phantom{\le C_{\ref{uglobalbound}}(\omega)c_1\int_{t'-(\delta\wedge(t'-t))}^{t'} }\times[2^{-N}\vee|x-x'|+\sqrt{t'-s}+|y-x'|]^{2\gamma\xi}\\
\nonumber&\phantom{\le C_{\ref{uglobalbound}}(\omega)c_1\int_{t'-(\delta\wedge(t'-t))}^{t'} }\times\Bigl\{[\sqrt{a_n}\vee 2^{-N}\vee|x-x'|+\sqrt{t'-s}+|y-x'|]^{2\gamma(\tilde\gamma_m-1)}+a_n^{2\beta\gamma}\Bigr\}dyds\\
\nonumber&\le C_{\ref{uglobalbound}}(\omega)c_2(K)\int_{t'-(\delta\wedge(t'-t))}^{t'}(t'-s)^{-1/2}[(2^{-N}\vee|x-x'|)^{2\gamma\xi}+(t'-s)^{\gamma\xi}]\\
\label{hQTbnd1}&\phantom{\le C_{\ref{uglobalbound}}(\omega)c_2(K)\int}\times\Bigr\{((\sqrt{a_n}\vee 2^{-N}\vee|x-x'|)^{2\gamma(\tilde\gamma_m-1)}+(t'-s)^{\gamma(\tilde\gamma_m-1)})+a_n^{2\beta\gamma}\Bigr\}ds.
\end{align}
For $t\le s\le t'$ and $c,p\ge 0$, 
\[(c\vee |x-x'|)^p+(t'-s)^{p/2}\le 2(c\vee d)^p.\]
Use this with $c=2^{-N}$ or $\sqrt{a_n}\vee 2^{-N}$ to bound \eqref{hQTbnd1} by 
\begin{align*}
&C_{\ref{uglobalbound}}(\omega)c_3(K)((t'-t)\wedge\delta)^{{1\over 2}}[2^{-N}\vee d]^{2\gamma\xi}\Bigl\{(\sqrt{a_n}\vee 2^{-N}\vee d)^{2\gamma(\tilde\gamma_m-1)}+a_n^{2\beta\gamma}\Bigr\}\\
&\le C_{\ref{uglobalbound}}(\omega)c_3(K)(\sqrt{t'-t}\wedge\sqrt\delta)^{1-{\eta_1\over 2}}\bar d_N^{2\gamma\xi+{\eta_1\over 2}}\Bigl\{\bar d_{n,N}^{2\gamma(\tilde\gamma_m-1)}+a_n^{2\beta\gamma}\Bigr\}.
\end{align*}
The conditions on $\eta_1$ and definition of $\xi$ imply $2\gamma\xi+{\eta_1\over 2}\ge 2\gamma$, and so the bound on $\hat Q_{T,2,\delta}$ is established.
 
 Turning to $\hat Q_{T,1,\delta}$, we may assume $t'>\delta$, or else $\hat Q_{T,1,\delta}=0$.  Argue as in 
the derivation of \eqref{hQTbnd1} to see that for $\omega$, $t$, $t'$ and $x$ as in \eqref{omegacond11}, 
\begin{align}
\label{hQTbnd2}\hat Q_{T,1,\delta}(t,t',x)&\le C_{\ref{uglobalbound}}(\omega)c_4(K)\int_{(t-\delta)^+}^{t\wedge(t'-\delta)}(t-s)^{-1/2}[(2^{-2N\gamma\xi}+(t-s)^{\gamma\xi}]\\
\nonumber&\phantom{\le C_{\ref{uglobalbound}}(\omega)c_4(K)\int_{(t-\delta)^+}^{t\wedge(t'-\delta)}}\times\Bigr\{(\sqrt{a_n}\vee 2^{-N})^{2\gamma(\tilde\gamma_m-1)}+(t-s)^{\gamma(\tilde\gamma_m-1)}+a_n^{2\beta\gamma}\Bigr\}ds.
\end{align}
Elementary calculations give
\begin{equation}\label{integral1}
\int_{t-\delta}^{t\wedge(t'-\delta)}(t-s)^{-1/2}ds\le 2(\sqrt{t'-t}\wedge\sqrt\delta),
\end{equation}
and for $p\ge 0$,
\begin{equation}\label{integral2}
\int_{t-\delta}^{t\wedge(t'-\delta)}(t-s)^{p}ds\le \delta^p((t'-t)\wedge\delta).
\end{equation}
For the integral in \eqref{hQTbnd2} consider separately the cases (i) $\sqrt{t-s}<2^{-N}$, (ii) $\sqrt{t-s}\ge \sqrt{a_n}\vee 2^{-N}$, and (iii) $2^{-N}\le \sqrt{t-s}<\sqrt{a_n}\vee 2^{-N}$, the latter implying $\sqrt{a_n}\vee2^{-N}=\sqrt{a_n}$, to bound $\hat Q_{T,1,\delta}(t,t',x)$ by
\begin{align*}
&C_{\ref{uglobalbound}}(\omega)c_4(K)\Bigl\{\int_{(t-\delta)^+}^{t\wedge (t'-\delta)}(t-s)^{-1/2}ds2^{-2N\gamma\xi}\Bigl[(\sqrt{a_n}\vee 2^{-N})^{2\gamma(\tilde\gamma_m-1)}+a_n^{2\beta\gamma}\Bigr]\\
&\phantom{C_{\ref{uglobalbound}}(\omega)c_4(K)\Bigl\{ }+\int_{(t-\delta)^+}^{t\wedge (t'-\delta)} (t-s)^{\gamma(\tilde\gamma_m+\xi-1)-{1\over 2}}+a_n^{2\beta\gamma}(t-s)^{\gamma\xi-{1\over 2}}ds\\
&\phantom{C_{\ref{uglobalbound}}(\omega)c_4(K)\Bigl\{ }+\int_{(t-\delta)^+}^{t\wedge (t'-\delta)}(t-s)^{\gamma\xi-{1\over 2}}ds\Bigl[a_n^{\gamma(\tilde\gamma_m-1)}+a_n^{2\beta\gamma}\Bigr]\Bigr\}\\
&\le C_{\ref{uglobalbound}}(\omega)c_5(K)\Bigl\{(\sqrt{\delta}\wedge\sqrt{t'-t})\bar d_N^{-2\gamma\xi}\Bigl[\bar d_{n,N}^{2\gamma(\tilde\gamma_m-1)}+a_n^{2\beta\gamma}\Bigr]\\
&\phantom{\le C_{\ref{uglobalbound}}(\omega)c_4(K)\Bigl\{}+(\delta\wedge(t'-t))\Bigl[\delta^{\gamma(\tilde\gamma_m+\xi-1)-{1\over 2}}+a_n^{2\beta\gamma}\delta^{\gamma\xi-{1\over 2}}\Bigr]\\
&\phantom{\le C_{\ref{uglobalbound}}(\omega)c_4(K)\Bigl\{}+(\delta\wedge(t'-t))\Bigl[\delta^{\gamma\xi-{1\over 2}}a_n^{\gamma(\tilde\gamma_m-1)}+\delta^{\gamma\xi-{1\over 2}}a_n^{2\beta\gamma}\Bigr]\Bigr\}.
\end{align*}
In the last we have used \eqref{integral1} and \eqref{integral2}, and the fact that our choice of $\xi$ implies $\gamma\xi>1/2$ and hence our choices of $p$ are indeed non-negative when applying \eqref{integral2}.  Since $\delta\ge a_n$, the third term above is dominated by the second term.  Therefore
\begin{align*}
\hat Q_{T,1,\delta}(t,t',x)&\le C_{\ref{uglobalbound}}(\omega)c_6(K)\Bigl\{(\sqrt{\delta}\wedge\sqrt{t'-t})^{1-{\eta_1\over 2}}\bar d_N^{2\gamma\xi+{\eta_1\over 2}}\Bigl[\bar d_{n,N}^{2\gamma(\tilde\gamma_m-1)}+a_n^{2\beta\gamma}\Bigr]\\
&\phantom{\le C_{\ref{uglobalbound}}(\omega)c_6(K)}+(\delta\wedge(t'-t))^{1-{\eta_1\over 2}}\Bigl[\delta^{\gamma(\tilde\gamma_m+\xi-1)-{1\over 2}+{\eta_1\over 2}}+a_n^{2\beta\gamma}\delta^{\gamma\xi-{1\over 2}+{\eta_1\over 2}}\Bigr]\Bigr\}.
\end{align*}
Our choice of $\xi$ and conditions on $\eta_1$ imply that 
\[2\gamma\xi+{\eta_1\over 2}\ge 2\gamma\hbox{ and }\gamma(\xi-1)+{\eta_1\over 2}\ge 0,\]
and the required bound on $\hat Q _{T,1,\delta}$ follows.
\gdm

The above square function bounds suggest we will need a modified form of Lemma~\ref{dyadicexp} to  obtain our modulus of continuity for $u_{2,a_n^\alpha}$.  The proof of the following result is almost identical to that of Lemma~\ref{dyadicexp} and so is omitted.  

\begin{lemma}\label{dyadicexp2}
Let $c_0,c_1,c_2, k_0$ be positive (universal constants), $\eta\in(0,1/2)$, and 
$\Delta_i:\N\times(0,1]\to\Rp$, $i=1,2$ satisfy $\Delta_i(n,2^{-N+1})\le k_0\Delta_i(n,2^{-N})$ for all $n,N\in\N$ and $i=1,2$.  For $n\in\N$ and $\tau$ in a set $S$ assume $\{Y_{\tau,n}(t,x):(t,x)\in\Rp\times\Re\}$ is a real-valued continuous process.  Assume for each $(n,\tau)\in\N\times S$, $K\in\N$, and $\beta\in[0,1/2]$, there is an $N_0(\omega)=N_0(n,\eta,K,\tau, \beta)(\omega)\in\N$ a.s., stochastically bounded uniformly in $(n,\tau,\beta)$, such that for any $N\in\N$, $(t,x)\in \Rp\times\Re$, if $ d=d((t,x),(t',x'))\le 2^{-N}$, then
\begin{align}\label{dyadicexphyp2}
P(&|Y_{\tau,n}(t,x)-Y_{\tau,n}(t',x')|>d^{{1\over 2}(1-\eta)}\Delta_{1}(n,2^{-N})+{ d} ^{1-\eta}\Delta_2(n,2^{-N}), \\
\nonumber&\phantom{P(|Y_{\tau,n}(t,x)-Y_{\tau,n}(t',x')|>d^{{1\over 2}(1-\eta)}}(t,x)\in Z(N,n,K,\beta),N\ge N_0,t'\le T_K)\\
\nonumber&\le c_0\exp(-c_1{ d}^{-\eta c_2}).
\end{align}
Then there is an $N'_0(\omega)=N'_0(n,\eta,K,\tau,\beta)(\omega)\in\N$ a.s., also stochastically bounded uniformly in $(n,\tau,\beta)$, such that for all $N\ge N'_0(\omega)$, $(t,x)\in Z(N,n,K,\beta)(\omega)$, $d=d((t,x),(t',x'))\le 2^{-N}$, and $t'\le T_K$, 
\begin{equation*}
|Y_{\tau,n}(t,x)-Y_{\tau,n}(t',x')|\le 2^7k_0^2\Bigl[d^{{1\over 2}(1-\eta)}\Delta_1(n,2^{-N})+{ d}^{1-\eta}\Delta_2(n,2^{-N})\Bigr].
\end{equation*}
\end{lemma}

\paragraph{Proof of Proposition~\ref{u2modulus}.} The proof follows closely that of Proposition~\ref{Fmodulus}, using Lemma~\ref{dyadicexp2} in place of  Lemma~\ref{dyadicexp}. Let $R={25\over \eta_1}$ and choose $\eta_0\in({1\over R},{\eta_1\over 24})$.  Define $\bar d_N=d\vee 2^{-N}$, as usual, and set
\[\hat Q_{a_n^\alpha}(t,x,t',x')=\sum_{i=1}^2\hat Q_{S,i,a_n^\alpha,\eta_0}(t,x,t',x')+\hat Q_{T,1,a_n^\alpha}(t,t',x)+\hat Q_{T,2,a_n^\alpha}(t,t',x').\]
By Lemmas \ref{hQS1bound}, \ref{hQS2bound} and \ref{hQTbound}, for all $K\in \N$ (the restriction $K\ge K_1$ is illusory as these results only strengthen as $K$ increases) there is a constant $c_1(K,\eta_1)$ and $N_2(m,n,\eta_1,\ve_0,K,\beta)\in\N$ a.s., stochastically bounded uniformly in $(n,\beta)$, such that for all $N\in\N$, $(t,x)\in \Rp\times\Re$, 
\[\hbox{on }\{\omega:(t,x)\in Z(N,n,K+1,\beta), N\ge N_2\}, \]
\begin{align}
\nonumber&R_0^\gamma\hat Q_{a_n^\alpha}(t,x,t',x')^{1/2}\\
\nonumber&\le c_1(K,\eta_1)[a_n^{-\ve_0}+2^{2N_2}]\Bigl\{(d\wedge a_n^{{\alpha\over 2}})^{{1\over 2}(1-{\eta_1\over 2})}{\bar d_N}^\gamma\Bigl[(\bar d_N\vee a_n^{1/2})^{\gamma(\tilde\gamma_m-1)}+a_n^{\beta\gamma}\Bigl]\\
\nonumber&\phantom{\le c_1(K,\eta_1)[a_n^{-\ve_0}+2^{2N_2}]\Bigl\{}+(d\wedge a_n^{\alpha\over 2})^{1-{\eta_1\over 2}}\Bigl[ a_n^{{\alpha\over 2}(\gamma\tilde\gamma_m-{1\over 2})}+a_n^{\beta\gamma}a_n^{{\alpha\over 2}(\gamma-{1\over 2})}\Bigr]\Bigr\}\\
\label{hQbound2} &\phantom{\le c_1(K,\eta_1)[a_n^{-\ve_0}+2^{2N_2}]\Bigl\{(d\wedge a_n^{{\alpha\over 2}})^{{1\over 2}(1-{\eta_1\over 2})}{\bar d_N}^\gamma\Bigl[(\bar d_N\vee a_n^{1/2})^{\gamma(\tilde\gamma_m-1)}}\hbox{ for all }t\le t'\le T_K,\ |x'|\le K+2.
\end{align}
Let $N_3={25\over \eta_1}[N_2+N_4(K,\eta_1)]$, where $N_4(K,\eta_1)$ is 
chosen  large enough so that 
\begin{align}
\nonumber c_1(K,\eta_1)[a_n^{-\ve_0}+2^{2N_2}]2^{-N_3
\eta_1/8}&\le c_1(K,\eta_1)[a_n^{-\ve_0}+2^{2N_2}]2^{-6N_2}2^{-6N_4(K,\eta_1)}\\
\label{rancoeffbnd}&\le a_n^{-\ve_0}2^{-104}.
\end{align}
Let $\Delta_{i,u_2}=2^{-100}\Bar\Delta_{i,u_2}$, $i=1,2$.  Assume $d\le 2^{-N}$.  Use \eqref{rancoeffbnd} in \eqref{hQbound2} to see that for all $(t,x)$, $N$, on 
\[\{\omega:(t,x)\in Z(N,n,K+1,\beta),\ N\ge N_3\}\hbox{ (which implies $|x'|\le K+2$)},\]
\begin{align*}
R_0^\gamma\hat Q_{a_n^\alpha}(t,x,t',x')^{1/2}
&\le (d\wedge a_n^{{\alpha\over 2}})^{{1\over 2}
(1-{3\eta_1\over 4})}\Delta_{1,u_2}(m,n,2^{-N})/16\\
&\quad +(d\wedge a_n^{{\alpha\over 2}})^{1-
{5\eta_1\over 8}}\Delta_{2,u_2}(m,n)/16\quad\hbox{ for all }t\le t'\le T_K,\ x'\in\Re.
\end{align*}
Combine this with \eqref{u2decompa}, \eqref{Dbound}, the definition of $\hat Q_{a_n^\alpha}$ and the Dubins-Schwarz theorem to conclude that for $t\le t'$, $x'\in\Re$, and $d((t,x),(t',x'))\le 2^{-N}$, 
\begin{align}
\nonumber &P\Bigl(|u_{2,a_n^\alpha}(t,x)-u_{2,a_n^\alpha}(t',x')|\ge d^{{1\over 2}(1-\eta_1)}\Delta_{1,u_2}(m,n,2^{-N})/4+d^{1-\eta_1}\Delta_{2,u_2}(m,n)/4,\\
\nonumber&\phantom{P(|u_{2,a_n^\alpha}(t,x)-u_{2,a_n^\alpha}(t',x')|\ge d^{{1\over 2}(1-\eta_1)}\Delta_{1,u_2}}(t,x)\in Z(N,n,K+1,\beta),\ N\ge N_3,\ t'\le T_K\Bigr)\\
\nonumber&\le 3P\Bigl(\sup\{|B(u)|:u\le [d^{{1\over 2}(1-{3\eta_1\over 4})}\Delta_{1,u_2}(m,n,2^{-N})/16+d^{1-{5\eta_1\over 8}}\Delta_{2,u_2}(m,n)/16]^2\}\\
\nonumber&\phantom{\le 3P\Bigl(\sup\{|B(u)|:u\le [d^{1-{3\eta_1\over 4}}}\ge (d^{{1\over 2}(1-\eta_1)}\Delta_{1,u_2}(m,n,2^{-N})+d^{1-\eta_1}\Delta_{2,u_2}(m,n))/12\Bigr)\\
\label{Gaussu2a}&\le 3P\Bigl(\sup_{u\le 1}|B(u)|\ge d^{-{\eta_1\over 8}}\Bigr)\le c_0\exp\Bigl(-{d^{-{\eta_1\over 4}}\over 2}\Bigr).
\end{align}
Here $B(u)$ is a one-dimensional Brownian motion.  

If $(t,x)\in Z(N,n,K,\beta)$, $d\le 2^{-N}$, and $t'\le t$, then as in the proof of Proposition~\ref{Fmodulus}, $(t',x')\in Z(N-1,n,K+1,\beta)$ and one can interchange the roles of $(t,x)$ and $(t',x')$ and replace $N$ with $N-1$ in the above to conclude (as in \eqref{gaussbound3}),
\begin{align}
\nonumber P&(|u_{2,a_n^\alpha}(t,x)-u_{2,a_n^\alpha}(t',x')|\ge d^{{1\over 2}(1-\eta_1)}\Delta_{1,u_2}(m,n,2^{-N})+d^{1-\eta_1}\Delta_{2,u_2}(m,n),\\ 
\nonumber&\phantom{P(|u_{2,a_n^\alpha}(t,x)-u_{2,a_n^\alpha}(t',x')|\ge d^{{1\over 2}(1-\eta_1)}\Delta_{1,u_2}}(t,x)\in Z(N,n,K,\beta),\ N\ge N_3+1)\\
\label{Gaussu2b}&\le c_0\exp\Bigl(-{d^{-{\eta_1\over 4}}\over 2}\Bigr).
\end{align}
\eqref{Gaussu2a} and \eqref{Gaussu2b} allow us to apply Lemma~\ref{dyadicexp2} with $\tau=\alpha$, $Y_{\tau,n}=u_{2,a_n^\alpha}$, and $k_0=4$. The result is then immediate once one recalls that $2^7k_0^2\Delta_{i,u_2}=2^{-89}\bar\Delta_{i,u_2}$. 
\gdm

\section{Incorporating Drifts}\label{secdrifts}
\setcounter{equation}{0}
\setcounter{theorem}{0} 

Beginning in 
Section~\ref{sec:verifiction_Hyp}
we assumed that the drift $b$ is zero.  Here we point out what additional reasoning is needed to include a drift $b$ satisfying \eqref{bLip}. If $B(s,y)=b(s,y,X^1(s,y))-b(s,y,X^2(s,y))$, then \eqref{ueq} becomes
\begin{align}\label {bueq}
u(t,x)&=\int_0^t\int p_{t-s}(y-x)D(s,y)W(ds,dy)+\int_0^t\int p_{t-s}(y-x)B(s,y)dyds\hbox { a.s. for all }(t,x)\\
\nn&\equiv u_D(t,x)+u_B(t,x),
\end{align}
and by \eqref{bLip},
\begin{equation}\label{bubnd}
|B(s,y)|\le B|u(s,y)|.
\end{equation}
If $u_{1,\delta},u_{2,\delta}$ and $G_\delta$ are defined as at the beginning of 
Section~\ref{sec:verifiction_Hyp}, then as for \eqref{u1eq}, \eqref{u2eq} and Lemma~\ref{G'}, but now using the ordinary Fubini theorem for $u_B$, we get
\begin{align*}
u_{1,\delta}(t,x)&=\int_0^{(t-\delta)^+}\int p_{t-s}(y-x)D(s,y)W(ds,dy)+\int_0^{(t-\delta)^+}\int p_{t-s}(y-x)B(s,y)dyds\\
&\equiv u_{1,D,\delta}(t,x)+u_{1,B,\delta}(t,x),\\
u_{2,\delta}(t,x))&=\int_{(t-\delta)^+}^t\int p_{t-s}(y-x)D(s,y)W(ds,dy)+\int_{(t-\delta)^+}^t\int p_{t-s}(y-x)B(s,y)dyds\\
&\equiv u_{2,D,\delta}(t,x)+u_{2,B,\delta}(t,x),\\
\end{align*}
and
\begin {align*}
-G'_\delta(s,t,x)&\equiv F_\delta(s,t,x)\\
&=
\int_0^{(s-\delta)^+}\int p'_{(t\vee s)-r}(y-x)D(r,y)W(dr,dy)+\int_0^{(s-\delta)^+}\int p'_{(t\vee s)-r}(y-x)B(r,y)dydr\\
&\equiv F_{D,\delta}(s,t,x)+F_{B,\delta}(s.t,x).
\end{align*}
In addition, no changes are required in the verification of $(P_0)$ (including the refinement noted in Remark~\ref{m0case}) or the proof of Lemma~\ref{uglobalbound}.

The theorems in Section~\ref{sec4} apply directly to quantities like $u_{i,D,\delta}$ and $F_{D,\delta}$.  The corresponding expressions $u_{i,B,\delta}$ and $F_{B,\delta}$ are in fact much easier to handle because we are dealing with a deterministic integral and so regularity properties 
are easy to  
read off directly from the bounds in Lemma~\ref{uglobalbound}.  Furthermore, the Lipschitz condition on $b$ effectively sets $\gamma=1$ for these calculations. To illustrate this, we now prove a simple result which includes both Propositions~\ref{Fmodulus} and \ref{Fmodulus2} for $F_{B,a_n^\alpha}$ and only requires $(P_0)$, a consequence of the ``crude" modulus Theorem~\ref{theorem:lipmod}, already noted above.  

\begin{proposition}\label{bprop}
For any $\eta_1\in(0,{1\over 2})$ and $K\in\N^{\ge K_1}$ there is an $N_{\ref{bprop}}(\eta_1,K)(\omega)\in\N$ a.s. such that for all $n\in\N$, $\alpha\in[0,1]$ and $\beta\in[0,{1\over 2}]$, $N\ge N_{\ref{bprop}}(\eta_1,K)$, $(t,x)\in Z(N,n,K,\beta)$, and $t'\le T_K$,
\begin{align*}&d=d((t,x),(t',x'))\le 2^{-N}\hbox{ and }|s'-s|\le N^{-1}\hbox{ imply }\\
&|F_{B,a_n^\alpha}(s,t,x)-F_{B,a_n^\alpha}(s',t',x')|\le 2^{-88}\Bigl[|s'-s|^{{1\over 2}}2^{-N(1-\eta_1)}+|s'-s|^{1-{\eta_1\over 2}}\\
&\phantom{F_{B,a_n^\alpha}(s,t,x)-F_{B,a_n^\alpha}(s',t',x')|\le 2^{-88}\Bigl[}+d^{1-\eta_1}(1+a_n^{-3\alpha/4}2^{-2N})\Bigr].
\end{align*}
\end{proposition}
\paragraph{Proof.} Let $\xi=1-{\eta_1\over 4}$ and assume first that 
\begin{equation} N\ge N_1(0,\xi,K+1)+1,\label{N1cond}
\end{equation}
where $N_1$ is as in $(P_0)$ and Remark~\ref{m0case}. Assume
\begin{equation}\label{bcond1}
(t,x)\in Z(N,n,K,\beta), t'\le T_K, d\le 2^{-N}\hbox{ and }|s'-s|\le N^{-1}.
\end{equation}
One easily checks that $|t\wedge s-t'\wedge s'|\le N^{-1}$ and so, by replacing $(s,s')$ with $(t\wedge s,t'\wedge s')$, we may assume that $s\le t$ and $s'\le t'$.  Define $\overline s=s\vee s'$ and $\underline s=s\wedge s'$.
As before, $(t',x')\in Z(N-1,n,K+1,\beta)$, and again, by interchanging $(t,x)$ with $(t',x')$, we may assume $t\le t'$ (this is the reason for the $K+1$ and adding $1$ to $N_1$ in \eqref{N1cond}).  As for \eqref{Fdecomp}, \eqref{bubnd} implies
\begin{align*}
&|F_{B,a_n^\alpha}(s,t,x)-F_{B,a_n^\alpha}(s',t',x')|\\
&\le\int_{(\underline s-a_n^\alpha)^+}^{(\overline s-a_n^\alpha)^+}\int |p'_{t'-r}(y-x')|B|u(r,y)|dydr+\int _0^{(s-a_n^\alpha)^+}\int |p'_{t'-r}(y-x')-p'_{t-r}(y-x)|B|u(r,y)|dydr\\
&\equiv T_1+T_2.
\end{align*}

To bound $T_1$, we may assume $\overline s\ge a_n^\alpha$.  Elementary inequalities using $t'\le T_K\le K$, show that for $p\ge0$,
\begin{align}\label{belem}
\int|p'_{t'-r}(y-x')||y-x'|^pe^{|y-x'|}dy&\le (t'-r)^{-1}\int |y-x'|^{p+1}e^{|y-x'|}p_{t'-r}(y-x')dy\\
\nn&\le c_1(K,p)(t'-r)^{{p-1\over 2}}.
\end{align}
Now apply Lemma~\ref{uglobalbound} with $m=0$ to see that 
\begin{align}
\nn T_1&\le B\sqrt{C_{\ref{uglobalbound}}}\int_{(\underline s-a_n^\alpha)^+}^{\overline s-a_n^\alpha}\int |p'_{t'-r}(y-x')|e^{|y-x'|}e^{|x-x'|}((\sqrt{t-r}+|y-x|)\vee 2^{-N})^\xi dydr\\
\nn&\le c_2(K)\sqrt{C_{\ref{uglobalbound}}}\int_{(\underline s-a_n^\alpha)^+}^{\overline s-a_n^\alpha}\int |p'_{t'-r}(y-x')|e^{|y-x'|}[(|x-x'|\vee 2^{-N})^\xi+\sqrt{t'-r}^\xi+|y-x'|^\xi]dydr\\
\nn&\le c_3(K)\sqrt{C_{\ref{uglobalbound}}}\int_{(\underline s-a_n^\alpha)^+}^{\overline s-a_n^\alpha}(t'-r)^{-1/2}[2^{-N\xi}+(t'-r)^{\xi/2}]dr\quad\hbox{ (by \eqref{belem} and }|x-x'|\le 2^{-N})\\
\nn&\le c_4(K)
2^{2N_1(\omega)}[|s'-s|^{1/2}2^{-N(1-{\eta_1\over 4})}+|s'-s|^{1-{\eta_1\over 8}}]\\
\nn&\le c_4(K)
2^{2N_1(\omega)}[2^{-N\eta_1/2}+|s'-s|^{\eta_1/4}][|s'-s|^{1/2}2^{-N(1-\eta_1)}+|s-s'|^{1-{\eta_1\over 2}}]\\
\label{bT1bound}&\le c_4(K)2^{2N_1(\omega)}[2^{-N\eta_1/2}+N^{-\eta_1/4}][|s'-s|^{1/2}2^{-N(1-\eta_1)}+|s-s'|^{1-{\eta_1\over 2}}].
\end{align}
Recall in the above that we may set $\ve_0=0$ in the definition of $C_{\ref{uglobalbound}}$ (see Remark~\ref{m0case}).  

For $T_2$ we use both $|u(r,y)|\le K e^{|y|}$ for $r\le t'\le T_K$ ($K\ge K_1$) and Lemma~\ref{uglobalbound} with $m=0$ to write (we may assume $s>a_n^\alpha$),
\begin{align*}
T_2\le &\int_0^{s-a_n^\alpha}\int |p'_{t-r}(y-x')-p'_{t-r}(y-x)|BKe^{|y|}1(|y-x|>(t'-r)^{{1\over 2}-{\eta_1\over 4}}\vee (2|x'-x|))dydr\\
&+\int_0^{s-a_n^\alpha}\int |p'_{t-r}(y-x')-p'_{t-r}(y-x)|B\sqrt{C_{\ref{uglobalbound}}}e^{|y-x|}[(\sqrt{t-r}+|y-x|)\vee 2^{-N}]^\xi\\
&\phantom{+\int_0^{s-a_n^\alpha}\int |p'{t--r}(y-x')}\times1(|y-x|\le (t'-r)^{{1\over 2}-{\eta_1\over 4}}\vee (2|x'-x|))dydr\\
\le& c_5(K)\sqrt{C_{\ref{uglobalbound}}}\Bigl[\int_0^{s-a_n^\alpha}\int |p'_{t-r}(y-x')-p'_{t-r}(y-x)|e^{2|y-x|}\\
&\phantom{c_5(K)\sqrt{C_{\ref{uglobalbound}}}\Bigl[\int_0^{s-a_n^\alpha}\int |p'_{t-r}(y-x')}\times1(|y-x|>(t'-r)^{{1\over 2}-{\eta_1\over 4}}\vee (2|x'-x|))e^{-|y-x|}dydr\\
&+\int _0^{s-a_n^\alpha}\int  |p'_{t-r}(y-x')-p'_{t-r}(y-x)|1(|y-x|\le 2(2K+1))dy((t'-r)^{{1\over 2}-{\eta_1\over 4}}\vee |x'-x|\vee 2^{-N})^\xi dr\Bigr]\\
\le &c_6(K)\sqrt{C_{\ref{uglobalbound}}}\Bigl[\int_0^{s-a_n^\alpha}\int (p'_{t-r}(y-x')-p'_{t-r}(y-x))^2e^{4|y-x|}\\
&\phantom{c_6(K)\sqrt{C_{\ref{uglobalbound}}}\Bigl[\int_0^{s-a_n^\alpha}}\times1(|y-x|>(t-r)^{{1\over 2}-{\eta_1\over 4}}\vee (2|x'-x|))dy^{1/2}\Bigl(\int e^{-2|y-x|}dy\Bigr)^{1/2}dr\\
&+\int_0^{s-a_n^\alpha}\int (p'_{t-r}(y-x')-p'_{t-r}(y-x))^2dy^{1/2}(2(2K+1))^{1/2}[(t-r)^{1/2}\vee 2^{-N})^{\xi(1-{\eta_1\over 2})}dr\Bigr],
\end{align*}
where in the last line for the second term we use $d\le 2^{-N}$ and 
\[\sqrt{t'-r}^{1-{\eta_1\over 2}}\le \sqrt{t'-t}^{1-{\eta_1\over 2}}+\sqrt{t-r}^{1-{\eta_1\over 2}}\le d^{1-{\eta_1\over 2}}+\sqrt{t-r}^{1-{\eta_1\over 2}}.\]
Now apply Lemma~\ref{pt'bnds} and conclude
\begin{align}
\nn T_2\le & c_7(K,\eta_1)\sqrt{C_{\ref{uglobalbound}}}\Bigl\{\int_0^{s-a_n^\alpha}(t-r)^{-3/4}
\exp\Bigl\{{-\eta_1(t-r)^{-\eta_1/2}\over 128}\Bigr\}\Bigl(1\wedge {d^2\over t-r}\Bigr)^{{1\over 2}-{\eta_1\over 4}}dr\\
\nn&\phantom{c_7(K,\eta_1)\sqrt{C_{\ref{uglobalbound}}}\Bigl\{}+\int_0^{s-a_n^\alpha}(t-r)^{-3/4}\Bigl(1\wedge{d^2\over t-r}\Bigr)^{1/2}[(t-r)^{1/2}\vee 2^{-N}]^{\xi(1-{\eta_1\over 2})}dr\Bigr\}\\
\nn\le & c_7(K,\eta_1)\sqrt{C_{\ref{uglobalbound}}}\Bigl\{\int_0^{s-a_n^\alpha}(t-r)^{-3/4}
\exp\Bigl\{{-\eta_1(t-r)^{-\eta_1/2}\over 128}\Bigr\}\Bigl(1\wedge {d^2\over t-r}\Bigr)^{{1\over 2}-{\eta_1\over 4}}dr\\
\nn&\phantom{c_7(K,\eta_1)\sqrt{C_{\ref{uglobalbound}}}\Bigl\{}+\int_0^{s-a_n^\alpha}(t-r)^{-{3\over 4}+{\xi\over 2}(1-{\eta_1\over 2})}\Bigl(1\wedge {d^2\over t-r}\Bigr)^{1/2}dr\\
\nn&\phantom{c_7(K,\eta_1)\sqrt{C_{\ref{uglobalbound}}}\Bigl\{}+\int_0^{s-a_n^\alpha}1(r\ge t-2^{-2N})(t-r)^{-{3\over 4}}\Bigl(1\wedge {d^2\over t-r}\Bigr)^{1/2}dr2^{-N\xi(1-{\eta_1\over 2})}\Bigr\}\\
\nn\le &c_8(K,\eta_1)\sqrt{C_{\ref{uglobalbound}}}\Bigl\{\int _0^{s-a_n^\alpha}\Bigl(1\wedge {d^2\over t-r}\Bigr)^{{1\over 2}-{\eta_1\over 4}}dr+d\\
\nn&\phantom{le c_8(K,\eta_1)\sqrt{C_{\ref{uglobalbound}}}\Bigl\{}+1(a_n^\alpha<2^{-2N})d(d^2\vee a_n^\alpha)^{-1/4}2^{-N\xi(1-{\eta_1\over 2})}\Bigr\}\\
\nn\le &c_9(K,\eta_1)\sqrt{C_{\ref{uglobalbound}}}\Bigl\{d^{1-{\eta_1\over 2}}+1(a_n^\alpha<2^{-2N})da_n^{-3\alpha/4}a_n^{\alpha/2}2^{-N\xi(1-{\eta_1\over 2})}\Bigr\}\\
\nn\le &c_9(K,\eta_1)\sqrt{C_{\ref{uglobalbound}}}d^{1-{3\eta_1\over 4}}\Bigl[1+a_n^{-3\alpha\over 4}2^{-{3N\eta_1\over 4}}2^{-N}2^{-N(1-{3\eta_1\over 4})}\Bigr]\\
\label{bT2bound}\le &c_{10}(K,\eta_1)2^{2N_1(\omega)}2^{-N\eta_1/4}d^{1-\eta_1}[1+a_n^{-{3\alpha\over 4}}2^{-2N}].
\end{align}
We have used Lemma~\ref{Jbnd} in the above with a bit of algebra to see which case applies, and in the last two lines again used $d\le 2^{-N}$.   \eqref{bT1bound} and \eqref{bT2bound} together show there is an $N_{\ref{bprop}}(\eta_1,K)(\omega)\in\N$ a.s. such that $N_{\ref{bprop}}(\eta_1,K)\ge N_1(0,\xi,K+1)+1$ and, if $N\ge N_{\ref{bprop}}$, then 
\begin{align*}
T_1+T_2&\le 2^{-88}\Bigl[|s'-s|^{1/2}2^{-N(1-\eta_1)}+|s'-s|^{1-{\eta_1\over 2}}+d^{1-
\eta_1}[1+a_n^{-{3\alpha\over 4}}2^{-2N}]\Bigr].
\end{align*}
\gdm

A bit of arithmetic shows that the above bound in the contexts of Propositions~\ref{Fmodulus} and \ref{Fmodulus2} lead to upper bounds that are bounded by the ones obtained there for increments of $F_{D,a_n^\alpha}$.  For Proposition~\ref{Fmodulus} one only needs the first two terms in $\bar{\Delta}_{u'_1}$ and we leave this easy check for the reader.  For Proposition~\ref{Fmodulus2} we may set $(s',t',x')=(t,t,x)$ in the above so that the upper bound becomes
\begin{align*}
&2^{-88}[|t-s|^{1/2}2^{-N(1-\eta_1)}+|t-s|^{1-{\eta_1\over 2}}]\\
&\le 2^{-88}[2^{-N(1-\eta_1)}+(t-s)^{1-\eta_1\over 2}(\sqrt{t-s}\vee \sqrt{a_n})^{\gamma\tilde \gamma_m-{3\over 2}}],
\end{align*}
which is bounded by two of the terms on the right-hand side of the upper bound in Proposition~\ref{Fmodulus2}.  Hence we may combine these bounds for $F_{B,a_n^\alpha}$ with those derived in Propositions \ref{Fmodulus} and \ref{Fmodulus2}  for $F_{D,a_n^\alpha}$ and hence complete the proofs of Propositions \ref{Fmodulus} and \ref{Fmodulus2} (and hence also Corollary~\ref{u1'modulus}) for solutions with Lipschitz drifts. 
 
We omit the analogues of the above for Propositions~\ref{Gmodulus1} and \ref{u2modulus} as they are even simpler.  
The proof of Propositions~\ref{u1modulus} and \ref{Pminduction} now proceed as before.  With Propositions ~\ref{Pminduction}, \ref{u1modulus}, \ref{u2modulus}, and 
\ref{Fmodulus2},
 and Corollary~\ref{u1'modulus} in hand, the proof of 
Proposition~\ref{tildeJ}
may now be completed for Lipschitz drifts $b$, exactly as in Section~\ref{sec5}. 
Then verification of the hypotheses of Proposition~\ref{prop:Inbound} may now be completed for Lipschitz drifts $b$ exactly as in Section~\ref{sec:verifiction_Hyp} and this finishes the proof of Theorem~\ref{theorem:unique}.

\end{document}